\documentclass[a4paper,11pt]{article}
\usepackage{hyperref,appendix,amssymb,amsmath,amsthm,enumerate, graphicx,pstricks,fullpage}

\newtheorem{thm}{Theorem}[section]
\newtheorem{lem}[thm]{Lemma}

\newtheorem{cor}[thm]{Corollary}

\newtheorem{propn}[thm]{Proposition}

\newtheorem{assu}{Assumption}

\begin{document}

\title{Scaling limit for the random walk on the largest connected component of the critical random graph}
\author{David A. Croydon\footnote{D. A. Croydon: Dept of Statistics,
University of Warwick, Coventry, CV4 7AL, UK;
{d.a.croydon@warwick.ac.uk}}}
\date{5 July 2011}

\maketitle

\begin{abstract}
In this article, a scaling limit for the simple random walk on the largest connected component of the Erd\H{o}s-R\'{e}nyi random graph $G(n,p)$ in the critical window, $p=n^{-1}+\lambda n^{-4/3}$, is deduced. The limiting diffusion is constructed using resistance form techniques, and is shown to satisfy the same quenched short-time heat kernel asymptotics as the Brownian motion on the continuum random tree.
\end{abstract}

\section{Introduction}

\begin{figure}[t]
\begin{center}
\vspace{-38pt}
\scalebox{0.6}{\includegraphics{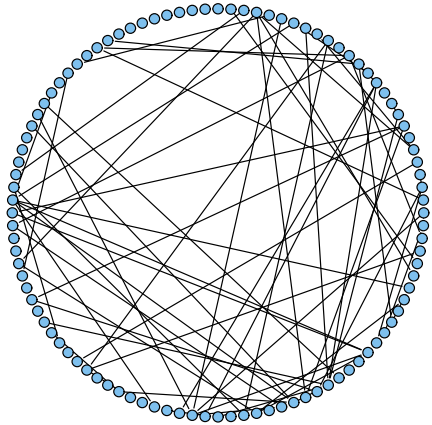}}
\vspace{-55pt}
\end{center}
\caption{A typical realisation of $G(n,p)$ with $n=100$ and $p=0.01$.}\label{erfig}
\vspace{-10pt}
\end{figure}

It is known that the asymptotic behaviour of the Erd\H{o}s-R\'{e}nyi random graph $G(n,p)$, in which every edge of the complete graph on $n$ labelled vertices $\{1,\dots,n\}$ is present with probability $p$ independently of the other edges (see Figure \ref{erfig}), exhibits a phase transition at $p\sim n^{-1}$. On the one hand, for $p\sim cn^{-1}$ with $c>1$, the largest connected component $\mathcal{C}_1^n$ of $G(n,p)$ incorporates a non-trivial proportion of the $n$ original vertices asymptotically. By contrast, if $c<1$, then $\mathcal{C}_1^n$ consists of only $O(\ln n)$ vertices. A third kind of behaviour is seen at criticality itself; for when $p= n^{-1}$, the number of vertices of $\mathcal{C}_1^n$ is of order $n^{2/3}$ (all these results can be found in \cite{ER}). Under a finer scaling $p=n^{-1}+\lambda n^{-4/3}$, where $\lambda\in\mathbb{R}$ is fixed -- the so-called critical window, it is also possible to describe the asymptotic structure of $\mathcal{C}_1^n$. Specifically, in this regime, when graph distances are rescaled by $n^{-1/3}$, the largest connected component $\mathcal{C}_1^n$ converges to a random fractal metric space, $\mathcal{M}$ say, whose distribution depends on the particular value of $\lambda$ chosen (\cite{ABG}). In this article, our goal is to add a further level of detail to this picture: we will consider the discrete time simple random walk on $\mathcal{C}_1^n$ in the critical window, and show that it converges, when rescaled appropriately, to a diffusion on $\mathcal{M}$.

The precise characterisation of $\mathcal{M}$ in \cite{ABG} involves ``glued'' real trees, and this picture turns out to be especially useful in constructing the limiting diffusion $X^\mathcal{M}$, which we will call the Brownian motion on $\mathcal{M}$. In particular, in defining $\mathcal{M}$, the authors of \cite{ABG} start by introducing a random compact real tree $\mathcal{T}$. On the same probability space, a point process is described, and this gives a recipe for selecting a finite number of pairs $\{u_i,v_i\}$, $i=1,\dots ,J$, of vertices of $\mathcal{T}$. (Note that $J$ is a random variable.) Once the joint distribution of $\mathcal{T}$ and the point process is chosen appropriately, what is shown in \cite{ABG} is that the scaling limit of the largest connected component $\mathcal{C}_1^n$ in the critical window is simply the metric space $\mathcal{M}$ we arrive at from $\mathcal{T}$ by identifying $u_i$ and $v_i$ for each $i=1,\dots,J$. To build a process on $\mathcal{M}$ in this setting, we first note that results of \cite{Kigamidendrite} let us view $\mathcal{T}$ as an electrical network by equipping it with a corresponding resistance form $(\mathcal{E}_\mathcal{T},\mathcal{F}_\mathcal{T})$. From this, we obtain a related resistance form $(\mathcal{E}_\mathcal{M},\mathcal{F}_\mathcal{M})$ on $\mathcal{M}$ by ``fusing'' the vertices $u_i$ and $v_i$ together for $i=1,\dots, J$. Finally, letting $\mu_\mathcal{M}$ be the natural measure on $\mathcal{M}$ -- the scaling limit of the uniform measure on $\mathcal{C}_1^n$, we can further prove that $(\mathcal{E}_\mathcal{M},\mathcal{F}_\mathcal{M})$ is a local, regular Dirichlet form on $L^2(\mathcal{M},\mu_\mathcal{M})$.  It is by applying the standard association between such quadratic forms and Markov processes that the Brownian motion $X^\mathcal{M}$ will be defined.

In proving the main convergence result for the discrete time simple random walk $X^{\mathcal{C}_1^n}$ on $\mathcal{C}_1^n$, an argument developed for demonstrating convergence to the Brownian motion on the continuum random tree is adapted. In \cite{Croydoncbp}, it was established that for any family of graph trees $(T_n)_{n\geq 1}$ that converge in a suitable manner to the continuum random tree of Aldous (see \cite{Aldous2}, for example), the associated simple random walks converge to a Markov process called the Brownian motion on the continuum random tree. The key to proving this result was first considering approximations to the simple random walks and the limiting diffusion on subtrees spanning a fixed number of leaves, where proving convergence from the discrete to continuous models was straightforward, and then letting the number of leaves go to infinity. By constructing the sets $\mathcal{C}_1^n$ and $\mathcal{M}$ from trees as in \cite{ABG}, we are able to define similarly useful approximations for $X^{\mathcal{C}_1^n}$ and $X^\mathcal{M}$ that take values in subsets spanning a fixed number of vertices of $\mathcal{C}_1^n$ and $\mathcal{M}$, respectively. Once these subsets and approximating processes have been chosen appropriately, an almost identical argument to the one applied for trees in \cite{Croydoncbp} yields that
\begin{equation}\label{rwcresult}
\left(n^{-1/3}X^{\mathcal{C}_1^n}_{\lfloor tn\rfloor}\right)_{t\geq 0} \rightarrow \left(X^{\mathcal{M}}_{t}\right)_{t\geq 0}
\end{equation}
in distribution, where we postpone the precise statement of this result until Section \ref{rwc}. Of course, being an example of a random walk in a random environment-type problem, there are two kinds of results such as (\ref{rwcresult}) that we could prove. Firstly, a quenched result, where we fix a sequence of typical realisations of $\mathcal{C}_1^n$ that converge to a typical realisation of $\mathcal{M}$, and secondly, an annealed result, where the distributions of the processes $X^{\mathcal{C}_1^n}$ and $X^\mathcal{M}$ are averaged over the randomness of the environments $\mathcal{C}_1^n$ and $\mathcal{M}$. We will prove both; see Theorems \ref{quench} and \ref{anne}, respectively, and we appeal to these results to justify our referring to $X^\mathcal{M}$ as Brownian motion on $\mathcal{M}$. Note that the time scaling factor $n$ for the simple random walk on the largest connected component $\mathcal{C}_1^n$ in the critical window has previously been observed in the mixing time asymptotics of $X^{\mathcal{C}_1^n}$, see \cite{NachPer}, Theorem 1.1.

Once the Brownian motion $X^\mathcal{M}$ has been demonstrated to be the scaling limit of $X^{\mathcal{C}_1^n}$, it is natural to investigate further properties of the process. To make a first step in this direction we again turn to the link between $\mathcal{M}$ and $\mathcal{T}$, which immediately allows us to transfer some known results about the Brownian motion on the continuum random tree to our setting. Specifically, modulo random mass and distance scaling factors, the distribution of the random compact real tree $\mathcal{T}$ is absolutely continuous with respect to the distribution of the continuum random tree. Thus, when suitably rescaled by constants, a typical realisation of $\mathcal{T}$ looks exactly like a typical realisation of the continuum random tree. Moreover, since $\mu_\mathcal{M}$ is a non-atomic measure and there are only a finite number of pairs $\{u_i,v_i\}$, $i=1,\dots,J$, it is clear that $\mu_\mathcal{M}$-a.e. every point of $\mathcal{M}$ admits a neighbourhood that is isometric to a neighbourhood in $\mathcal{T}$. With the local geometry of $\mathcal{M}$ and $\mathcal{T}$ being the same and $\mathcal{T}$ looking like the continuum random tree, it is easy to check that the short-time behaviour of $X^\mathcal{M}$ is the same as the short-time behaviour of the Brownian motion on the continuum random tree. For example, if we let $(p_t^\mathcal{M}(\bar{x},\bar{y}))_{\bar{x},\bar{y}\in\mathcal{M},t>0}$ be the transition density of $X^\mathcal{M}$, then it is possible to show that every typical realisation of $\mathcal{M}$ satisfies
\begin{equation}\label{specdim}
\lim_{t\rightarrow 0}\frac{2\ln p^\mathcal{M}_t(\bar{x},\bar{x})}{-\ln t}= \frac{4}{3},\hspace{20pt}\forall \bar{x}\in\mathcal{M}.
\end{equation}
Thus the spectral dimension of (the Brownian motion on) $\mathcal{M}$ is almost-surely $4/3$, which is identical to that of the continuum random tree, see \cite{Croydoncrt}. More detailed short-time transition density asymptotics are discussed in Section \ref{propsec}.

To conclude the introduction, let us draw attention to the connections between our results and those thought to hold for critical edge percolation clusters on the lattice $\mathbb{Z}^d$. In high dimensions, $d>6$, large critical edge percolation clusters on $\mathbb{Z}^d$, when rescaled, are thought to look like a random structure known as the integrated super-Brownian excursion (see \cite{Slade}, Section 16.5, for example). The integrated super-Brownian excursion can simply be thought of as a random embedding of the continuum random tree into $\mathbb{R}^d$, and one might therefore hypothesise that by mapping the Brownian motion on the continuum random tree into $\mathbb{R}^d$ with the same embedding we might be able to define the scaling limit of the random walks on the critical percolation clusters. Although  in the percolation setting a result along these lines has not yet been proved, this procedure does yield the scaling limit of the random walks on the trees generated by related branching random walk models, at least in the case when $d>7$ (see \cite{CroydonH}). Of course, the complete graph on $n$ points does not look very much like the integer lattice. However, the relationship between $\mathcal{M}$ and the continuum random tree discussed above suggests that, at criticality, the asymptotic structure of the largest connected component of the Erd\H{o}s-R\'{e}nyi random graph, which can be thought of as an edge percolation model, closely corresponds to the asymptotic structure of the large critical clusters in the $\mathbb{Z}^d$ edge percolation setting, at least locally. What we conjecture on the basis of the results of this article is that in the critical regime the dynamical properties are also similar. In fact, work already exists that indicates that the spectral dimension of the scaling limit of the simple random walks on large critical percolation clusters in $\mathbb{Z}^d$ in high dimensions agrees with the $4/3$ seen at (\ref{specdim}), see \cite{NK}.

The remainder of the article is arranged as follows. In Section \ref{sectree} we present details of the construction of $X^\mathcal{M}$ on typical realisations of $\mathcal{M}$. We also introduce there some approximations for $X^\mathcal{M}$ on finite length subsets of $\mathcal{M}$, and in Section \ref{seccont} describe continuity properties of the laws of these. The distribution of $\mathcal{M}$ and largest connected component scaling limit result of \cite{ABG} are recalled more thoroughly than above in Section \ref{crg}. In order to formalise the main simple random walk convergence result and check the measurability of our construction of $X^\mathcal{M}$, we need to introduce a space of paths on compact length spaces, and this is done in Section \ref{contsec}. In Section \ref{subsetsec} we characterise what constitutes the quenched convergence of $\mathcal{C}_1^n$ to $\mathcal{M}$ in a way that will be most helpful to us, and present some other preliminary results about subsets of $\mathcal{C}_1^n$. Finally, we establish our precise versions of (\ref{rwcresult}) in Section \ref{rwc} and conclude the article by outlining some properties of $X^\mathcal{M}$ in Section \ref{propsec}.

\section{Brownian motion on fused real trees}\label{sectree}

The goal of this section is to construct, at least for typical realisations of $\mathcal{M}$, the Dirichlet form $(\mathcal{E}_\mathcal{M},\mathcal{F}_\mathcal{M})$ on $L^2(\mathcal{M},\mu_\mathcal{M})$ that will be used to define the limiting diffusion in our main result. We also introduce some approximations to the process on spaces built from a finite number of line segments.

To do this, it will be enough to work in the deterministic setting. In particular, until noted otherwise we fix a compact real tree $\mathcal{T}=(\mathcal{T},d_\mathcal{T})$ (\cite{LegallDuquesne}, Definition 2.1, for example) containing more than one point, and a distinguished vertex $\rho\in\mathcal{T}$, which we call the root. Let $\mu_\mathcal{T}$ be a finite Borel measure on $\mathcal{T}$ satisfying
\begin{equation}\label{lowervol}
\liminf_{r\rightarrow 0} \frac{\inf_{x\in\mathcal{T}}\mu_\mathcal{T}(B_\mathcal{T}(x,r))}{r^\kappa}>0
\end{equation}
for some $\kappa>0$, where $B_\mathcal{T}(x,r)$ is the open ball of radius $r$ (with respect to the metric $d_\mathcal{T}$) centred at $x\in\mathcal{T}$. Note that this condition guarantees that $\mu_\mathcal{T}$ is of full support. Suppose $(u_i)_{i=1}^J$ and $(v_i)_{i=1}^J$ are finite collections of vertices of $\mathcal{T}$, and write the collection of sets $\{u_i,v_i\}$ as $E=\{\{u_i,v_i\}:i=1,\dots,J\}$. Define $\sim_E$ by setting $x\sim_E y$ if and only if
\begin{eqnarray}
x=y&\mbox{ or }&\exists \{x_i,y_i\}\in E,i=1,\dots,k,\mbox{ such that }x_1=x,y_i=x_{i+1},y_k=y,\label{equiv}
\end{eqnarray}
which is easily checked to be an equivalence relation on $\mathcal{T}$. The canonical projection from $\mathcal{T}$ into the quotient space $\mathcal{M}:=\mathcal{T}/\sim_E$ will be denoted $\phi$, and we will also use the notation $\bar{x}:=\phi(x)$ for $x\in\mathcal{T}$. The space $\mathcal{M}$ becomes a metric space when equipped with the quotient metric
\[d_{\mathcal{M}}(\bar{x},\bar{y})=
\inf\left\{\sum_{i=1}^kd_\mathcal{T}(x_i,y_i):\bar{x}_1=\bar{x},\bar{y}_i=\bar{x}_{i+1},\bar{y}_k=\bar{y},k\in\mathbb{N}\right\}\]
(\cite{BBI}, Exercise 3.1.13), and we will also define $\mu_\mathcal{M}:=\mu_\mathcal{T}\circ\phi^{-1}$. Note that we allow the possibility that $J=0$. In this case, we simply have that $(\mathcal{M},d_\mathcal{M},\mu_\mathcal{M})$ and $(\mathcal{T},d_\mathcal{T},\mu_\mathcal{T})$ are identical as metric-measure spaces.

Although it would be possible to define the form of interest, $(\mathcal{E}_\mathcal{M},\mathcal{F}_\mathcal{M})$, as the limit of resistance forms (\cite{Kigami}, Definition 2.3) on an increasing sequence of finite approximations to the space $\mathcal{M}$ by applying results of \cite{Kigamidendrite}, a more concise construction is provided by following the steps briefly described in the introduction, i.e. starting with the natural resistance form $(\mathcal{E}_\mathcal{T},\mathcal{F}_\mathcal{T})$ on the space $\mathcal{T}$, which is easily obtained by viewing $\mathcal{T}$ as a resistance network, and then fusing the vertices at $u_i$ and $v_i$ together for each $i=1,\dots,J$.  Note how the latter description neatly complements the understanding of the quotient metric space $(\mathcal{M},d_\mathcal{M})$ as the space $(\mathcal{T},d_\mathcal{T})$ glued along the relation $\sim_E$. More specifically, by \cite{Kigamidendrite}, Theorem 5.4, there exists a unique resistance form $(\mathcal{E}_\mathcal{T},\mathcal{F}_\mathcal{T})$ on $\mathcal{T}$ that satisfies
\begin{equation}\label{dtchar}
d_\mathcal{T}(x,y)=\inf\{\mathcal{E}_\mathcal{T}(f,f):f\in\mathcal{F}_\mathcal{T},f(x)=0,f(y)=1\}^{-1}
\end{equation}
for $x\neq y\in \mathcal{T}$. Furthermore, the same result implies that $(\mathcal{E}_\mathcal{T},\mathcal{F}_\mathcal{T})$ is a local, regular Dirichlet form on $L^2(\mathcal{T},\mu_\mathcal{T})$. Given this form, let $\mathcal{F}_\mathcal{M}:=\{f:\mathcal{M}\rightarrow\mathbb{R}\::\:f_\phi\in\mathcal{F}_\mathcal{T}\}$,
where, given a function $f:\mathcal{M}\rightarrow\mathbb{R}$, we define $f_\phi:=f\circ\phi$, and set
\[\mathcal{E}_\mathcal{M}(f,f):=\mathcal{E}_\mathcal{T}(f_\phi,f_\phi),\hspace{20pt}\forall f\in\mathcal{F}_\mathcal{M}.\]
We will eventually show that $(\mathcal{E}_\mathcal{M},\mathcal{F}_\mathcal{M})$ is  local, regular Dirichlet form on the space $L^2(\mathcal{M},\mu_\mathcal{M})$, but first prove that it is a resistance form on $\mathcal{M}$ and define an associated resistance metric.

{\propn $(\mathcal{E}_\mathcal{M},\mathcal{F}_\mathcal{M})$ is a resistance form on $\mathcal{M}$.}
\begin{proof} The following properties of a resistance form are easily checked from the definition of $(\mathcal{E}_\mathcal{M},\mathcal{F}_\mathcal{M})$ (and the fact that $(\mathcal{E}_\mathcal{T},\mathcal{F}_\mathcal{T})$ is itself resistance form): $\mathcal{F}_\mathcal{M}$ is a linear subspace of functions on $\mathcal{M}$ containing constants; $\mathcal{E}_\mathcal{M}(f,f)=0$ if and only if $f$ is constant on $\mathcal{M}$; after constant functions are quotiented out, then $\mathcal{E}_\mathcal{M}$ is an inner product on $\mathcal{F}_\mathcal{M}$; for any finite subset $V\subseteq \mathcal{M}$ and function $f:V\rightarrow\mathbb{R}$, there exists a function $g\in\mathcal{F}_\mathcal{M}$ such that $g|_V=f$; if $\bar{f}:=(f\vee 0)\wedge 1$ for some $f\in \mathcal{F}_\mathcal{M}$, then   $\bar{f}\in \mathcal{F}_\mathcal{M}$ and $\mathcal{E}_\mathcal{M}(\bar{f},\bar{f})\leq \mathcal{E}_\mathcal{M}(f,f)$. Here, we will merely establish the remaining properties: the inner product space $(\mathcal{F}_\mathcal{M},\mathcal{E}_\mathcal{M})$ is complete and also
\begin{equation}\label{supeq}
\sup\left\{\frac{(f(\bar{x})-f(\bar{y}))^2}{\mathcal{E}_\mathcal{M}(f,f)}:f\in\mathcal{F}_\mathcal{M},\mathcal{E}_\mathcal{M}(f,f)>0\right\}<\infty
\end{equation}
for every $\bar{x},\bar{y}\in\mathcal{M}$.

First, if $(f_n)_{n\geq 1}$ is a Cauchy sequence in $(\mathcal{F}_\mathcal{M},\mathcal{E}_\mathcal{M})$, then  $(f_{n\phi})_{n\geq 1}$ is a Cauchy sequence in $(\mathcal{F}_\mathcal{T},\mathcal{E}_\mathcal{T})$. Thus, since $(\mathcal{E}_\mathcal{T},\mathcal{F}_\mathcal{T})$ is a resistance form, there exists a function $f'\in\mathcal{F}_\mathcal{T}$ such that $\mathcal{E}_\mathcal{T}(f_{n\phi}-f',f_{n\phi}-f')\rightarrow 0$. Noting that $f_{n\phi}(u_i)=f_{n\phi}(v_i)$ for each $i=1,\dots,J$, applying (\ref{dtchar}) yields $(f'(u_i)-f'(v_i))^2\leq d_{\mathcal{T}}(u_i,v_i)\mathcal{E}_\mathcal{T}(f_{n\phi}-f',f_{n\phi}-f')\rightarrow0$. In particular, it follows that $f'(u_i)=f'(v_i)$ for each $i=1,\dots,J$, and so $f'=f_\phi$ for some $f\in\mathcal{F}_\mathcal{M}$. Moreover, $\mathcal{E}_\mathcal{M}(f_{n}-f,f_{n}-f)=\mathcal{E}_\mathcal{T}(f_{n\phi}-f_\phi,f_{n\phi}-f_\phi)\rightarrow0$, which confirms that $(\mathcal{F}_\mathcal{M},\mathcal{E}_\mathcal{M})$ is complete.

For any $\bar{x},\bar{y}\in\mathcal{M}$, we can rewrite the supremum in (\ref{supeq}) as
\begin{eqnarray}
\lefteqn{\sup\left\{\frac{(f_\phi({x})-f_\phi({y}))^2}{\mathcal{E}_\mathcal{T}(f_\phi,f_\phi)}:f\in\mathcal{F}_\mathcal{M},\mathcal{E}_\mathcal{T}(f_\phi,f_\phi)>0\right\}}\nonumber\\
&\leq&\sup\left\{\frac{(f({x})-f({y}))^2}{\mathcal{E}_\mathcal{T}(f,f)}:f\in\mathcal{F}_\mathcal{T},\mathcal{E}_\mathcal{T}(f,f)>0\right\}\label{rmdt}\\
&=&d_\mathcal{T}(x,y),\nonumber
\end{eqnarray}
which is finite. Note that, to deduce the above equality, we have again applied the characterisation of $(\mathcal{E}_\mathcal{T},\mathcal{F}_\mathcal{T})$ at (\ref{dtchar}).
\end{proof}

Given that $(\mathcal{E}_\mathcal{M},\mathcal{F}_\mathcal{M})$ is a resistance form, by \cite{Kigami}, Theorem 2.3.4, if we define a function $R_\mathcal{M}:\mathcal{M}\times \mathcal{M}\rightarrow \mathbb{R}$ by setting $R_\mathcal{M}(\bar{x},\bar{y})$ to be equal to the supremum at (\ref{supeq}), then $R_\mathcal{M}$ is a metric on $\mathcal{M}$. We call this metric the resistance metric on $\mathcal{M}$, and the next result shows that it is equivalent to the quotient metric $d_\mathcal{M}$. In the proof we will use the notation $b^\mathcal{T}(x,y,z)$ to represent the branch-point of $x,y,z\in \mathcal{T}$, which is the unique point satisfying
\begin{equation}\label{bt}
\left\{b^\mathcal{T}(x,y,z)\right\}=[[x,y]]\cap[[y,z]]\cap[[z,x]],
\end{equation}
where, for two vertices $x,y\in\mathcal{T}$, $[[x,y]]$ is the unique (non-self intersecting) path from $x$ to $y$ in the real tree $\mathcal{T}$. Furthermore, for a form $(\mathcal{E},\mathcal{F})$ defined on a set $A$, the trace onto $B\subseteq A$, which will be denoted ${\rm Tr}(\mathcal{E}|B)$, satisfies
\begin{equation}\label{trace}
{\rm Tr}(\mathcal{E}|B)(f,f)=\inf\{\mathcal{E}(g,g):g\in\mathcal{F},g|_B=f\},
\end{equation}
with the domain of ${\rm Tr}(\mathcal{E}|B)$ being the collection of functions $f:B\rightarrow\mathbb{R}$ such that the right-hand side above is finite.

{\lem \label{compmet} There exists a strictly positive constant $c$, depending only on $J$, such that
\[cd_{\mathcal{M}}(\bar{x},\bar{y})\leq R_\mathcal{M}(\bar{x},\bar{y})\leq d_{\mathcal{M}}(\bar{x},\bar{y}),\]
for every $\bar{x},\bar{y}\in\mathcal{M}$.}
\begin{proof} For the upper bound, let $\bar{x},\bar{y}\in\mathcal{M}$ and $k\in\mathbb{N}$, and suppose $x_{i},y_i\in\mathcal{T}$, $i=1,\dots,k$, are vertices satisfying $\bar{x}_1=\bar{x}$, $\bar{y}_i=\bar{x}_{i+1}$, $\bar{y}_k=\bar{y}$. Applying the triangle inequality for $R_\mathcal{M}$, we have that
\[R_\mathcal{M}(\bar{x},\bar{y})\leq  \sum_{i=1}^k R_{\mathcal{M}}(\bar{x}_i,\bar{y}_i)
\leq \sum_{i=1}^k d_{\mathcal{T}}({x}_i,{y}_i),\]
where the second inequality follows from (\ref{rmdt}). Taking the infimum over the sets of sequences satisfying the assumptions yields that $R_\mathcal{M}(\bar{x},\bar{y})\leq d_\mathcal{M}(\bar{x},\bar{y})$, as required.

We now establish the lower bound. Fix $x,y\in\mathcal{T}$, and set
\[V=\left\{b^\mathcal{T}(u,v,w):u,v,w\in \{x,y,u_1,\dots,u_J,v_1,\dots,v_J\}\right\},\]
where we note that $\{x,y,u_1,\dots,u_J,v_1,\dots,v_J\}\subseteq V$. Moreover, let $U=V/\sim_E$, where $\sim_E$ is the equivalence defined at (\ref{equiv}). Writing $u\leftrightarrow v$ to signify that $u\neq v\in V$ satisfy $[[u,v]]\cap V=\{u,v\}$, which is a formalisation of the notion that $u$ and $v$ are neighbours in $V$, we define an electrical network with vertex set $U$ by supposing that vertices $\bar{u}$ and $\bar{v}$ are connected by wires with resistances $(d_{\mathcal{T}}(u,v))_{u\in\bar{u},v\in\bar{v},u\leftrightarrow v}$ (and not directly connected by a wire if there are no $u,v\in V$ satisfying $u\in\bar{u}$, $v\in\bar{v}$, $u\leftrightarrow v$). If the vertices in this network are held at potential $f:U\rightarrow \mathbb{R}$, then the total energy dissipation (see \cite{DS}, Section 1.3.5) is given by
\begin{equation}\label{finsum}
\mathcal{E}_U(f,f)=\frac{1}{2}\sum_{\bar{u},\bar{v}\in U}\sum_{\substack{u\in\bar{u},v\in\bar{v},\\u\leftrightarrow v}}\frac{(f(\bar{u})-f(\bar{v}))^2}{d_\mathcal{T}(u,v)}=\frac{1}{2}\sum_{\substack{u,v\in V,\\u\leftrightarrow v}}\frac{(f_\phi({u})-f_\phi({v}))^2}{d_\mathcal{T}(u,v)}.
\end{equation}
Using ideas from \cite{Kigamidendrite} (in particular, see proof of Proposition 5.1 and Corollary 1.8), we observe that the final sum here can be rewritten ${\rm Tr}(\mathcal{E}_\mathcal{T}|V)(f_\phi,f_\phi)$. Hence, if $R_U:U\times U\rightarrow \mathbb{R}$ is the effective resistance between vertices in the electrical network described above, then
\begin{eqnarray}
&&\label{finsum2}\\
\lefteqn{R_U(\bar{x},\bar{y})^{-1}}\nonumber\\
&=&\inf\{\mathcal{E}_U(f,f):f:U\rightarrow\mathbb{R},f(\bar{x})=0,f(\bar{y})=1\}\nonumber\\
&=&\inf\{{\rm Tr}(\mathcal{E}_\mathcal{T}|V)(f,f):f:V\rightarrow\mathbb{R},f({x})=0,f({y})=1,f(u_i)=f(v_i),i=1,\dots,J\}\nonumber\\
&=&\inf\{\mathcal{E}_\mathcal{T}(f,f):f\in\mathcal{F}_\mathcal{T},f({x})=0,f({y})=1,f(u_i)=f(v_i),i=1,\dots,J\}\nonumber\\
&=&\inf\{\mathcal{E}_\mathcal{M}(f,f):f\in\mathcal{F}_\mathcal{M},f(\bar{x})=0,f(\bar{y})=1\}\nonumber\\
&=&R_\mathcal{M}(\bar{x},\bar{y})^{-1},\nonumber
\end{eqnarray}
where the first equality is an application of the Dirichlet principle for electrical networks (see \cite{DS}, Exercise 1.3.11, for example), the third equality follows from (\ref{trace}), and the final two equalities are consequences of the definitions of $(\mathcal{E}_\mathcal{M},\mathcal{F}_\mathcal{M})$ and $R_\mathcal{M}$ respectively. Thus, to complete the proof, it suffices to show that $R_U(\bar{x},\bar{y})$ is bounded below by $cd_\mathcal{M}(\bar{x},\bar{y})$ for some strictly positive constant $c$ depending only on $J$. Combining the general resistance lower bound for finite electrical networks of Lemma \ref{reslower} and the definition of $d_\mathcal{M}$, it is straightforward to demonstrate that this is the case with $c=1/(4J+1)!$.
\end{proof}

We are now ready to demonstrate that $(\mathcal{E}_\mathcal{M},\mathcal{F}_\mathcal{M})$ is a Dirichlet form. In proving that this form is regular, we will need to consider the collection of continuous functions $\mathcal{M}$ with respect to the metric $d_\mathcal{M}$. Note that, by the previous result, this is the same as the collection of continuous functions with respect to $R_\mathcal{M}$. We will denote the relevant set by $C(\mathcal{M})$, and observe that, by the definition of $R_\mathcal{M}$, we have that
\begin{equation}\label{holder}
\left(f(\bar{x})-f(\bar{y})\right)^2\leq R_{\mathcal{M}}(\bar{x},\bar{y})\mathcal{E}_\mathcal{M}(f,f),\hspace{20pt}\forall f\in\mathcal{F}_\mathcal{M},\bar{x},\bar{y}\in\mathcal{M},
\end{equation}
which implies $\mathcal{F}_\mathcal{M}\subseteq C(\mathcal{M})$. Moreover, the compactness of $\mathcal{M}$ and finiteness of $\mu_\mathcal{M}$ yield that $C(\mathcal{M})\subseteq L^2(\mathcal{M},\mu_\mathcal{M})$.

{\propn  $(\mathcal{E}_\mathcal{M},\mathcal{F}_\mathcal{M})$ is  local, regular Dirichlet form on $L^2(\mathcal{M},\mu_\mathcal{M})$.}
\begin{proof} We start by showing that $\mathcal{F}_\mathcal{M}$ is dense in $C(\mathcal{M})$ with respect to the supremum metric $\|\cdot\|_\infty$. If $f\in C(\mathcal{M})$, then $f_\phi$ is a continuous function on $\mathcal{T}$ (with respect to $d_\mathcal{T}$), and so, by the regularity of $(\mathcal{E}_\mathcal{T},\mathcal{F}_\mathcal{T})$, there exists a sequence $(g_n)_{n\geq 1}$ in $\mathcal{F}_\mathcal{T}$ such that $\|g_n-f_\phi\|_\infty\rightarrow 0$. To continue, set $V:=\{u_1,\dots,u_J,v_1,\dots,v_J\}$ and, for each $x\in V$, let $h_x$ be the harmonic extension of $\mathbf{1}_{\{x\}}$ from $V$ to $\mathcal{T}$. More precisely, $h_x$ is the unique function in $\mathcal{F}_\mathcal{T}$ satisfying
$h_x|_V=\mathbf{1}_{\{x\}}$ and also $\mathcal{E}_\mathcal{T}(h_x,h_x)={\rm Tr}(\mathcal{E}_\mathcal{T}|V)(\mathbf{1}_{\{x\}},\mathbf{1}_{\{x\}})$.
The existence of $h_x$ is guaranteed by \cite{Kigami}, Lemma 3.5, and moreover, \cite{Kigami}, Theorem 1.4, demonstrates that $0\leq h_x\leq 1$ everywhere on $\mathcal{T}$. Now define
\[f_n':=g_n+\sum_{x\in V}(f_\phi(x)-g_n(x))h_x.\]
Clearly $f_n'\in\mathcal{F}_\mathcal{T}$ and $f_n'(u_i)=f_\phi(u_i)=f_\phi(v_i)=f_n'(v_i)$ for each $i=1,\dots,J$. Hence there exists an $f_n\in\mathcal{F}_\mathcal{M}$ such that $f_n'=f_{n\phi}$. Furthermore,
\begin{eqnarray*}
\|f_n-f\|_\infty&=&\|f_{n\phi}-f_\phi\|_\infty\\
&\leq&\|g_{n}-f_\phi\|_\infty+\sum_{x\in V}|f_\phi(x)-g_n(x)|\\
&\leq&(1+\#V)\|g_{n}-f_\phi\|_\infty\\
&\rightarrow& 0,
\end{eqnarray*}
from which we obtain the desired conclusion.

It readily follows from the result of the previous paragraph that $\mathcal{F}_\mathcal{M}$ is also a dense subset of $L^2(\mathcal{M},\mu_\mathcal{M})$. Given this fact and that $(\mathcal{E}_\mathcal{M},\mathcal{F}_\mathcal{M})$ is a resistance form, to establish that the latter is also a Dirichlet form is relatively straightforward. One point that does require checking is that $(\mathcal{F}_\mathcal{M},\mathcal{E}_\mathcal{M}^1)$ is a Hilbert space, where
\begin{equation}\label{hilbnorm}
\mathcal{E}_\mathcal{M}^1(f,f):=\mathcal{E}_\mathcal{M}(f,f)+\int_\mathcal{M}|f|^2d\mu_\mathcal{M}
\end{equation}
for $f\in\mathcal{F}_\mathcal{M}$, but this can be done by following the proof of \cite{Kigami}, Theorem 2.4.1. Since we already know that $\mathcal{F}_\mathcal{M}$ is a dense subset of $C(\mathcal{M})$, it is further the case that the Dirichlet form $(\mathcal{E}_\mathcal{M},\mathcal{F}_\mathcal{M})$ is regular.

Finally, to prove $(\mathcal{E}_\mathcal{M},\mathcal{F}_\mathcal{M})$ is local it is necessary to show that $\mathcal{E}_\mathcal{M}(f,g)=0$ for any two functions $f,g\in\mathcal{F}_\mathcal{M}$ with disjoint support. However, it is clear that if $f$ and $g$ have disjoint support in $\mathcal{M}$, then $f_\phi$ and $g_\phi$ have disjoint support in $\mathcal{T}$. Since $(\mathcal{E}_\mathcal{T},\mathcal{F}_\mathcal{T})$ is local, it follows that $\mathcal{E}_\mathcal{M}(f,g)=\mathcal{E}_\mathcal{T}(f_\phi,g_\phi)=0$, as desired.
\end{proof}

Now we have the local, regular Dirichlet form $(\tfrac{1}{2}\mathcal{E}_\mathcal{M},\mathcal{F}_\mathcal{M})$ on $L^2(\mathcal{M},\mu_\mathcal{M})$, standard arguments (\cite{FOT}, Theorem 7.2.2) imply the existence of a corresponding $\mu_\mathcal{M}$-symmetric Markov diffusion
\[X^\mathcal{M}=\left(\left(X^\mathcal{M}_t\right)_{t\geq 0},\mathbf{P}^\mathcal{M}_{\bar{x}},\bar{x}\in\mathcal{M}\right),\]
which we will call Brownian motion on $\mathcal{M}$. The factor of $\tfrac{1}{2}$ here appears a little awkward, but will later ensure the correct time scaling of the process. Note that, since every point $\bar{x}\in\mathcal{M}$ has strictly positive capacity (see \cite{KigRes}, Theorem 8.8), the process $X^\mathcal{M}$ is uniquely determined (\cite{FOT}, Theorem 4.2.1 and Theorem 4.2.7). Moreover, since $(\mathcal{E}_\mathcal{M},\mathcal{F}_\mathcal{M})$ is a resistance form, it is irreducible and recurrent as a Dirichlet form, and so we can apply \cite{FOT}, Theorem 4.6.6, to deduce that $X^\mathcal{M}$ hits points in the sense that
$\mathbf{P}^\mathcal{M}_{\bar{x}}\left(h(\bar{y})<\infty\right)=1$,
for every $\bar{x},\bar{y}\in\mathcal{M}$, where $h(\bar{y}):=\inf\{t\geq 0:X^\mathcal{M}_t=\bar{y}\}$ is the hitting time of $\bar{y}$. So far, we have not fully applied the lower volume asymptotics of (\ref{lowervol}). The importance of this bound is that it will allow us to apply Theorem 6.3 of \cite{MarcusRosen} to deduce the existence of jointly continuous local times for $X^\mathcal{M}$. Since the following result can be proved exactly as in the real tree case (cf. \cite{Croydoncbp}, Lemma 2.5, and \cite{Croydoninf}, Lemmas 2.2 and 2.3), we omit its proof. Note that to use \cite{MarcusRosen}, Theorem 6.3, it is required that the process $X^\mathcal{M}$ is strongly symmetric (in the sense of \cite{MarcusRosen}), but this is an easy consequence of the fact that $X^\mathcal{M}$ admits a transition density $(p^\mathcal{M}_t(\bar{x},\bar{y}))_{t>0,\bar{x},\bar{y}\in\mathcal{M}}$ that is jointly continuous in $(t,\bar{x},\bar{y})$ (see \cite{KigRes}, Theorem 10.4).

{\lem \label{ltcont} The process $X^\mathcal{M}$ admits local times $(L^\mathcal{M}_t(\bar{x}))_{t\geq0, \bar{x}\in\mathcal{M}}$ that are $\mathbf{P}^\mathcal{M}_{\bar{y}}$-a.s. jointly continuous in $t$ and $\bar{x}$, for every $\bar{y}\in\mathcal{M}$. Furthermore, $\mathbf{P}^\mathcal{M}_{\bar{x}}$-a.s. we have
\[\lim_{t\rightarrow\infty}\inf_{\bar{y}\in\mathcal{M}}L^\mathcal{M}_t(\bar{y})=\infty,\]
for every $\bar{x}\in\mathcal{M}$.}
\bigskip

To complete this section, we transfer ideas developed in \cite{Croydoncbp} and \cite{Croydoninf} for approximating the Brownian motion on a compact real tree by Brownian motions on subtrees with a finite number of branches to the current setting. First, extend $(u_i)_{i=1}^J$ to a dense sequence $(u_i)_{i=1}^\infty$ in $\mathcal{T}$. Without loss of generality, we assume that $u_{J+1}\neq \rho$. For each $k\geq J+1$, define
\begin{equation}\label{tkdef}
\mathcal{T}(k):=\left(\cup_{i=1}^k[[\rho,u_i]]\right)\cup\left(\cup_{i=1}^J[[\rho,v_i]]\right),
\end{equation}
and set
\begin{equation}\label{mkdef}
\mathcal{M}(k):=\phi(\mathcal{T}(k)).
\end{equation}
We will consider two measures on $\mathcal{M}(k)$. One will be the one-dimensional Hausdorff measure, $\lambda_{\mathcal{M}(k)}$. The other will be image under $\phi$ of the natural projection of $\mu_\mathcal{T}$ onto $\mathcal{T}(k)$. More specifically, set
\begin{equation}\label{mumkdef}
\mu_{\mathcal{M}(k)}:=\mu_{\mathcal{T}(k)}\circ \phi^{-1},
\end{equation}
where $\mu_{\mathcal{T}(k)}:=\mu_\mathcal{T}\circ\phi_{\mathcal{T},\mathcal{T}(k)}^{-1}$ and, as in \cite{Croydoncbp}, the projection map $\phi_{\mathcal{T},\mathcal{T}(k)}$ is defined by setting, for $x\in\mathcal{T}$, $\phi_{\mathcal{T},\mathcal{T}(k)}(x)$  to be the unique point in $\mathcal{T}(k)$ satisfying
\begin{equation}\label{ttk}
d_\mathcal{T}(x,\phi_{\mathcal{T},\mathcal{T}(k)}(x))=\inf_{y\in\mathcal{T}(k)}d_\mathcal{T}(x,y).
\end{equation}
To check that $\lambda_{\mathcal{M}(k)}$ and $\mu_{\mathcal{M}(k)}$ are well-defined finite Borel (with respect to $d_\mathcal{M}$) measures with support equal to $\mathcal{M}(k)$ is straightforward from their construction. Moreover, since $\mu_{\mathcal{T}}\circ\phi_{\mathcal{T},\mathcal{T}(k)}^{-1}\rightarrow\mu_\mathcal{T}$ weakly as measures on $\mathcal{T}$ (cf. \cite{Croydoninf}, Section 2), the continuity of $\phi$ implies that $\mu_{\mathcal{M}(k)}\rightarrow \mu_\mathcal{M}$ weakly as measures on $\mathcal{M}$.

Since $\mathcal{M}(k)$ is a non-empty subset of $\mathcal{M}$, it follows from \cite{KigRes}, Theorem 8.4 that if $(\mathcal{E}_{\mathcal{M}(k)},\mathcal{F}_{\mathcal{M}(k)})$ is defined by setting
\[\mathcal{E}_{\mathcal{M}(k)}:={\rm Tr}(\mathcal{E}_{\mathcal{M}}|\mathcal{M}(k))\]
and taking $\mathcal{F}_{\mathcal{M}(k)}$ to be the domain of $\mathcal{E}_{\mathcal{M}(k)}$, then $(\mathcal{E}_{\mathcal{M}(k)},\mathcal{F}_{\mathcal{M}(k)})$ is a resistance form on $\mathcal{M}(k)$ and the associated resistance metric is  $R_\mathcal{M}|_{\mathcal{M}(k)\times\mathcal{M}(k)}$. We have the following alternative description of $(\mathcal{E}_{\mathcal{M}(k)},\mathcal{F}_{\mathcal{M}(k)})$.

{\lem Fix $k\geq J+1$. If $(\mathcal{E}_{\mathcal{T}(k)},\mathcal{F}_{\mathcal{T}(k)})$ is the resistance form associated with the real tree $(\mathcal{T}(k),d_\mathcal{T})$ by \cite{Kigamidendrite}, Theorem 5.4, then $\mathcal{F}_{\mathcal{M}(k)}=\{f:\mathcal{M}(k)\rightarrow\mathbb{R}\::\:f_\phi\in\mathcal{F}_{\mathcal{T}(k)}\}$ and
\begin{equation}\label{newchar}
\mathcal{E}_{\mathcal{M}(k)}(f,f)=\mathcal{E}_{\mathcal{T}(k)}(f_\phi,f_\phi),\hspace{20pt}\forall f\in\mathcal{F}_{\mathcal{M}(k)}.
\end{equation}}
\begin{proof} Fix $k\geq J+1$. Since $u_i,v_i\in\mathcal{T}(k)$, $i=1,\dots,J$, it is an elementary exercise to prove that, for a function $f:\mathcal{M}(k)\rightarrow \mathbb{R}$,
\[\{g_\phi:g\in\mathcal{F}_{\mathcal{M}},g|_{\mathcal{M}(k)}=f\}= \{h:h\in\mathcal{F}_{\mathcal{T}},h|_{\mathcal{T}(k)}=f_\phi\}.\]
Thus, for $f\in\mathcal{F}_{\mathcal{M}(k)}$,
\begin{eqnarray*}
\mathcal{E}_{\mathcal{M}(k)}(f,f)&=&\inf\{\mathcal{E}_{\mathcal{M}}(g,g):g\in\mathcal{F}_{\mathcal{M}},g|_{\mathcal{M}(k)}=f\}\\
&=&\inf\{\mathcal{E}_{\mathcal{T}}(h,h):h\in\mathcal{F}_{\mathcal{T}},h|_{\mathcal{T}(k)}=f_\phi\}\\
&=&\mathcal{E}_{\mathcal{T}(k)}(f_\phi,f_\phi),
\end{eqnarray*}
where the final equality is a consequence of the fact that $\mathcal{E}_{\mathcal{T}(k)}={\rm Tr}(\mathcal{E}_{\mathcal{T}}|\mathcal{T}(k))$, which can be deduced by first observing that ${\rm Tr}(\mathcal{E}_{\mathcal{T}}|\mathcal{T}(k))$ is a resistance form with resistance metric $d_\mathcal{T}|_{\mathcal{T}(k)\times\mathcal{T}(k)}$ (by \cite{KigRes}, Theorem 8.4) and then noting that resistance forms are uniquely specified by their resistance metrics (\cite{Kigami}, Theorem 2.3.6).
\end{proof}

As with $(\tfrac{1}{2}\mathcal{E}_{\mathcal{M}},\mathcal{F}_{\mathcal{M}})$, we can check that, for any $k\geq J+1$, $(\tfrac{1}{2}\mathcal{E}_{\mathcal{M}(k)},\mathcal{F}_{\mathcal{M}(k)})$ is a Dirichlet form on both $L^2(\mathcal{M}(k),\lambda_{\mathcal{M}(k)})$ and $L^2(\mathcal{M}(k),\mu_{\mathcal{M}(k)})$, and we can use the fact that $\mathcal{M}(k)$ is closed to establish that the form is regular (\cite{KigRes}, Theorem 8.4).  We will denote the (unique) associated Hunt processes by $X^{\lambda_{\mathcal{M}(k)}}$, $X^{\mu_{\mathcal{M}(k)}}$ and their laws starting from $\bar{x}\in\mathcal{M}(k)$ by $\mathbf{P}^{\lambda_{\mathcal{M}(k)}}_{\bar{x}}$, $\mathbf{P}^{\mu_{\mathcal{M}(k)}}_{\bar{x}}$ respectively. Note that the above characterisation of $(\mathcal{E}_{\mathcal{M}(k)},\mathcal{F}_{\mathcal{M}(k)})$ readily yields that $(\mathcal{E}_{\mathcal{M}(k)},\mathcal{F}_{\mathcal{M}(k)})$ is local, exactly as in the proof of the corresponding result for $(\mathcal{E}_{\mathcal{M}},\mathcal{F}_{\mathcal{M}})$, and so these processes are actually diffusions. For the laws of the processes $X^{\mu_{\mathcal{M}(k)}}$ we are able to prove the following convergence result. Since it is a relatively simple adaptation of the proof of \cite{Croydoncbp}, Lemma 3.1, we only sketch the proof.

{\propn \label{kinfty} As $k\rightarrow \infty$,
\[\mathbf{P}_{\bar{\rho}}^{\mu_{\mathcal{M}(k)}}\rightarrow \mathbf{P}_{\bar{\rho}}^{\mathcal{M}}\]
weakly as probability measures on $C(\mathbb{R}_+,\mathcal{M})$.}
\begin{proof} Applying the weak convergence of  $\mu_{\mathcal{M}(k)}$ to $\mu_\mathcal{M}$ and the joint continuity of the local times of $X^{\mathcal{M}}$ (see Lemma \ref{ltcont}), we obtain for every $t\geq 0$ that, $\mathbf{P}_{\bar{\rho}}^{\mathcal{M}}$-a.s.,
\[\tilde{A}^{\mathcal{M}(k)}_t:=\int_{\mathcal{T}(k)}L_t^\mathcal{M}(\bar{x})\mu_{\mathcal{M}(k)}(d\bar{x})\rightarrow t.\]
Moreover, an elementary monotonocity argument yields this convergence result uniformly on compact intervals. As a consequence of this, $\tilde{\tau}^{\mathcal{M}(k)}(t):=\inf\{s:\tilde{A}^{\mathcal{M}(k)}_s>t\}\rightarrow t$ uniformly on compacts, $\mathbf{P}_{\bar{\rho}}^{\mathcal{M}}$-a.s. Now, the trace theorem for Dirichlet forms (see \cite{FOT}, Theorem 6.2.1, for example) allows one to check that the law of
\[(X_{\tilde{\tau}^{\mathcal{M}(k)}(t)}^{\mathcal{M}})_{t\geq0}\]
under $\mathbf{P}_{\bar{\rho}}^{\mathcal{M}}$ is precisely $\mathbf{P}_{\bar{\rho}}^{\mu_{\mathcal{M}(k)}}$ (cf. \cite{Croydoncbp}, Lemma 2.6), and hence the result follows.
\end{proof}

We now describe how the two processes $X^{\lambda_{\mathcal{M}(k)}}$ and  $X^{\mu_{\mathcal{M}(k)}}$ can be coupled. Similarly to above, it is possible to deduce the existence of jointly continuous local times $(L^{\mathcal{M}(k)}_t(\bar{x}))_{t\geq0, \bar{x}\in\mathcal{M}(k)}$ for $X^{\lambda_{\mathcal{M}(k)}}$. Use these to define a continuous additive functional $\hat{A}^{\mathcal{M}(k)}=(\hat{A}^{\mathcal{M}(k)}_t)_{t\geq0}$ by setting
\begin{equation}\label{hata}
\hat{A}^{\mathcal{M}(k)}_t:=\int_{\mathcal{M}(k)}L^{\mathcal{M}(k)}_t(\bar{x})\mu_{\mathcal{M}(k)}(d\bar{x}),
\end{equation}
and its inverse by $\hat{\tau}^{\mathcal{M}(k)}(t):=\inf\{s:\hat{A}^{\mathcal{M}(k)}_s>t\}$. As with the time-change employed in the proof of the previous result, the following lemma is a straightforward consequence of the trace theorem for Dirichlet forms, and so will be stated without proof.

{\lem \label{timechange} Fix $k\geq J+1$. If the process $X^{\lambda_{\mathcal{M}(k)}}$ has law  $\mathbf{P}^{\lambda_{\mathcal{M}(k)}}_{\bar{\rho}}$, then the process
\[\left(X^{\lambda_{\mathcal{M}(k)}}_{\hat{\tau}^{\mathcal{M}(k)}(t)}\right)_{t\geq 0}\]
has law $\mathbf{P}_{\bar{\rho}}^{\mu_{\mathcal{M}(k)}}$.}
\bigskip

Finally, we state a result about the paths of $\hat{A}^{\mathcal{M}(k)}$. Applying the continuity and uniform divergence of the local times of $X^\mathcal{M}$, as stated in Lemma \ref{ltcont}, this can be proved identically to \cite{Croydoninf}, Lemma 2.5, and so we once again omit the proof. The scaling factor
\begin{equation}\label{Lambdak}
\Lambda^{(k)}:=\lambda_{\mathcal{M}(k)}(\mathcal{M}(k))
\end{equation}
arises here as a result of the fact that we have not normalised the measure $\lambda_{\mathcal{M}(k)}$.

{\lem \label{akprops} For each $k\geq J+1$, $\mathbf{P}_{\bar{\rho}}^{\lambda_{\mathcal{M}(k)}}$-a.s., the function $\hat{A}^{\mathcal{M}(k)}$ is continuous and strictly increasing. Moreover, for every $t_0\geq0$,
\[\lim_{t\rightarrow\infty}\limsup_{k\rightarrow\infty}\mathbf{P}_{\bar{\rho}}^{\lambda_{\mathcal{M}(k)}}\left(\hat{A}^{\mathcal{M}(k)}_{t_0\Lambda^{(k)}}>t\right)=0.\]
}

\section{Continuity under perturbations}\label{seccont}

The aim of this section is to show how the process $X^{\lambda_{\mathcal{M}(k)}}$ and related continuous additive functional $\hat{A}^{\mathcal{M}(k)}$ are affected continuously by perturbations of the metric $d_{\mathcal{T}}|_{\mathcal{T}(k)\times\mathcal{T}(k)}$ and measure $\mu_{\mathcal{M}(k)}$, where we continue to work in the deterministic framework introduced in the previous section. Throughout, we suppose that $k\geq J+1$ is fixed.

Our first main assumption is that $(d_{\mathcal{T}}^n)_{n\geq 1}$ is a sequence of metrics on $\mathcal{T}(k)$ for which $(\mathcal{T}(k),d_\mathcal{T}^n)$ is a real tree that also satisfy
\begin{equation}\label{compare}
\delta_n d_{\mathcal{T}}^n(x,y)\leq d_\mathcal{T}(x,y)\leq \delta_n^{-1}d^n_\mathcal{T}(x,y),\hspace{20pt}\forall x,y\in\mathcal{T}(k),
\end{equation}
where $(\delta_n)_{n\geq 1}$ is a sequence in $(0,1]$ that converges to 1. For each $n$, the metric $d_{\mathcal{T}}^n$ induces a new quotient metric $d_\mathcal{M}^n$ and one-dimensional Hausdorff measure $\lambda_{\mathcal{M}(k)}^{n}$ on $\mathcal{M}(k)$. Clearly $d_\mathcal{M}^n$ and $d_\mathcal{M}$ satisfy a comparability property analogous to (\ref{compare}). By this equivalence of metrics, $\lambda_{\mathcal{M}(k)}^{n}$ is a finite, non-zero Borel measure on $(\mathcal{M}(k),d_\mathcal{M})$ and also satisfies
\begin{equation}\label{mucompare}
\delta_n\lambda_{\mathcal{M}(k)}^n(A)\leq \lambda_{\mathcal{M}(k)}(A)\leq \delta_n^{-1}\lambda_{\mathcal{M}(k)}^n(A),
\end{equation}
for any measurable $A\subseteq\mathcal{M}$. Finally, since $(\mathcal{T}(k),d_\mathcal{T}^n)$ is a compact real tree, we can again apply \cite{Kigamidendrite}, Theorem 5.4, to define a corresponding resistance form $(\mathcal{E}_{\mathcal{T}(k)}^n,\mathcal{F}_{\mathcal{T}(k)}^n)$, and we use this to characterise $(\mathcal{E}_{\mathcal{M}(k)}^n,\mathcal{F}_{\mathcal{M}(k)}^n)$ through a relation similar to (\ref{newchar}). From this construction, it is clear that (\ref{compare}) yields $\mathcal{F}_{\mathcal{M}(k)}^n=\mathcal{F}_{\mathcal{M}(k)}$ and also
\begin{equation}\label{ecompare}
\delta_n\mathcal{E}_{\mathcal{M}(k)}^n(f,f)\leq \mathcal{E}_{\mathcal{M}(k)}(f,f)\leq\delta_n^{-1}\mathcal{E}_{\mathcal{M}(k)}^n(f,f),\hspace{20pt}\forall f\in\mathcal{F}_{\mathcal{M}(k)},
\end{equation}
for each $n$ (one way of checking this is to apply the finite approximation result of \cite{Kigamidendrite}, Lemma 3.7, in combination with the expression for the Laplacian on a ``fine'' finite subset of a dendrite used in the proof of \cite{Kigamidendrite}, Proposition 5.1, for example).

We will denote by $X^{\lambda_{\mathcal{M}(k)}^n}$ the process associated with  $(\tfrac{1}{2}\mathcal{E}_{\mathcal{M}(k)}^n,\mathcal{F}_{\mathcal{M}(k)})$ considered as a regular Dirichlet form on $L^2(\mathcal{M}(k),\lambda_{\mathcal{M}(k)}^n)$, and its law started from $\bar{x}\in\mathcal{M}(k)$ will be written $\mathbf{P}^{\lambda_{\mathcal{M}(k)}^n}_{\bar{x}}$. Note that, from (\ref{ecompare}) and the fact that the form $(\mathcal{E}_{\mathcal{M}(k)},\mathcal{F}_{\mathcal{M}(k)})$ is local, we have that the Dirichlet form $(\mathcal{E}_{\mathcal{M}(k)}^n,\mathcal{F}_{\mathcal{M}(k)})$ is local, and so $X^{\lambda_{\mathcal{M}(k)}^n}$ is a diffusion. Its jointly continuous local times, guaranteed by the same argument as was used for the local times of $X^\mathcal{M}$, will be written $(L^{\mathcal{M}(k),n}_t(\bar{x}))_{t\geq0, \bar{x}\in\mathcal{M}(k)}$. Now we have introduced the most significant notation used in this section, we can state the main result that will be proved here.

{\propn \label{convfin} Let $k\geq J+1$. If the process $X^{\lambda_{\mathcal{M}(k)}^n}$ has law  $\mathbf{P}^{\lambda_{\mathcal{M}(k)}^n}_{\bar{\rho}}$ and the process $X^{\lambda_{\mathcal{M}(k)}}$ has law  $\mathbf{P}^{\lambda_{\mathcal{M}(k)}}_{\bar{\rho}}$, then
\[\left(X^{\lambda_{\mathcal{M}(k)}^n},L^{\mathcal{M}(k),n}\right)\rightarrow
\left(X^{\lambda_{\mathcal{M}(k)}},L^{\mathcal{M}(k)}\right)\]
in distribution as $n\rightarrow\infty$ in $C(\mathbb{R}_+,(\mathcal{M}(k),d_\mathcal{M}))\times C(\mathbb{R}_+\times (\mathcal{M}(k),d_\mathcal{M}),\mathbb{R}_+)$.}
\bigskip

We start by proving tightness.

{\lem \label{tightness} Let $k\geq J+1$. If the process $X^{\lambda_{\mathcal{M}(k)}^n}$ has law  $\mathbf{P}^{\lambda_{\mathcal{M}(k)}^n}_{\bar{\rho}}$, then the collection
\[\left\{\left(X^{\lambda_{\mathcal{M}(k)}^n},L^{\mathcal{M}(k),n}\right):n\geq 1\right\}\]
is tight in the space $C(\mathbb{R}_+,(\mathcal{M}(k),d_\mathcal{M}))\times C(\mathbb{R}_+\times (\mathcal{M}(k),d_\mathcal{M}),\mathbb{R}_+)$.}
\begin{proof}
First, note that (\ref{mucompare}) implies that there exists constants $c_1,c_2\in(0,\infty)$ such that
\begin{eqnarray*}
 \inf_{n\geq 1}\inf_{\bar{x}\in\mathcal{M}(k)}\lambda_{\mathcal{M}(k)}^n\left(B_{(\mathcal{M}(k),d_{\mathcal{M}})}(\bar{x},r)\right)&\geq&c_1r,\\ \sup_{n\geq 1}\sup_{\bar{x}\in\mathcal{M}(k)}\lambda_{\mathcal{M}(k)}^n\left(B_{(\mathcal{M}(k),d_{\mathcal{M}})}(\bar{x},r)\right)&\leq& c_2r,
\end{eqnarray*}
for every $r\in(0,1]$.  By applying the argument of \cite{Kumagai}, Lemma 4.2, this implies that
\[\limsup_{n\rightarrow\infty}\sup_{\bar{x}\in\mathcal{M}(k)}
\mathbf{P}_{\bar{x}}^{\lambda_{\mathcal{M}(k)}^n}\left(\inf
\{s:d_{\mathcal{M}}\left(\bar{x},X^{\lambda_{\mathcal{M}(k)}^n}_s\right)>r\}<t\right)\leq c_3e^{-\frac{c_4r^2}{t}}\]
for every $r\in(0,1]$, $t\in(0,t_1]$, for some constants $c_3,c_4,t_1\in(0,\infty)$. Consequently
\[\lim_{t\rightarrow0}\limsup_{n\rightarrow\infty}t^{-1}\sup_{\bar{x}\in\mathcal{M}(k)}\mathbf{P}_{\bar{x}}^{\lambda_{\mathcal{M}(k)}^n}\left(\inf\{s:d_{\mathcal{M}}\left(x,X^{\lambda_{\mathcal{M}(k)}^n}_s\right)>r\}<t\right)=0\]
for any $r>0$, which implies the tightness of $(X^{\lambda_{\mathcal{M}(k)}^n})_{n\geq 1}$ in $C(\mathbb{R}_+,(\mathcal{M}(k),d_\mathcal{M}))$, as required (cf. the corollary to Theorem 7.4 of \cite{Bill2}).

Our argument for local times is an adaptation of the proof of \cite{Croydoninf}, Lemma 3.5, and involves observing $X^{\lambda_{\mathcal{M}(k)}^n}$ on particularly simple subsets of $\mathcal{M}(k)$. Define a finite subset $V:=\left\{b^\mathcal{T}(u,v,w):u,v,w\in \{\rho,u_1,\dots,u_k,v_1,\dots,v_J\}\right\}\subseteq \mathcal{T}(k)$,
where the branch-point function $b^\mathcal{T}$ is defined as at (\ref{bt}), and set
\[\varepsilon_0:=\frac{1}{2}\inf_{n\geq 1}\min_{\substack{x,y\in V:\\x\neq y}}d_\mathcal{T}^n(x,y),\]
which, by (\ref{compare}), is strictly positive. By the definition of $\mathcal{T}(k)$, for this choice of $\varepsilon_0$, it is possible to deduce that for each $n$ there exists a collection of paths $([[a_i,b_i]])_{i\in\mathcal{I}_n}$ covering $\mathcal{T}(k)$ such that $d_{\mathcal{T}}^n(a_i,b_i)=\varepsilon_0$ and
\begin{equation}\label{start}
[[a_i,b_i]]\cap V\subseteq \{a_i,b_i\}.
\end{equation}
Moreover, the collections can be chosen in such a way that $\#\mathcal{I}_n$ is bounded uniformly in $n$. We write $U_i:=\phi([[a_i,b_i]])$, so that $(U_i)_{i\in\mathcal{I}_n}$ is a cover for $\mathcal{M}(k)$ for each $n$.

Applying the notation of the previous paragraph, define, for $i\in\mathcal{I}_n$,
\[A_t^i:=\int_{U_i}L^{\mathcal{M}(k),n}_t(\bar{x})\lambda^n_{\mathcal{M}(k)}(d\bar{x}),\]
and $\tau^i(t):=\inf\{s:{A}^i_s>t\}$. By the trace theorem for Dirichlet forms, if we characterise $X^i$ by setting $X^i_t:=X^{\lambda_{\mathcal{M}(k)}^n}_{\tau^i(t)}$, then, under $\mathbf{P}_{\bar{\rho}}^{\lambda_{\mathcal{M}(k)}^n}$, $X^i$ is the Markov process associated with $(\tfrac{1}{2}\mathcal{E}^i,\mathcal{F}^i)$, where $\mathcal{E}^i:={\rm Tr}(\mathcal{E}_{\mathcal{M}(k)}^n|U_i)$, considered as a Dirichlet form on $L^2(U_i,\lambda^n_{\mathcal{M}(k)}(U_i\cap\cdot))$. Furthermore, the local times of $X^i$ are given by $L^i_t(\bar{x}):=L^{\mathcal{M}(k),n}_{\tau^i(t)}(\bar{x})$
for $t\geq 0$ and $\bar{x}\in U_i$ (cf. \cite{Croydoncbp}, Lemma 3.4). Similarly to the proof of \cite{Croydoninf}, Lemma 3.5, this construction yields the following upper bound, for $\delta\in(0,\varepsilon_0)$ and $t_0<\infty$,
\begin{eqnarray}
\lefteqn{\sup_{\substack{\bar{x},\bar{y}\in\mathcal{M}(k):\\d_\mathcal{M}^n(\bar{x},\bar{y}) \leq \delta}}\sup_{\substack{s,t\in[0,t_0]:\\|s-t|\leq \delta}}\left|L^{\mathcal{M}(k),n}_s(\bar{x})-L^{\mathcal{M}(k),n}_t(\bar{y})\right|}\nonumber\\
&\leq&\sum_{i\in\mathcal{I}_n}\sup_{\substack{\bar{x},\bar{y}\in U_i:\\d_\mathcal{M}^n(\bar{x},\bar{y})\leq \delta}}\sup_{\substack{s,t\in[0,t_0]:\\|s-t|\leq \delta}}\left|L^i_s(\bar{x})-L^i_t(\bar{y})\right|.\label{upperbound23}
\end{eqnarray}

Now, fix $i\in\mathcal{I}_n$, and set $V':=V\cup\{{a}_i,{b}_i\}$ and $U':=\phi(V')$. Supposing $R^n_{\mathcal{M}(k)}$ is the resistance metric associated with $(\mathcal{E}^n_{\mathcal{M}(k)},\mathcal{F}_{\mathcal{M}(k)})$, by taking steps similar to (\ref{finsum}) and (\ref{finsum2}), we can deduce
\begin{eqnarray*}
R^n_{\mathcal{M}(k)}(\bar{a}_i,\bar{b}_i)^{-1}&=&\inf_{\substack{f:U'\rightarrow \mathbb{R},\\f(\bar{a}_i)=0,f(\bar{b}_i)=1}}
\frac{1}{2}\sum_{\substack{u,v\in V',\\u\leftrightarrow v}}\frac{(f_\phi({u})-f_\phi({v}))^2}{d^n_{\mathcal{T}}(u,v)}\\
&=&\inf_{\substack{f:U'\rightarrow \mathbb{R},\\f(\bar{a}_i)=0,f(\bar{b}_i)=1}}
\frac{1}{2}\sum_{\substack{u,v\in V',\:u\leftrightarrow v,\\ \{u,v\}\neq \{a_i,b_i\}}}\frac{(f_\phi({u})-f_\phi({v}))^2}{d^n_{\mathcal{T}}(u,v)}+\varepsilon_0^{-1},
\end{eqnarray*}
where we note that $a_i\leftrightarrow b_i$ (in $V'$) is a consequence of (\ref{start}). In particular, the infimum in the final line here, which represents the conductance from $\bar{a}_i$ to $\bar{b}_i$ in the network $\mathcal{M}(k)$ with the segment $U_i\backslash\{\bar{a}_i,\bar{b}_i\}$ removed, is equal to the difference $R^n_{\mathcal{M}(k)}(\bar{a}_i,\bar{b}_i)^{-1}-\varepsilon_0^{-1}$. Thus, if $V'':=V'\cup\{x,y\}$ and $U'':=\phi(V'')$, where $x,y\in [[a_i,b_i]]$ are such that $d_\mathcal{T}^n(a_i,x)<d_\mathcal{T}^n(a_i,y)$, one can check by a simple rescaling of this infimum that ${R^n_{\mathcal{M}(k)}(\bar{x},\bar{y})^{-1}}$ is equal to
\begin{eqnarray*}
\lefteqn{\inf_{\substack{f:U''\rightarrow \mathbb{R},\\f(\bar{x})=0,f(\bar{y})=1}}\left\{
\frac{1}{2}\sum_{\substack{u,v\in V',\:u\leftrightarrow v\\ \{u,v\}\neq \{a_i,b_i\}}}\frac{(f_\phi({u})-f_\phi({v}))^2}{d^n_{\mathcal{T}}(u,v)}+\frac{f_\phi({a}_i)^2}{d^n_{\mathcal{T}}(a_i,x)}+
\frac{1}{d^n_{\mathcal{T}}(x,y)}+\frac{(1-f_\phi({b}_i))^2}{d^n_{\mathcal{T}}(y,b_i)}\right\}}\\
&=&\inf_{\alpha,\beta\in\mathbb{R}}\left\{ (R^n_{\mathcal{M}(k)}(\bar{a}_i,\bar{b}_i)^{-1}-\varepsilon_0^{-1})(\beta-\alpha)^2
+\frac{\alpha^2}{d^n_{\mathcal{T}}(a_i,x)}+
\frac{1}{d^n_{\mathcal{T}}(x,y)}+\frac{(1-\beta)^2}{d^n_{\mathcal{T}}(y,b_i)}\right\}\\
&=&\frac{1}{r_i-d^n_{\mathcal{T}}(x,y)}+\frac{1}{d^n_{\mathcal{T}}(x,y)},
\end{eqnarray*}
where $r_i:=(R^n_{\mathcal{M}(k)}(\bar{a}_i,\bar{b}_i)^{-1}-\varepsilon_0^{-1})^{-1}+\varepsilon_0$. Moreover, the choice of $\varepsilon_0$ allows $d^n_{\mathcal{T}}(x,y)$ to be replaced by $d_\mathcal{M}^n(\bar{x},\bar{y})$ in the above expression. Hence we have shown that either of the isometries from $(U_i,d_\mathcal{M}^n)$ to $[0,\varepsilon_0]$ (equipped with the Euclidean metric) also map $R^n_{\mathcal{M}(k)}$ to $\delta^{(r_i)}$, where, for $r\in(\varepsilon_0,\infty]$, $\delta^{(r)}$ is a metric on $[0,\varepsilon_0]$ defined by
\[\delta^{(r)}(x,y):=\left(|x-y|^{-1}+(r-|x-y|)^{-1}\right)^{-1},\hspace{20pt}\forall x,y\in[0,\varepsilon_0].\]

Observing that $\delta^{(r)}$ is the resistance metric of (\ref{resr}), general results regarding the uniqueness of resistance forms (for example, \cite{Kigami}, Theorem 2.3.6) can therefore be used to deduce that $X^i$ behaves identically to the trace of a Brownian motion on a circle of length $r_i$ on an arc of length $\varepsilon_0$ (see Section \ref{ltsec}). Noting further that the definition of $\varepsilon_0$ and (\ref{start}) readily imply $r_i\geq 2\varepsilon_0$ and $X^i_0\in \{\bar{a}_i,\bar{b}_i\}$, it follows from this and the bound at (\ref{upperbound23}) that, for $\varepsilon>0$,
\begin{eqnarray*}
\lefteqn{\mathbf{P}_{\bar{\rho}}^{\lambda_{\mathcal{M}(k)}^n}\left(\sup_{\substack{\bar{x},\bar{y}\in\mathcal{M}(k):\\d_\mathcal{M}^n(\bar{x},\bar{y}) \leq \delta}}\sup_{\substack{s,t\in[0,t_0]:\\|s-t|\leq \delta}}\left|L^{\mathcal{M}(k),n}_s(\bar{x})-L^{\mathcal{M}(k),n}_t(\bar{y})\right|>\varepsilon\right)}\\
&\leq &\#\mathcal{I}_n\sup_{r\geq 2\varepsilon_0}\mathbf{P}
\left(\#\mathcal{I}_n\sup_{\substack{{x},{y}\in [0,\varepsilon_0]:\\|x-y|\leq \delta}}\sup_{\substack{s,t\in[0,t_0]:\\|s-t|\leq \delta}}\left|L^{r,\varepsilon_0}_s({x})-L^{r,\varepsilon_0}_t({y})\right|>\varepsilon\right),
\end{eqnarray*}
where $(L^{r,\varepsilon_0}_t(x))_{t\geq 0,x\in[0,\varepsilon_0]}$ denote the jointly continuous local times associated with the trace of a Brownian motion on a circle of length $r$ on an arc of length $\varepsilon_0$ started from $0$, which we assume are built on a probability space with probability measure $\mathbf{P}$ (for a precise definition, see Section \ref{ltsec}). Hence, applying (\ref{compare}), the uniform boundedness of $\#\mathcal{I}_n$ and Lemma \ref{ltlem}, for $\varepsilon>0$,
\[\lim_{\delta\rightarrow 0}\sup_{n\geq 1}\mathbf{P}_{\bar{\rho}}^{\lambda_{\mathcal{M}(k)}^n}\left(\sup_{\substack{\bar{x},\bar{y}\in\mathcal{M}(k):\\d_\mathcal{M}(\bar{x},\bar{y}) \leq \delta}}\sup_{\substack{s,t\in[0,t_0]:\\|s-t|\leq \delta}}\left|L^{\mathcal{M}(k),n}_s(\bar{x})-L^{\mathcal{M}(k),n}_t(\bar{y})\right|>\varepsilon\right)=0,\]
and so the family $(L^{\mathcal{M}(k),n})_{n\geq 1}$ is tight in $C(\mathbb{R}_+\times (\mathcal{M}(k),d_\mathcal{M}),\mathbb{R}_+)$ as required.
\end{proof}

We now consider a time-changed version of $X^{\lambda_{\mathcal{M}(k)}^n}$. For $t\geq 0$, let
\[{A}^{\mathcal{M}(k),n}_t:=\int_{\mathcal{M}(k)}L^{\mathcal{M}(k),n}_t(\bar{x})\lambda_{\mathcal{M}(k)}(d\bar{x}),\]
and set ${\tau}^{\mathcal{M}(k),n}(t):=\inf\{s:{A}^{\mathcal{M}(k),n}_s>t\}$. If $X^{\lambda_{\mathcal{M}(k)}^n}$ has law $\mathbf{P}^{\lambda_{\mathcal{M}(k)}^n}_{\bar{\rho}}$, then by the trace theorem we obtain that
\[\tilde{X}^{\lambda_{\mathcal{M}(k)}^n}_t:=X^{\lambda_{\mathcal{M}(k)}^n}_{{\tau}^{\mathcal{M}(k),n}(t)}\]
defines a version of the process associated with the $(\tfrac{1}{2}\mathcal{E}_{\mathcal{M}(k)}^n,\mathcal{F}_{\mathcal{M}(k)})$ considered as a regular Dirichlet form on $L^2(\mathcal{M}(k),\lambda_{\mathcal{M}(k)})$, started from $\bar{\rho}$. As in \cite{Croydoncbp}, Lemma 3.4, we can check that the local times of this process are given by
\[\tilde{L}^{{\mathcal{M}(k),n}}_t(\bar{x}):={L}^{{\mathcal{M}(k),n}}_{{\tau}^{\mathcal{M}(k),n}(t)}(\bar{x}).\]
The following result confirms $(\tilde{X}^{\lambda_{\mathcal{M}(k)}^n},\tilde{L}^{\mathcal{M}(k),n})$ is a good approximation of the pair $(X^{\lambda_{\mathcal{M}(k)}^n},L^{\mathcal{M}(k),n})$ for large $n$.

{\lem \label{approx} Let $k\geq J+1$, and choose $\varepsilon>0$ and $t_0<\infty$. As $n\rightarrow\infty$,
\[\mathbf{P}^{\lambda_{\mathcal{M}(k)}^n}_{\bar{\rho}}\left(\sup_{t\in[0,t_0]}
d_{\mathcal{M}}\left({X}^{\lambda_{\mathcal{M}(k)}^n}_t,\tilde{X}^{\lambda_{\mathcal{M}(k)}^n}_t\right)+\sup_{t\in[0,t_0]}\sup_{\bar{x}\in\mathcal{M}(k)}
\left|L^{\mathcal{M}(k),n}_t(\bar{x})-\tilde{L}^{\mathcal{M}(k),n}_t(\bar{x})\right|>\varepsilon\right)\rightarrow 0.\]}
\begin{proof} By (\ref{mucompare}) and the definition of ${A}^{\mathcal{M}(k),n}_t$, we have that
\begin{eqnarray*}
\delta_nt& =& \delta_n\int_{\mathcal{M}(k)}L^{\mathcal{M}(k),n}_t(\bar{x})\lambda_{\mathcal{M}(k)}^n(d\bar{x})\\
&\leq &{A}^{\mathcal{M}(k),n}_t\leq \delta_n^{-1}\int_{\mathcal{M}(k)}L^{\mathcal{M}(k),n}_t(\bar{x})\lambda_{\mathcal{M}(k)}^n(d\bar{x})=\delta_n^{-1}t,
\end{eqnarray*}
and so $\delta_nt\leq {\tau}^{\mathcal{M}(k),n}(t)\leq \delta_n^{-1}t$. Thus
\[\sup_{t\in[0,t_0]}
d_{\mathcal{M}}\left({X}^{\lambda_{\mathcal{M}(k)}^n}_t,\tilde{X}^{\lambda_{\mathcal{M}(k)}^n}_t\right)
\leq
\sup_{\substack{s\geq0,t\in[0,t_0]:\\|s-t|\leq t_0 (\delta_n^{-1}-1)}}
d_{\mathcal{M}}\left({X}^{\lambda_{\mathcal{M}(k)}^n}_s,{X}^{\lambda_{\mathcal{M}(k)}^n}_t\right),\]
and a similar bound exists for
\[\sup_{t\in[0,t_0]}\sup_{\bar{x}\in\mathcal{M}(k)}
\left|L^{\mathcal{M}(k),n}_t(\bar{x})-\tilde{L}^{\mathcal{M}(k),n}_t(\bar{x})\right|.\]
The result is therefore a consequence of Lemma \ref{tightness} and the fact that $\delta_n\rightarrow 1$.
\end{proof}

In view of this result, it is clear that to complete the proof of Proposition \ref{convfin} it will suffice to prove the same limit with $(X^{\lambda_{\mathcal{M}(k)}^n},L^{\mathcal{M}(k),n})$ replaced by
$(\tilde{X}^{\lambda_{\mathcal{M}(k)}^n},\tilde{L}^{\mathcal{M}(k),n})$, as we do in the subsequent lemma.

{\lem Let $k\geq J+1$. If the process $X^{\lambda_{\mathcal{M}(k)}^n}$ has law  $\mathbf{P}^{\lambda_{\mathcal{M}(k)}^n}_{\bar{\rho}}$ and the process $X^{\lambda_{\mathcal{M}(k)}}$ has law  $\mathbf{P}^{\lambda_{\mathcal{M}(k)}}_{\bar{\rho}}$, then
\begin{equation}\label{twoc}
\left(\tilde{X}^{\lambda_{\mathcal{M}(k)}^n},\tilde{L}^{\mathcal{M}(k),n}\right)\rightarrow
\left(X^{\lambda_{\mathcal{M}(k)}},L^{\mathcal{M}(k)}\right)
\end{equation}
in distribution as $n\rightarrow\infty$ in  $C(\mathbb{R}_+,(\mathcal{M}(k),d_\mathcal{M}))\times C(\mathbb{R}_+\times (\mathcal{M}(k),d_\mathcal{M}),\mathbb{R}_+)$.}
\begin{proof} We start by adapting the proof of \cite{HM}, Theorem 6.1, to show that the forms $(\mathcal{E}_{\mathcal{M}(k)}^n,\mathcal{F}_{\mathcal{M}(k)})$ Mosco-converge to $(\mathcal{E}_{\mathcal{M}(k)},\mathcal{F}_{\mathcal{M}(k)})$ on $L^2(\mathcal{M}(k),\lambda_{\mathcal{M}(k)})$, by which it is meant that:\\
(a) for every sequence $(f_n)_{n\geq 1}$ converging weakly to $f$ in $L^2(\mathcal{M}(k),\lambda_{\mathcal{M}(k)})$,
\[\liminf_{n\rightarrow\infty}\mathcal{E}_{\mathcal{M}(k)}^n(f_n,f_n)\geq \mathcal{E}_{\mathcal{M}(k)}(f,f);\]
(b) for every $f\in L^2(\mathcal{M}(k),\lambda_{\mathcal{M}(k)})$, there exists a sequence $(f_n)_{n\geq 1}$ converging strongly to $f$ in $L^2(\mathcal{M}(k),\lambda_{\mathcal{M}(k)})$ such that
\[\limsup_{n\rightarrow\infty}\mathcal{E}_{\mathcal{M}(k)}^n(f_n,f_n)\leq \mathcal{E}_{\mathcal{M}(k)}(f,f).\]

Applying (\ref{ecompare}) and taking $f_n=f$ for every $n$, property (b) is immediate. For (a), since $\mathcal{E}_{\mathcal{M}(k)}^n(f_n,f_n)=\infty$ for $f_n\not\in\mathcal{F}_{\mathcal{M}(k)}$, it will suffice to consider a sequence $(f_n)_{n\geq 1}$ in $\mathcal{F}_{\mathcal{M}(k)}$ that satisfies
\[\liminf_{n\rightarrow \infty}\mathcal{E}_{\mathcal{M}(k)}^n(f_n,f_n)<\infty\] and converges weakly to some $f$ in $L^2(\mathcal{M}(k),\lambda_{\mathcal{M}(k)})$. Applying the uniform boundedness principle, this final condition implies that $(f_n)_{n\geq 1}$ is bounded in the space $L^2(\mathcal{M}(k),\lambda_{\mathcal{M}(k)})$. Hence, appealing to (\ref{ecompare}) again,
\[\liminf_{n\rightarrow \infty}\left(\mathcal{E}_{\mathcal{M}(k)}(f_n,f_n)+\int_{\mathcal{M}(k)}f_n(\bar{x})^2\lambda_{\mathcal{M}(k)}(d\bar{x})\right)<\infty.\]
It follows that $(f_n)_{n\geq 1}$ admits a weakly convergent subsequence with respect to the inner product $\mathcal{E}_{\mathcal{M}(k)}^1$ defined from $\mathcal{E}_{\mathcal{M}(k)}$ and $\lambda_{\mathcal{M}(k)}$ similarly to (\ref{hilbnorm}). The limit of this sequence is necessarily identical to $f$ (up to $\lambda_{\mathcal{M}(k)}$-a.e. equivalence), and so we can assume that $f$ is also contained in $\mathcal{F}_{\mathcal{M}(k)}$.  We now proceed in three steps. Firstly, observe that (\ref{ecompare}) implies that
\begin{equation}\label{first}
\liminf_{n\rightarrow\infty}\mathcal{E}_{\mathcal{M}(k)}^n(f_n,f_n)\geq \liminf_{n\rightarrow\infty}\mathcal{E}_{\mathcal{M}(k)}(f_n,f_n).
\end{equation}
Secondly, suppose that $(V_m)_{m\geq 1}$ is an increasing sequence of finite subsets of $\mathcal{M}(k)$ such that $\cup_{m\geq 1}V_m$ is dense in $(\mathcal{M}(k),d_\mathcal{M})$ and let $(f_{n_i})_{i\geq 1}$ be a subsequence for which $\mathcal{E}_{\mathcal{M}(k)}(f_{n_i},f_{n_i})\rightarrow \liminf_{n\rightarrow\infty}\mathcal{E}_{\mathcal{M}(k)}(f_n,f_n)<\infty$. Note that (\ref{holder}) and Lemma \ref{compmet} yields that the collection of functions $(f_{n_i})_{i\geq 1}$ is equicontinuous with respect to $d_\mathcal{M}$. Similarly, $f$ is continuous. Hence the weak convergence of $(f_{n_i})_{i\geq 1}$ also implies pointwise convergence. It follows from this, the definition of the trace and the choice of subsequence that, for every $m\geq 1$,
\begin{equation}\label{second}
\liminf_{n\rightarrow\infty}\mathcal{E}_{\mathcal{M}(k)}(f_n,f_n)\geq \liminf_{i\rightarrow\infty}{\rm Tr}(\mathcal{E}_{\mathcal{M}(k)}|V_m)(f_{n_i},f_{n_i})={\rm Tr}(\mathcal{E}_{\mathcal{M}(k)}|V_m)(f,f),
\end{equation}
where to deduce the equality we also use the fact that ${\rm Tr}(\mathcal{E}_{\mathcal{M}(k)}|V_m)$ is a bilinear form on a finite-dimensional space and is therefore continuous. Thirdly, again applying Lemma \ref{compmet}, this time in conjunction with the resistance form limit result of \cite{Kigami}, Lemma 2.3.8, we find that
\begin{equation}\label{third}
\lim_{m\rightarrow\infty}{\rm Tr}(\mathcal{E}_{\mathcal{M}(k)}|V_m)(f,f)=
\mathcal{E}_{\mathcal{M}(k)}(f,f).
\end{equation}
Combining (\ref{first}), (\ref{second}) and (\ref{third}) yields (a), and completes the proof of Mosco-convergence.

With Mosco-convergence of forms, from \cite{Mosco}, Corollary 2.6.1, we obtain that the associated semigroups also converge in the strong operator topology of the space $L^2(\mathcal{M}(k),\lambda_{\mathcal{M}(k)})$. It follows that, for any finite collection of times $0<t_1<\cdots<t_m<\infty$ and bounded functions $f_1,\dots,f_m\in L^2(\mathcal{M}(k),\lambda_{\mathcal{M}(k)})$,
\begin{eqnarray*}
\lefteqn{\mathbf{P}^{\lambda_{\mathcal{M}(k)}^n}_{\bar{x}}\left(f_1\left(\tilde{X}^{\lambda_{\mathcal{M}(k)}^n}_{t_1}\right)\dots f_m\left(\tilde{X}^{\lambda_{\mathcal{M}(k)}^n}_{t_m}\right)\right)}\\
&\rightarrow & \mathbf{P}^{\lambda_{\mathcal{M}(k)}}_{\bar{x}}\left(f_1\left({X}^{\lambda_{\mathcal{M}(k)}}_{t_1}\right)\dots f_m\left(X^{\lambda_{\mathcal{M}(k)}}_{t_m}\right)\right),
\end{eqnarray*}
as functions of $\bar{x}$ in $L^2(\mathcal{M}(k),\lambda_{\mathcal{M}(k)})$. Since the right-hand side above is continuous in $\bar{x}$ and the left-hand side is equicontinuous in $\bar{x}$ as $n$ varies (this can be proved using (\ref{holder}), (\ref{ecompare}), Lemma \ref{compmet} and \cite{FOT}, Lemma 1.3.3), it follows that the same convergence holds pointwise. In conjunction with the tightness of
$\tilde{X}^{\lambda_{\mathcal{M}(k)}^n}$, which is a consequence of Lemmas \ref{tightness} and \ref{approx}, this implies the distributional convergence of the first coordinate of (\ref{twoc}).

By the separability of $C(\mathbb{R}_+,(\mathcal{M}(k),d_\mathcal{M}))$, it is possible to assume that we are considering versions of the processes $\tilde{X}^{\lambda_{\mathcal{M}(k)}^n}$ and ${X}^{\lambda_{\mathcal{M}(k)}}$, each starting from $\bar{\rho}$, built on a common probability space such that the convergence of $\tilde{X}^{\lambda_{\mathcal{M}(k)}^n}$ to ${X}^{\lambda_{\mathcal{M}(k)}}$ occurs almost-surely (using the Skorohod coupling of, for example, \cite{Kallenberg}, Theorem 4.30). Supposing that the jointly continuous local times are also defined on this space, we find that, almost-surely, for any continuous bounded function $f$ on $(\mathcal{M}(k),d_\mathcal{M})$ and $t\geq 0$,
\begin{eqnarray*}
\int_{\mathcal{M}(k)}\tilde{L}^{\mathcal{M}(k),n}_t(\bar{x})f(\bar{x})\lambda_{\mathcal{M}(k)}(d\bar{x})
&=&\int_0^t f(\tilde{X}^{\lambda_{\mathcal{M}(k)}^n}_s)ds\\
&\rightarrow&\int_0^t f({X}^{\lambda_{\mathcal{M}(k)}}_s)ds\\
&=&\int_{\mathcal{M}(k)}L^{\mathcal{M}(k)}_t(\bar{x})f(\bar{x})\lambda_{\mathcal{M}(k)}(d\bar{x}).
\end{eqnarray*}
Applying the tightness of $\tilde{L}^{\mathcal{M}(k),n}$ (readily deduced from Lemmas \ref{tightness} and \ref{approx}) and the almost-sure continuity of $L^{\mathcal{M}(k)}$, it follows that $\tilde{L}^{\mathcal{M}(k),n}$ converges in distribution to $L^{\mathcal{M}(k)}$ on the probability space under consideration simultaneously with the convergence of processes, thereby completing the proof.
\end{proof}

The second main assumption of this section is that we have a sequence $(\mu_{\mathcal{M}(k)}^n)_{n\geq 1}$ of finite Borel measures on $(\mathcal{M}(k),d_\mathcal{M})$ that converges weakly to $\mu_{\mathcal{M}(k)}$. Under this assumption, we are able to show that the continuous additive functionals  $\hat{A}^{\mathcal{M}(k),n}$ defined by, for $t\geq 0$,
\[\hat{A}^{\mathcal{M}(k),n}_t:=\int_{\mathcal{M}(k)}L^{\mathcal{M}(k),n}_t(\bar{x})\mu^n_{\mathcal{M}(k)}(d\bar{x}),\]
converge to $\hat{A}^{\mathcal{M}(k)}$, as defined at (\ref{hata}). Since the result is an easy corollary of Proposition \ref{convfin} (and the continuous mapping theorem), we state it without proof.

{\cor\label{twoproc} Let $k\geq J+1$. If the process $X^{\lambda_{\mathcal{M}(k)}^n}$ has law  $\mathbf{P}^{\lambda_{\mathcal{M}(k)}^n}_{\bar{\rho}}$ and the process $X^{\lambda_{\mathcal{M}(k)}}$ has law  $\mathbf{P}^{\lambda_{\mathcal{M}(k)}}_{\bar{\rho}}$, then
\[\hat{A}^{\mathcal{M}(k),n}\rightarrow\hat{A}^{\mathcal{M}(k)}\]
in distribution as $n\rightarrow\infty$ in the space $C(\mathbb{R}_+,\mathbb{R}_+)$, simultaneously with the convergence statements of Proposition \ref{convfin}.}

\section{Critical random graph scaling limit}\label{crg}

In this section, for the purposes of introducing notation, we describe the scaling limit of the largest connected component of the critical random graph, as constructed in \cite{ABG}. As noted in the introduction, the basic ingredient in the definition of the random metric space $\mathcal{M}$ is a tilted version of the continuum random tree, and we start by presenting the excursion framework for this.

We will denote by $\mathcal{W}$ the space of continuous excursions, or more precisely the set
\[\left\{f\in C(\mathbb{R}_+,\mathbb{R}_+):\exists \sigma_f\in[0,\infty)\mbox{ such that $f(t)>0$ if and only if $t\in(0,\sigma_f)$}\right\}.\]
Throughout, we will reserve the notation $e^{(\sigma)}$ to represent a Brownian excursion of length $\sigma>0$, and set $e:=e^{(1)}$. The usual Brownian scaling applies to excursions, so that $e^{(\sigma)}$ has the same distribution as $(\sqrt{\sigma}e(t/\sigma))_{t\geq 0}$. As in \cite{ABG}, define a tilted excursion of length $\sigma$, $\tilde{e}^{(\sigma)}$ say, to be a random variable taking values in $\mathcal{W}$ whose distribution is characterised by
\[\mathbf{P}\left(\tilde{e}^{(\sigma)}\in{W}\right)=\frac{\mathbf{E}\left(\mathbf{1}_{\{e^{(\sigma)}\in W\}}e^{\int_0^\infty e^{(\sigma)}(t)dt}\right)}{\mathbf{E}\left(e^{\int_0^\infty e^{(\sigma)}(t)dt}\right)},\]
for measurable $W\subseteq \mathcal{W}$, where the $\sigma$-algebra we consider on $\mathcal{W}$ is that induced by the supremum norm on $C(\mathbb{R}_+,\mathbb{R}_+)$.

Let us now briefly outline the well-known map from $\mathcal{W}$ to the space of real trees. For $f\in\mathcal{W}$, define a distance on the interval $[0,\sigma_f]$ by setting $d_f(s,t):=f(s)+f(t)-2\inf\{f(r):\:r\in[s\wedge t,s\vee t]\}$, and then use the equivalence $s\sim_f t$ if and only if $d_f(s,t)=0$, to define $\mathcal{T}_f:=[0,\sigma_f]/\sim_f$. Denoting the canonical projection (with respect to $\sim_f$) from $[0,\sigma_f]$ to $\mathcal{T}_f$ by $\hat{f}$, it is possible to check that
$d_{\mathcal{T}_f}(\hat{f}(s),\hat{f}(t)):=d_f(s,t)$ defines a metric on
$\mathcal{T}_f$, and also that with this metric $\mathcal{T}_f$ is a compact real tree (see \cite{LegallDuquesne}, Theorem 2.1). In this article, the root of the tree $\mathcal{T}_f$ will always be defined to be the equivalence class $\hat{f}(0)$. Finally, although we will not need to refer to it in the remainder of this section, let us remark that the natural Borel measure on $\mathcal{T}_f$ can be constructed by setting $\mu_f:=\lambda_{[0,\sigma_f]}\circ \hat{f}^{-1}$, where $\lambda_{[0,\sigma_f]}$ is the usual one-dimensional Lebesgue measure on $[0,\sigma_f]$. This measure has full support and total mass equal to $\sigma_f$.

\begin{figure}[t]
\begin{center}
\vspace{5pt}
\scalebox{0.45}{\includegraphics{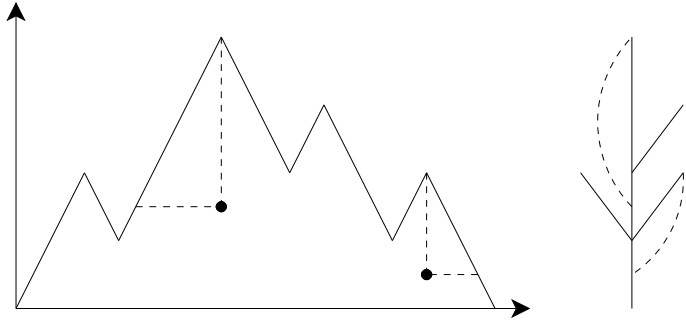}}
\vspace{-5pt}
\end{center}
\caption{Construction of $\mathcal{M}_{f,\mathcal{Q}}$. The left-hand side shows an example excursion $f$ and two point $\mathcal{Q}\cap f$. The right-hand side shows the associated $\mathcal{T}_f$ in solid lines. The space $\mathcal{M}_{f,\mathcal{Q}}$ is obtained from this by identifying the points at the end of the dotted lines.}\label{gluefig}
\vspace{-10pt}
\end{figure}

Given an excursion $f\in\mathcal{W}$ and a point-set $\mathcal{Q}\subseteq \mathbb{R}_+\times\mathbb{R}_+$ that only contains finitely many points in any compact set, the following procedure for defining a glued real tree is introduced in \cite{ABG}. First, define $\mathcal{Q}\cap f:=\left\{(t,x)\in\mathcal{Q}:0<x\leq f(t)\right\}$. For each point $(t,x)\in \mathcal{Q}\cap f$, let $u_{(t,x)}$ be the vertex $\hat{f}(t)\in\mathcal{T}_f$ and $v_{(t,x)}$ be the unique vertex on the path from the root $\hat{f}(0)$ to $u_{(t,x)}$ at a distance $x$ from the root. As at the beginning of Section \ref{sectree}, from the finite collection $E_\mathcal{Q}=\left\{(u_{(t,x)},v_{(t,x)}):(t,x)\in\mathcal{Q}\cap f\right\}$ of pairs of vertices of $\mathcal{T}_f$, we can define \[\mathcal{M}_{f,\mathcal{Q}}:=\mathcal{T}_f/\sim_{E_\mathcal{Q}},\]
which is a metric space when equipped with the quotient metric corresponding to $d_{\mathcal{T}_f}$, $d_{\mathcal{M}_{f,\mathcal{Q}}}$ say. (Figure \ref{gluefig} provides an  illustration of this construction.) The particular random point set of interest to us will be a Poisson process $\mathcal{P}$ on $\mathbb{R}_+\times \mathbb{R}_+$ of unit intensity with respect to Lebesgue measure, and we will write $\mathcal{M}^{(\sigma)}=(\mathcal{M}^{(\sigma)},d_{\mathcal{M}^{(\sigma)}})$ to be a random (non-empty) compact metric space with the distribution of
\[\left(\mathcal{M}_{\tilde{e}^{(\sigma)},\mathcal{P}},2d_{\mathcal{M}_{\tilde{e}^{(\sigma)},\mathcal{P}}}\right),\]
where $\tilde{e}^{(\sigma)}$ and $\mathcal{P}$ are assumed to be independent, which we can alternatively write as the quotient metric space
$(\mathcal{T}_{\tilde{e}^{(\sigma)}},2d_{\mathcal{T}_{\tilde{e}^{(\sigma)}}})/\sim_{E_\mathcal{P}}$. Note that formalisation of the above random variables is achieved by assuming the space of point sets is endowed with the topology induced by the usual Hausdorff convergence of non-empty compact subsets of $\mathbb{R}_+\times \mathbb{R}_+$, and the collection of non-empty compact metric spaces is endowed with the Gromov-Hausdorff topology, which will be generalised in the next section.

Using the above notions, it is possible to write the scaling limit of the largest connected component $\mathcal{C}_1^n$ as $\mathcal{M}^{(Z_1)}$, where the excursion length $Z_1$ is a random variable whose distribution we now define. Fix $\lambda\in\mathbb{R}$ to be the parameter in the description of the critical window (i.e. $p=n^{-1}+\lambda n^{-4/3}$), and let $B^\lambda=(B^\lambda_t)_{t\geq 0}$ be a Brownian motion with parabolic drift obtained by setting
\[B_t^\lambda:=B_t+\lambda t-\frac{t^2}{2},\]
where $B=(B_t)_{t\geq 0}$ is a standard Brownian motion started from the origin. Now, let $Z=(Z_{n})_{n\geq 1}$ be the lengths of the excursions of the reflected process $B^\lambda_t-\min_{s\in[0,t]}B^\lambda_s$ above 0 arranged in descending order, so that $Z_1$ is the length of the longest such excursion. (In fact, the reflected process is also related to the previous discussion, since each of its excursions, when conditioned to have length $\sigma$, is distributed as $\tilde{e}^{(\sigma)}$).

Although we will not apply it in its following form, for completeness, it seems appropriate to put all the above pieces together and state a (simplified version of) the main result of \cite{ABG}, which demonstrates the convergence of the rescaled $\mathcal{C}_1^n$. So as to make $\mathcal{C}_1^n$ into a metric space, it is assumed to be equipped with the usual shortest path graph distance $d_{\mathcal{C}_1^n}$. Note that the convergence of the first coordinate was originally proved as part of \cite{Aldouscrg}, Corollary 2.

{\thm [\cite{ABG}, Theorem 24] Let $\mathcal{C}_1^n$ be the largest connected component of the random graph $G(n,p)$ and $Z^n_1$ be the number of vertices of $\mathcal{C}_1^n$, where $p=n^{-1}+\lambda n^{-4/3}$, then
\[\left(n^{-2/3}Z^n_1,\left(\mathcal{C}_1^n,n^{-1/3}d_{\mathcal{C}_1^n}\right)\right)\rightarrow \left(Z_1,\mathcal{M}\right),\]
in distribution, where $\mathcal{M}=(\mathcal{M},d_\mathcal{M})$ is a random compact metric space such that, conditional on $Z_1$, $\mathcal{M}\buildrel{d}\over = \mathcal{M}^{(Z_1)}$.}
\bigskip

Finally, we observe that the absolute continuity of the law of $\tilde{e}^{(\sigma)}$ with respect to $e^{(\sigma)}$ and Brownian scaling easily imply that the canonical measure $\mu_\mathcal{T}:=\mu_{\tilde{e}^{(Z_1)}}$ on $\mathcal{T}=(\mathcal{T},d_\mathcal{T}):=(\mathcal{T}_{\tilde{e}^{(Z_1)}},2d_{\mathcal{T}_{\tilde{e}^{(Z_1)}}})$ satisfies up to constants the same $\mathbf{P}$-a.s. asymptotic results as were proved for the volume measure on the continuum random tree in \cite{Croydoncrt}. In particular, the condition at (\ref{lowervol}) is satisfied, $\mathbf{P}$-a.s. Consequently, if we write $\mu_\mathcal{M}:=\mu_\mathcal{T}\circ\phi^{-1}$, where $\phi$ is the canonical projection from $\mathcal{T}$ to $\mathcal{M}$, then, for $\mathbf{P}$-a.e. realisation of $(\tilde{e}^{(Z_1)},\mathcal{P})$, it is possible to define the Dirichlet form  $(\mathcal{E}_\mathcal{M},\mathcal{F}_\mathcal{M})$ on $L^2(\mathcal{M},\mu_\mathcal{M})$ and the associated diffusion law $\mathbf{P}_{\bar{\rho}}^\mathcal{M}$ precisely as was done in Section \ref{sectree}.

\section{Continuous paths on compact length spaces}\label{contsec}

In this section, we introduce the generalised Gromov-Hausdorff topology for continuous paths on compact length spaces in which our main result will be proved. Moreover, with respect to this topology, we will show that the construction of $\mathbf{P}_{\overline{\rho}}^\mathcal{M}$ from $(\tilde{e}^{(Z_1)},\mathcal{P})$ is measurable. The objects under consideration here will be of the form $\mathcal{K}=(K,d_K,X^K)$, where $(K,d_K)$ is a non-empty compact length space and $X^K$ is a path in $C([0,1],K)$ (for the definition of a length space, see \cite{BBI}, Definition 2.1.6). We will denote by $\mathbb{K}$ the set of path-preserving isometry classes of such triples, where by a path-preserving isometry between $\mathcal{K}$ and $\mathcal{K}'$, we mean an isometry $\psi:(K,d_K)\rightarrow (K',d_{K'})$ that satisfies $X^{K'}=\psi\circ X^K$. Define a distance $d_\mathbb{K}$ on $\mathbb{K}$ by setting
\begin{eqnarray*}
\lefteqn{d_\mathbb{K}(\mathcal{K},\mathcal{K}'):=}\\
&&\inf_{(M,d_M),\varphi,\varphi'}\left\{d^{(M,d_M)}_H(\varphi(K),\varphi'(K'))\vee\sup_{t\in[0,1]}d_M(\varphi(X^K_t),\varphi'(X^{K'}_t))\right\},
\end{eqnarray*}
where the infimum is taken over all choices of metric space $(M,d_M)$ and isometric embeddings $\varphi:(K,d_K)\rightarrow(M,d_M)$, $\varphi':(K',d_{K'})\rightarrow(M,d_M)$, and $d_H^{(M,d_M)}$ is the usual Hausdorff distance between compact subsets of $M$. To check that this defines a metric on $\mathbb{K}$ such that $(\mathbb{K},d_\mathbb{K})$ is separable is a simple extension of the corresponding result for metric spaces without paths, and so we will only sketch its proof. Note that $d_\mathbb{K}$ actually provides a metric on all (isometry classes of) triples of the form $\mathcal{K}=(K,d_K,X^K)$, where $(K,d_K)$ is a compact metric space and $X^K$ is a path in $C([0,1],K)$. The restriction to length spaces allows us to apply a graph approximation result that is useful in proving separability.

{\lem $(\mathbb{K},d_\mathbb{K})$ is a separable metric space.}
\begin{proof} That $d_\mathbb{K}$ is non-negative, symmetric and finite is easy to check. To prove that it satisfies the triangle inequality, an alternative characterisation is useful. First, for two metric spaces $(K,d_K),(K',d_{K'})$ define a correspondence $\mathfrak{C}$ between them to be a subset of $K\times K'$ such that: for every $x\in K$ there exists at least one $x'\in K'$ such that $(x,x')\in\mathfrak{C}$, and similarly for every $x'\in K'$ there exists at least one $x\in K$ such that $(x,x')\in\mathfrak{C}$. It is then the case (cf. \cite{BBI}, Theorem 7.3.25) that, for $\mathcal{K},\mathcal{K}'\in\mathbb{K}$,
\begin{equation}\label{dkexp}
d_{\mathbb{K}}(\mathcal{K},\mathcal{K}')=\frac{1}{2}\inf_{\mathfrak{C}:(X^K,X^{K'})\in\mathfrak{C}}{\rm dis} \mathfrak{C},
\end{equation}
where the infimum is taken over all correspondences $\mathfrak{C}$ between $(K,d_K)$ and $(K',d_{K'})$ such that $(X^K_t,X^{K'}_t)\in\mathfrak{C}$ for every $t\in[0,1]$, and ${\rm dis} \mathfrak{C}$ is the distortion of $\mathfrak{C}$, as defined by ${\rm dis} \mathfrak{C}:=\sup\{|d_K(x,y)-d_{K'}(x',y')|:(x,x'),(y,y')\in \mathfrak{C}\}$. Given the expression for $d_{\mathbb{K}}$ at (\ref{dkexp}), it is possible to check that $d_\mathbb{K}$ satisfies the triangle inequality by making the obvious adaptations to \cite{BBI}, Exercise 7.3.26. A second consequence of (\ref{dkexp}) is that if $d_{\mathbb{K}}(\mathcal{K},\mathcal{K}')<\varepsilon$, then there exists a $2\varepsilon$-isometry $f_\varepsilon$ from $(K,d_K)$ to $(K',d_{K'})$ such that $d_{K'}(f_\varepsilon(X^K_t),X^{K'}_t)<2\varepsilon$ for $t\in[0,1]$ (cf. \cite{BBI}, Corollary 7.3.28). Applying this fact, we can repeat the proof of \cite{BBI}, Theorem 7.3.30 to confirm that $d_{\mathbb{K}}$ is positive definite, choosing the countably dense set considered there to include $\{X^K_q:q\in [0,1]\cap\mathbb{Q}\}$. Thus $d_\mathbb{K}$ is a metric, as required.

For separability, we start by noting that if $(K,d_K)$ is a compact length space, then it can be approximated arbitrarily well by finite graphs. In particular, for every $n\geq 1$, by compactness we can choose a finite $n^{-1}$-net $\tilde{K}_n\subseteq K$. Make this into a graph by connecting points $x,y\in\tilde{K}_n$ by an edge of length $d_{K}(x,y)$ if and only if they satisfy $d_{K}(x,y)<\varepsilon_n$, where $\varepsilon_n:=(8n^{-1} {\rm diam}(K,d_K) )^{1/2}\vee 3n^{-1}$. It is then the case that
\[d_{K}(x,y)\leq d_{\tilde{K}_n}(x,y)\leq d_{K}(x,y)+\varepsilon_n,\hspace{20pt}\forall x,y\in\tilde{K}_n,\]
where $d_{\tilde{K}_n}$ is the shortest path graph distance on $\tilde{K}_n$ (see proof of \cite{BBI}, Proposition 7.5.5). We extend the space $(\tilde{K}_n,d_{\tilde{K}_n})$ into a compact length space $(K_n,d_{K_n})$ by including line segments along edges with lengths equal to the lengths of the edges in the graph. Now, if $X^K\in C([0,1],K)$, then there exists a $\delta_n\in(0,1)\cap \mathbb{Q}$ such that
\[\sup_{\substack{s,t\in[0,1]:\\|s-t|\leq\delta_n}}d_K(X^K_s,X^K_t)<n^{-1}.\]
For $k=0,1,\dots,\lfloor \delta_n^{-1}\rfloor$, choose $X^{K_n}_{k\delta_n}$ to be a vertex in $\tilde{K}_n$ such that \[d_K(X^{K_n}_{k\delta_n},X^K_{k\delta_n})<n^{-1}.\]
Observe that
\[d_K(X^{K_n}_{k\delta_n},X^{K_n}_{(k+1)\delta_n})<2n^{-1}+d_K(X^{K}_{k\delta_n},X^{K}_{(k+1)\delta_n})<3n^{-1}.\] Hence, $X^{K_n}_{k\delta_n}$ and $X^{K_n}_{(k+1)\delta_n}$ are connected by a graph edge, and so we can extend the definition of $X^{K_n}$ to a path in $C([0,1],K_n)$ by linearly interpolating along the relevant line segments in $K_n$ (after time $\lfloor \delta_n^{-1}\rfloor \delta_n$, set $X^{K_n}$ to be constantly equal to $X^{K_n}_{\lfloor \delta_n^{-1}\rfloor \delta_n}$). With this choice of $(K_n,d_{K_n},X^{K_n})$, we have that $d_{\mathbb{K}}((K,d_K,X^K),(K_n,d_{K_n},X^{K_n}))<2(n^{-1}+\varepsilon_n)$. Since $\tilde{K}_n$ is a finite set, there is no difficulty in perturbing the metric $d_{{K}_n}$ so that it takes rational values between any two points of $\tilde{K}_n$ and is linear along edges, but the triple $(K_n,d_{K_n},X^{K_n})$ still satisfies the same bound. Consequently, we have shown that for any $\mathcal{K}\in\mathbb{K}$ we can find a sequence $(\mathcal{K}_n)_{n\geq 1}$ drawn from a countable subset of $\mathbb{K}$ such that $d_\mathbb{K}(\mathcal{K},\mathcal{{K}}_{n})\rightarrow 0$.
\end{proof}

The relevance of the space $(\mathbb{K},d_\mathbb{K})$ to our setting depends on the following lemma.

{\lem \label{lengthspace} In the setting of Section \ref{sectree}, $(\mathcal{M},d_\mathcal{M})$ is a length space, as is the metric space $(\mathcal{M}(k),d_\mathcal{M})$ for every $k\geq {J+1}$.}
\begin{proof} By \cite{BBI}, Corollary 2.4.17, to prove that $(\mathcal{M},d_\mathcal{M})$ is a length space, it will suffice to show that for every $\bar{x},\bar{y}\in\mathcal{M}$, $\varepsilon>0$, there exists a finite sequence $\bar{x}=\bar{x}_0,\bar{x}_1,\dots,\bar{x}_k=\bar{y}$ such that
$\sum_{i=1}^k d_\mathcal{M}(\bar{x}_{i-1},\bar{x}_i)\leq d_\mathcal{M}(\bar{x},\bar{y})+\varepsilon$ and also $d_\mathcal{M}(\bar{x}_{i-1},\bar{x}_i)\leq \varepsilon$ for each $i=1,\dots,k$. To prove this, fix $\bar{x},\bar{y}\in\mathcal{M}$, $\varepsilon>0$. By the definition of $d_\mathcal{M}$, there exist vertices $x_i,y_i\in\mathcal{T}$, $i=1,\dots,k$, such that $\bar{x}_1=\bar{x}$, $\bar{y}_i=\bar{x}_{i+1}$, $\bar{y}_k=\bar{y}$, and also $\sum_{i=1}^kd_\mathcal{T}(x_i,y_i)\leq d_\mathcal{M}(\bar{x},\bar{y})+\varepsilon$. Define $t_0=0$ and $t_i=\sum_{j=1}^id_\mathcal{T}(x_j,y_j)$, $i=1,\dots,k$, and let $\gamma:[t_{i-1},t_{i}]\rightarrow\mathcal{T}$ be the path of unit speed from $x_i$ to $y_i$ in $\mathcal{T}$. This map is not well-defined on $[0,t_k]$ in general, since at the times $t_i$ it might be defined multiply. However, since $\bar{y}_i=\bar{x}_{i+1}$, this is not a problem when its image under $\phi$ is considered. In particular, the map $\phi\circ\gamma:[0,t_k]\rightarrow\mathcal{M}$ is well-defined, and by construction is easily checked to satisfy
\begin{equation}\label{frod}
d_\mathcal{M}(\phi\circ \gamma(s),\phi\circ \gamma(t))\leq |s-t|,
\end{equation}
for every $s,t\in[0,t_k]$. Let $n:=\lceil t_k \varepsilon^{-1}\rceil$ and set $\bar{z}_i:=\phi\circ\gamma(it_k/n)$, $i=0,1,\dots, n$, then $\sum_{i=1}^n d_\mathcal{M}(\bar{z}_{i-1},\bar{z}_i)\leq t_k \leq d_\mathcal{M}(\bar{x},\bar{y})+\varepsilon$, where the first inequality is an application of (\ref{frod}), and the second follows from the definition of $t_k$. Moreover, we also have that $d_\mathcal{M}(\bar{z}_{i-1},\bar{z}_i)\leq t_k/n\leq \varepsilon$ for each $i=1,\dots,n$. This completes the proof that $(\mathcal{M},d_\mathcal{M})$ is a length space, and the proof for $(\mathcal{M}(k),d_\mathcal{M})$ is identical.
\end{proof}

As remarked at the end of the previous section, for $\mathbf{P}$-a.e. realisation of $(\tilde{e}^{(Z_1)},\mathcal{P})$ we can construct what we will call the quenched law of the Brownian motion on $\mathcal{M}$ started from $\bar{\rho}$, $\mathbf{P}^\mathcal{M}_{\bar{\rho}}$, as a probability measure on $C([0,1], \mathcal{M})$. Clearly, by the above result, we can also consider this as a probability measure on $\mathbb{K}$, and we are able to deduce the following measurability result for it.

{\propn \label{meas} In the setting of Section \ref{crg}, with respect to the weak convergence of measures on $\mathbb{K}$,
$\mathbf{P}^\mathcal{M}_{\bar{\rho}}$ is an $(\tilde{e}^{(Z_1)},\mathcal{P})$-measurable random variable.}
\begin{proof} We will follow an approximation argument similar to that applied in \cite{Croydoncbp}, Lemma 8.1. To begin with, let $\xi=(\xi_i)_{i=1}^\infty$ be a sequence of independent $U(0,1)$ random variables, which is assumed to be independent of $(\tilde{e}^{(Z_1)},\mathcal{P})$. For $\mathbf{P}$-a.e. realisation of $(\tilde{e}^{(Z_1)},\mathcal{P},\xi)$, we can well-define sequences of vertices $(u_i)_{i=1}^\infty$ and $(v_i)_{i=1}^J$ of $\mathcal{T}=\mathcal{T}_{\tilde{e}^{(Z_1)}}$, where $J:= \#\mathcal{P}\cap \tilde{e}^{(Z_1)}$, as follows. For $i\leq J$, let $(u_i,v_i)$ be equal to $(u_{(t_i,x_i)},v_{(t_i,x_i)})$, where $(t_i,x_i)$ is the point of $\mathcal{P}\cap \tilde{e}^{(Z_1)}$ for which $t_i$ is the $i$th smallest element of $\{s:(s,y)\in\mathcal{P}\cap \tilde{e}^{(Z_1)}\mbox{ for some }y\geq0\}$. For $i\geq J+1$, set $u_i:=\hat{\tilde{e}}^{(Z_1)}(Z_1\xi_{i-J})$.

Now, suppose $\Gamma$ is the collection of realisations of $(\tilde{e}^{(Z_1)},\mathcal{P},\xi)$ such that: $\xi$ is dense in $[0,1]$; $J$ is finite; the ordered sequence $\{(t_i,x_i)\}_{i=1}^J$ of elements of $\mathcal{P}\cap \tilde{e}^{(Z_1)}$ is well-defined (i.e. $t_1<t_2<\dots<t_J$); the canonical measure $\mu_\mathcal{T}$ on $\mathcal{T}$ is non-atomic and satisfies the lower bound at (\ref{lowervol}); $\mathcal{T}\backslash\{x\}$ consists of no more than three connected components for any $x\in\mathcal{T}$; $(u_i)_{i=1}^\infty$ are distinct leaves of $\mathcal{T}$ (i.e. $\mathcal{T}\backslash \{u_i\}$ is connected for any $i$), not equal to $\rho$; $(v_i)_{i=1}^J$ are distinct and $v_i\in[[\rho,u_i]]\backslash(\{b^\mathcal{T}(\rho,u_i,u_j):i,j\geq 1\}\cup\{\rho\})$ for every $i$. That $\mathbf{P}((\tilde{e}^{(Z_1)},\mathcal{P},\xi)\in\Gamma)=1$ can be confirmed by applying: elementary properties of uniform random variables and Poisson processes; the volume bounds of \cite{Croydoncrt}, as discussed in the previous section; \cite{LegallDuquesne}, Theorem 4.6; the fact that $\mu^\mathcal{T}$ is non-atomic and supported on the leaves of $\mathcal{T}$, $\mathbf{P}$-a.s. (\cite{Aldous3} or \cite{LegallDuquesne}, Theorem 4.6); and a straightforward argument using the fact that, conditional on $\mathcal{T}$ and $u_i$, the distribution of $v_i$ on $[[\rho,u_i]]$ is the (normalised) one-dimensional Hausdorff measure on this path.

Assume that we have a sequence of realisations $(\tilde{e}^{(Z_{1n})}_n,\mathcal{P}_n,\xi^{n}=(\xi^n_i)_{i=1}^\infty)\in\Gamma$ such that
\begin{eqnarray*}
\lefteqn{\left(\tilde{e}^{(Z_{1n})}_n,Z_{1n},\mathcal{P}_n\cap \tilde{e}^{(Z_{1n})}_n ,J_n:=\#\mathcal{P}_n\cap\tilde{e}^{(Z_{1n})}_n,\xi^{n}\right)}\\
&\rightarrow& \left(\tilde{e}^{(Z_1)},Z_1,\mathcal{P}\cap\tilde{e}^{(Z_1)},J,\xi\right),\hspace{100pt}
\end{eqnarray*}
for some $(\tilde{e}^{(Z_1)},\mathcal{P},\xi)\in\Gamma$, where we recall the topology we are considering for the convergence of point sets is the usual Hausdorff convergence of non-empty compact subsets of $\mathbb{R}_+\times\mathbb{R}_+$. Since the integers $J_n\rightarrow J$, we must have that $J_n=J$ for large $n$. Fix $k\geq J+1$ and define $\mathcal{T}(k)$, as in Section \ref{sectree}, to be the subtree $(\cup_{i=1}^k [[\rho,u_i]])\cup (\cup_{i=1}^J [[\rho,v_i]])$, which is equal to simply $\cup_{i=1}^k [[\rho,u_i]]$ under our assumptions. Define $\mathcal{T}_n(k)$ similarly from the objects indexed by $n$. A simple adaptation of \cite{Croydoncbp}, Lemma 4.1 (cf. the proof of \cite{Aldous3}, Theorem 20), allows it to be deduced that, for large $n$, there exists a homeomorphism $\Upsilon_{n,k}$ from $\mathcal{T}_n(k)$ to $\mathcal{T}(k)$ such that $\Upsilon_{n,k}(\rho_n)=\rho$, $\Upsilon_{n,k}(u_{in})=u_i$ for $i=1,\dots,k$, $\Upsilon_{n,k}(v_{in})=v_i$ for $i=1,\dots,J$, and also if $d_\mathcal{T}^n$ is a metric on $\mathcal{T}(k)$ defined by
\[d^n_\mathcal{T}(x,y):=d_{\mathcal{T}_n}(\Upsilon^{-1}_{n,k}(x),\Upsilon^{-1}_{n,k}(y)),\]
for $x,y\in\mathcal{T}(k)$, then the condition at (\ref{compare}) is satisfied (at least once $n$ is large enough). Note that, it is for this argument that the condition on the number of components of $\mathcal{T}\backslash \{x\}$, $x\in\mathcal{T}$, is required, as if it did not hold then it would not necessarily be the case that $\mathcal{T}_n(k)$ was homeomorphic to $\mathcal{T}(k)$ for large $n$. Moreover, also by suitably modifying \cite{Croydoncbp}, Lemma 4.1, we can assume that if $\mu^n_{\mathcal{T}(k)}:=\mu_{\mathcal{T}_n(k)}\circ\Upsilon^{-1}_{n,k}$, where $\mu_{\mathcal{T}_n(k)}$ is the projection onto $\mathcal{T}_n(k)$ of the natural measure on $\mathcal{T}_n$, then $\mu^n_{\mathcal{T}(k)}\rightarrow\mu_{\mathcal{T}(k)}$ weakly as Borel measures on $\mathcal{T}(k)$ (recall the definition of $\mu_{\mathcal{T}(k)}$ from below (\ref{mumkdef})). Letting $\phi$ be the canonical projection from $\mathcal{T}$ to $\mathcal{M}$, it follows that $\mu^n_{\mathcal{M}(k)}:=\mu^n_{\mathcal{T}(k)}\circ\phi^{-1}\rightarrow\mu_{\mathcal{M}(k)}$ weakly as Borel measures on $\mathcal{M}(k)$. As a consequence of Lemmas \ref{timechange}, \ref{akprops} and Corollary \ref{twoproc}, we therefore have that
\begin{equation}\label{plim}
\mathbf{P}^{\mu^n_{\mathcal{M}(k)}}_{\bar{\rho}}\rightarrow \mathbf{P}^{\mu_{\mathcal{M}(k)}}_{\bar{\rho}}
\end{equation}
weakly as probability measures on $C([0,1],\mathcal{M}(k))$, where the left hand-side is the law of the process $X^{\mu^n_{\mathcal{M}(k)}}$ associated with $(\frac{1}{2}\mathcal{E}^n_{\mathcal{M}(k)},\mathcal{F}_{\mathcal{M}(k)})$ considered as a Dirichlet form on $L^2(\mathcal{M}(k),\mu^n_{\mathcal{M}(k)})$, started from $\bar{\rho}$ (see above the inequality at (\ref{ecompare}) for a definition of the resistance form), and right-hand side is defined as in Section \ref{sectree}.

Now, the properties of $\Upsilon_{n,k}$ listed above readily allow it to be deduced that this map from $\mathcal{T}_n(k)$ to $\mathcal{T}(k)$ induces a homeomorphism  $\Upsilon_{n,k}^\mathcal{M}:\mathcal{M}_n(k)\rightarrow \mathcal{M}(k)$ such that, for $x\in\mathcal{T}_n(k)$, \begin{equation}\label{induce}
\Upsilon_{n,k}^\mathcal{M}(\phi_n(x))=\phi(\Upsilon_{n,k}(x)),
\end{equation}
where $\phi_n$ is the canonical projection from $\mathcal{T}_n$ to $\mathcal{M}_n$. We claim that the left-hand side of (\ref{plim}) is equal to $\mathbf{P}^{\mu_{\mathcal{M}_n(k)}}_{\bar{\rho}}\circ({\Upsilon_{n,k}^\mathcal{M}})^{-1}$, i.e. the law of $\Upsilon_{n,k}^\mathcal{M}(X^{\mu_{\mathcal{M}_n(k)}})$ started from $\bar{\rho}$, where $X^{\mu_{\mathcal{M}_n(k)}}$ is the $\mathcal{M}_n(k)$-valued process defined analogously to $X^{\mu_{\mathcal{M}(k)}}$ using the objects indexed by $n$. To prove this, first observe that, since $(\mathcal{T}_n(k),d_{\mathcal{T}_n})$ is isometric to $(\mathcal{T}(k),d_{\mathcal{T}}^n)$, then the resistance form associated with $(\mathcal{T}(k),d_\mathcal{T}^n)$, $(\mathcal{E}_{\mathcal{T}(k)}^n, \mathcal{F}_{\mathcal{T}(k)})$, satisfies \[\mathcal{E}_{\mathcal{T}(k)}^n(f,f)=\mathcal{E}_{\mathcal{T}_n(k)}(f\circ\Upsilon_{n,k},f\circ\Upsilon_{n,k}),\]
for $f\in \mathcal{F}_{\mathcal{T}(k)}=\{f\circ\Upsilon_{n,k}^{-1}:f\in\mathcal{F}_{\mathcal{T}_n(k)}\}$, where $(\mathcal{E}_{\mathcal{T}_n(k)}, \mathcal{F}_{\mathcal{T}_n(k)})$ is the resistance form associated with $(\mathcal{T}_n(k),d_{\mathcal{T}_n})$. Hence, applying (\ref{newchar}) and (\ref{induce}), the corresponding resistance forms on $\mathcal{M}(k)$ and $\mathcal{M}_n(k)$ are related via \[\mathcal{E}_{\mathcal{M}(k)}^n(f,f)=\mathcal{E}_{\mathcal{M}_n(k)}(f\circ\Upsilon^\mathcal{M}_{n,k},f\circ\Upsilon^\mathcal{M}_{n,k}).\] This establishes the claim, and in conjunction with (\ref{plim}) demonstrates that
\[\mathbf{P}^{\mu_{\mathcal{M}_n(k)}}_{\bar{\rho}}\circ(\Upsilon_{n,k}^\mathcal{M})^{-1}\rightarrow \mathbf{P}^{\mu_{\mathcal{M}(k)}}_{\bar{\rho}}\]
weakly as probability measures on $C([0,1],\mathcal{M}(k))$. Taking into account Lemma \ref{lengthspace}, it is straightforward to obtain from this that $\mathbf{P}^{\mu_{\mathcal{M}_n(k)}}_{\bar{\rho}}\rightarrow \mathbf{P}^{\mu_{\mathcal{M}(k)}}_{\bar{\rho}}$ weakly as probability measures on $\mathbb{K}$.

For each $k\in\mathbb{N}$, define a probability measure on $\mathbb{K}$ by setting
\[Q^{(k)}:=\mathbf{P}^{\mu_{\mathcal{M}(k)}}_{\bar{\rho}}\mathbf{1}_{\{k\geq J+1\}}+\delta_{\mathcal{K}}\mathbf{1}_{\{k\leq J\}},\]
where $\delta_{\mathcal{K}}$ is a measure on $\mathbb{K}$ placing a unit mass on an arbitrary point $\mathcal{K}\in\mathbb{K}$. On the set of realisations of $(\tilde{e}^{(Z_1)},Z_1,\mathcal{P}\cap \tilde{e}^{(Z_1)},J,\xi)$ for which $(\tilde{e}^{(Z_1)},\mathcal{P},\xi)\in\Gamma$, by the conclusion of the previous paragraph, the map
\[(\tilde{e}^{(Z_1)},Z_1,\mathcal{P}\cap \tilde{e}^{(Z_1)},J,\xi)\mapsto Q^{(k)}\]
is continuous, and therefore measurable. Since it is the case that the quintuplet $(\tilde{e}^{(Z_1)},Z_1,\mathcal{P}\cap \tilde{e}^{(Z_1)},J,\xi)$ is $(\tilde{e}^{(Z_1)},\mathcal{P},\xi)$-measurable and $\mathbf{P}((\tilde{e}^{(Z_1)},\mathcal{P},\xi)\in\Gamma)=1$, it follows that $Q^{(k)}$ is $(\tilde{e}^{(Z_1)},\mathcal{P},\xi)$-measurable for each $k$. By Proposition \ref{kinfty}, we have that $Q^{(k)}\rightarrow \mathbf{P}^{\mathcal{M}}_{\bar{\rho}}$ on $\Gamma$, and so $\mathbf{P}^{\mathcal{M}}_{\bar{\rho}}$ is also $(\tilde{e}^{(Z_1)},\mathcal{P},\xi)$-measurable. The proof is completed on noting that integrating out the $\xi$ variable leaves the measure $\mathbf{P}^{\mathcal{M}}_{\bar{\rho}}$ unchanged.
\end{proof}

As an immediate consequence of this result, we can define the annealed law of the Brownian motion on $\mathcal{M}$ started from $\bar{\rho}$ by setting
\begin{equation}\label{ann}
\mathbb{P}^\mathcal{M}_{\bar{\rho}}(A):=\int \mathbf{P}^\mathcal{M}_{\bar{\rho}}(A)\mathbf{P}\left(d(\tilde{e}^{(Z_1)},\mathcal{P})\right),
\end{equation}
for measurable $A\subseteq\mathbb{K}$.

\section{Encoding $\mathcal{C}_1^n$ and subsets}\label{subsetsec}

We start this section by describing the construction of the largest connected component $\mathcal{C}_1^n$ from a random graph tree and a discrete point process, as presented in \cite{ABG}. We will then apply this to state our understanding of what properties are satisfied by a sequence of typical realisations of $\mathcal{C}_1^n$ (see Assumption \ref{as1} below). Finally, to complete the preparatory work for proving our precise versions of the convergence result at (\ref{rwcresult}), we define a collection of subsets $\mathcal{C}_1^n(k)\subseteq\mathcal{C}_1^n$ and prove a corresponding convergence result for a family of processes $X^{\mathcal{C}_1^n(k)}$, which take values in suitable modifications of $\mathcal{C}_1^n(k)$. Before we continue, however, note that for a graph $G=(V(G),E(G))$ we will often abuse notation by identifying $G$ and its vertex set $V(G)$. In particular, when we write $\#G$, we mean the number of vertices of the graph $G$. Similarly, $x\in G$ should be read as $x\in V(G)$. Moreover, for a graph $G$, we write $d_G$ to represent the usual shortest path metric on the vertices of $G$.

First, suppose that $T^n_1$ is a random ordered graph tree such that $\#T_1^n$ has the same distribution as $Z_1^n$ and, conditional on $\#T_1^n$,
\[\mathbf{P}(T_1^n=T)\propto (1-p)^{-a(T)},\]
where $T$ ranges over the set of ordered graph trees with $\#T_1^n$ vertices. Here, $a(T)$ is the number of edges ``permitted by the ordered depth-first search'' of $T$. More precisely, let $\mathcal{O}^T=(\mathcal{O}_m^T)_{m=0}^{\#T-1}$ be the ``stack'' process associated with the ordered depth-first search of $T$, i.e. for each $m$, $\mathcal{O}_m^T$ is the ordered subset of vertices of $T$ that have been seen but not yet explored by the depth-first search algorithm at time $m$ (see \cite{ABG}, Section 2, for details). Define the depth-first walk $D^T=(D^T_m)_{m=0}^{\#T-1}$ of $T$ by setting $D^T_m:=\#\mathcal{O}_m^T-1$, and then set
\[a(T):=\sum_{m=1}^{\#T-1}D^T_m,\]
which we observe represents the number of places that extra edges could be added to $T$ such that the ordered depth-first search of its vertices is preserved (\cite{ABG}, Lemma 7).

Secondly, let $\mathcal{Q}^n$ be a random subset of $\mathbb{N}\times \mathbb{N}$ in which each point is present independently with probability $p$, and write
\[\mathcal{Q}^n\cap D^n:=\left\{(m,j)\in\mathcal{Q}^n:m\leq \#T_1^n-1, j\leq D^n_m\right\},\]
where $D^n=(D^n_m)_{m=0}^{\#T_1^n-1}$ is the random depth-first walk of $T_1^n$. For each element $(m,j)\in \mathcal{Q}^n\cap D^n$, associate a pair of vertices $u_{(m,j)}^n,v_{(m,j)}^n\in T^n_1$ by the following. The vertex $u_{(m,j)}^n$ is that visited by the depth-first search of $T^n_1$ at time $m$. The vertex $v_{(m,j)}^n$ is that lying in position $\#\mathcal{O}^n_m-j+1$ of the corresponding random stack $\mathcal{O}_m^n$ at time $m$. Denote the collection of these pairs as $E^n:=\{\{u_{(m,j)}^n,v_{(m,j)}^n\}:(m,j)\in \mathcal{Q}^n\cap D^n\}$. By \cite{ABG}, Lemma 18, we then have that
\[\mathcal{C}_1^n=(V(T_1^n), E(T_1^n)\cup E^n)\]
as ordered graphs in distribution, where we write $V(T_1^n)$ and $E(T_1^n)$ to be the vertex and edges sets of $T_1^n$ respectively. In words, the largest connected  component $\mathcal{C}^n_1$ can be constructed by adding the random selection of edges $E^n$ to the random graph tree $T_1^n$. This result allows us to assume that $T_1^n$, $\mathcal{Q}^n$ and $\mathcal{C}_1^n$ are built on the same probability space in such a way that the above equality holds $\mathbf{P}$-a.s. In this case, we clearly have that $\#T_1^n=Z_1^n$, $\mathbf{P}$-a.s.

To formulate our quenched convergence assumption, we will appeal to the following theorem, which collects together several results proved in \cite{ABG} (see the proofs of Lemma 19 and Theorem 22 in particular). The function $H^n=(H^n_m)_{m=0}^{Z_1^n-1}$ is the height process of $T_1^n$, so that $H^n_m$ is the graph distance between the root (the first ordered vertex of $T^n_1$) and the vertex visited at time $m$ by the depth-first search algorithm. The function $C^n=(C^n_m)_{m=0}^{2(Z_1^n-1)}$ is the corresponding contour function, which is obtained by recording the distance from the root of a particle that traces the boundary of the tree $T^n_1$ in a clockwise fashion, starting from the root. We note that the contour function does not actually appear in \cite{ABG}, but the convergence result involving it that is stated below is a simple consequence of \cite{MM}, Theorem 2, where in a clash with the above terminology it is termed the depth-first walk. Finally, we introduce the notation $J_n:=\#\mathcal{Q}^n\cap D^n$.

{\thm[\cite{ABG}] \label{dhc} For the random objects defined above, we have that
\[\left(n^{-1/3}D^n_{\lfloor n^{2/3} t\rfloor},n^{-1/3}H^n_{\lfloor n^{2/3} t\rfloor},n^{-1/3}C^n_{\lfloor 2n^{2/3} t\rfloor}\right)_{t\geq 0}\rightarrow
\left(\tilde{e}^{(Z_1)},2\tilde{e}^{(Z_1)},2\tilde{e}^{(Z_1)}\right),\]
\[\{(n^{-2/3}m,n^{-1/3}j): (m,j)\in\mathcal{Q}^n\cap D^n\}\rightarrow \mathcal{P}\cap \tilde{e}^{(Z_1)},\]
\[n^{-2/3}Z_1^n\rightarrow Z_1,\hspace{20pt}J_n\rightarrow J,\]
simultaneously in distribution, where the processes $D^n$, $H^n$, $C^n$ are extended to all positive integers by setting them to be equal to 0 where they are not already defined, and $\tilde{e}^{(Z_1)}$, $Z_1$, $\mathcal{P}$, $J:=\#\mathcal{P}\cap \tilde{e}^{(Z_1)}$ are defined as in Section \ref{crg}. The first convergence statement is in the space $D(\mathbb{R}_+,\mathbb{R}_+)^3$. The second convergence statement is with respect to the Hausdorff convergence of compact sets.}
\bigskip

Based on this result, let us introduce an assumption we will henceforth commonly make for a sequence of realisations of the pairs $(T_1^n,\mathcal{Q}^n)$. The definition of $\Gamma$ should be recalled from the proof of Proposition \ref{meas}.

{\assu \label{as1} The deterministic sequence $\{(T_1^n,\mathcal{Q}^n)\}_{n\geq1}$ of ordered graph trees and discrete point sets satisfy the convergence statements of Theorem \ref{dhc} for some fixed realisation of $(\tilde{e}^{(Z_1)},\mathcal{P})$. Moreover, $(\tilde{e}^{(Z_1)},\mathcal{P},\xi)\in\Gamma$ for some sequence $\xi\in(0,1)^{\mathbb{N}}$.}
\bigskip

Note that, by Theorem \ref{dhc} and the definition of $\Gamma$, it is possible to construct versions of the random pairs $(T_1^n,\mathcal{Q}^n)$ introduced at the beginning of the section for which Assumption \ref{as1} holds, $\mathbf{P}$-a.s. Furthermore, whenever we refer to the quantities $\mathcal{T}$, $\mathcal{T}(k)$, $\mathcal{M}$, $\mathcal{M}(k)$, \dots, under Assumption \ref{as1}, we mean the quantities associated with $(\tilde{e}^{(Z_1)},\mathcal{P})$ as in Section \ref{crg} and the proof of Proposition \ref{meas}.

\begin{figure}[t]
\begin{center}
\vspace{0pt}
\scalebox{0.45}{\includegraphics{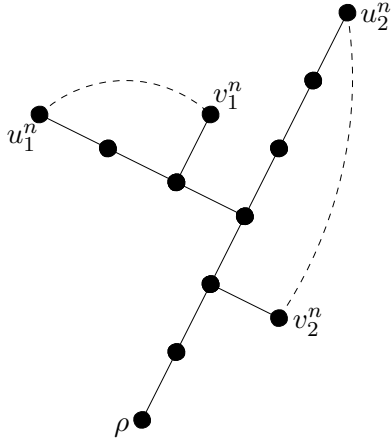}}
\rput(-4.5,3.85){$u_1^n$}
\rput(-1.8,4.4){$v_1^n$}
\rput(0.15,5.45){$u_2^n$}
\rput(-0.75,1.4){$v_2^n$}
\rput(-3.2,0){$\rho$}
\vspace{0pt}
\end{center}
\caption{A example realisation of $T_1^n(2)$ with $J_n=2$ (solid lines). The graph $\mathcal{C}_1^n(2)$ is obtained from this by adding the edges in $E^n$ (shown as dotted lines).}\label{t2nkfig}
\vspace{-15pt}
\end{figure}

Let us now suppose Assumption \ref{as1} holds and proceed to defining the subsets $\mathcal{C}_1^n(k)\subseteq \mathcal{C}_1^n$. For this purpose, it will be convenient to order the points of $\mathcal{Q}^n\cap D^n$, which will we do by supposing the elements of the sequence $\{(m_i,j_i)\}_{i=1}^{J_n}$ are the points of $\mathcal{Q}^n\cap D^n$ arranged in such a way that $m_1\leq m_2\leq \dots\leq m_{J_n}$. If there are multiple ways of doing this, we simply pick one arbitrarily. We will write $(u_i^n,v_i^n)=(u^n_{(m_i,j_i)},v^n_{(m_i,j_i)})$ for $i=1,\dots,J_n$. For $i\geq J_n+1$, let $u_i^n$ be the vertex of $T_1^n$ visited by the contour function at time $\lfloor 2 (Z_1^n-1) \xi_{i-J_n} \rfloor$. Analogously to the definitions of $\mathcal{T}(k)$ and $\mathcal{M}(k)$ at (\ref{tkdef}) and (\ref{mkdef}) respectively, we then set, for $k\geq J_n$,
\begin{equation}\label{ed}
T_1^n(k):=\left(\cup_{i=1}^k [[\rho,u^n_i]]\right)\cup\left(\cup_{i=1}^{J_n} [[\rho,v^n_i]]\right),
\end{equation}
where $\rho$ is the root of $T_1^n$, which is simply chosen to be the first ordered vertex of $T_1^n$, and $[[\rho,x]]$ is the unique injective path from $\rho$ to $x$ in $T_1^n$, and also
\begin{equation}\label{c1nk}
\mathcal{C}_1^n(k):=(V(T_1^n(k)), E(T_1^n(k))\cup E^n).
\end{equation}
(See Figure \ref{t2nkfig}.)

In the proof of the simple random walk convergence results of the next section, we will consider processes projected from $\mathcal{C}_1^n$ to $\mathcal{C}_1^n(k)$ by a map $\phi_{\mathcal{C}_1^n,\mathcal{C}_1^n(k)}$, which we define by setting, for $x$ a vertex in $\mathcal{C}_1^n$ (or $T_1^n$),
\begin{equation}\label{cproj}
\phi_{\mathcal{C}_1^n,\mathcal{C}_1^n(k)}(x)=\phi_{T_1^n,T_1^n(k)}(x),
\end{equation}
where $\phi_{T_1^n,T_1^n(k)}:T_1^n\rightarrow T_1^n(k)$ is defined as in the real tree case at (\ref{ttk}). Before continuing, let us collect together some first properties of the subsets $\mathcal{C}_1^n(k)$ and maps $\phi_{\mathcal{C}_1^n,\mathcal{C}_1^n(k)}$ in a lemma. In part (a) of the result, we consider the asymptotic properties of the quantities
\begin{equation}\label{deltankdef}
\Delta_n^{(k)}:=\sup_{x\in\mathcal{C}_1^n}d_{\mathcal{C}_1^n}\left(x,\phi_{\mathcal{C}_1^n,\mathcal{C}_1^n(k)}(x)\right),
\end{equation}
\begin{equation}\label{lambdankdef}
\Lambda_n^{(k)}:=\#E(\mathcal{C}_1^n(k)),
\end{equation}
under Assumption \ref{as1}. Part (b) is a simple structural result, the statement of which involves the projection of the uniform measure on $\mathcal{C}_1^n$ onto $\mathcal{C}_1^n(k)$. To be precise, let $\mu_{\mathcal{C}_1^n}$ be the measure placing mass 1 on each vertex of $\mathcal{C}_1^n$, so that it has total mass $Z_1^n$, and write
\begin{equation}\label{muc1nk}
\mu_{\mathcal{C}_1^n(k)}:=\mu_{\mathcal{C}_1^n}\circ \phi_{\mathcal{C}_1^n,\mathcal{C}_1^n(k)}^{-1}.
\end{equation}

{\lem \label{previous} Suppose Assumption \ref{as1} holds.\\
(a) The quantities $\Delta_n^{(k)}$, $\Lambda_n^{(k)}$, satisfy
\begin{equation}\label{a1}
\lim_{k\rightarrow\infty}\limsup_{n\rightarrow\infty}n^{-1/3}\Delta_n^{(k)}=0,
\end{equation}
and, for each $k\geq J_n+1$,
\begin{equation}\label{a2}
\lim_{n\rightarrow\infty}n^{-1/3}\Lambda_n^{(k)}=\Lambda^{(k)},
\end{equation}
where $\Lambda^{(k)}$ was defined at (\ref{Lambdak}).\\
(b) Let $k\geq J_n+1$, $x\in\mathcal{C}_1^n(k)$ and $E^x$ be those edges in $E(\mathcal{C}_1^n(k))$ that contain $x$. The component of $\mathcal{C}_1^n\backslash E^x$ containing $x$ is a graph tree on $\mu_{\mathcal{C}_1^n(k)}(\{x\})$ vertices, with $x$ being the only one of these in $\mathcal{C}_1^n(k)$.}
\begin{proof} First, fix $k\geq J_n+1$, and observe that
\begin{eqnarray}
\Delta^{(k)}_n&\leq&\sup_{x\in{T}_1^n}d_{{T}_1^n}\left(x,\phi_{{T}_1^n,{T}_1^n(k)}(x)\right)\nonumber\\
&\leq&\sup_{x\in{T}_1^n}d_{{T}_1^n}\left(x,\phi_{{T}_1^n,\tilde{{T}}_1^n(k)}(x)\right),\label{deltup}
\end{eqnarray}
where $\tilde{{T}}_1^n(k):=\cup_{i=1}^k[[\rho,u^n_i]]$ and the projection $\phi_{{T}_1^n,\tilde{{T}}_1^n(k)}:T_1^n\rightarrow\tilde{T}_1^n(k)$ is defined similarly to (\ref{ttk}). That the second inequality holds is a simple consequence of the fact that $\tilde{{T}}_1^n(k)\subseteq T_1^n(k)$.

Now, let $K^n(m):=2m-H^n_m$ for $m=0,1,\dots,Z_1^n-1$ and define a sequence of positive integers $\tilde{\xi}^n:=(\tilde{\xi}_i^n)_{i\geq 1}$ by setting
\[\tilde{\xi}_i^n:=\left\{\begin{array}{ll}
                          K^n(m_i), & {\mbox{for }}i\leq J_n,\\
                          \lfloor 2 (Z_1^n-1) \xi_{i-J_n} \rfloor, & \mbox{otherwise,}
                        \end{array}\right.
\]
so that, for every $i$, $u_i^n$ is the vertex of $T_1^n$ visited by the contour function at time $\tilde{\xi}_i^n$ (see \cite{LeGallDuquesne2}, Section 2.4, for the result when $i\leq J_n$). By Assumption \ref{as1}, we have that
\[\frac{2m_i^n}{n^{2/3}}\rightarrow 2t_i,\hspace{20pt}\frac{H^n_{m_i^n}}{n^{2/3}}\rightarrow 0,\]
where $t_i$ is defined as in the proof of Proposition \ref{meas}. Therefore $(2^{-1}n^{-2/3}\tilde{\xi}^n_i)_{i\geq 1}$ converges to the sequence $\tilde{\xi}:=(t_1,\dots,t_J,Z_1\xi_1,Z_1\xi_2,\dots)$. Applying this result, the denseness in $(0,1)$ of $\xi$, and the convergence of contour functions, by following a deterministic version of the proof of \cite{Aldous3}, Theorem 20 (which contains a distributional version of the same result), we find that
\[\lim_{k\rightarrow\infty}\limsup_{n\rightarrow\infty}n^{-1/3}\sup_{x\in{T}_1^n}d_{{T}_1^n}\left(x,\phi_{{T}_1^n,\tilde{{T}}_1^n(k)}(x)\right)=0.\]
Together with the upper bound at (\ref{deltup}), this completes the proof of (\ref{a1}). Another simple consequence of \cite{Aldous3}, Theorem 20, is that $n^{-1/3}\#E(\tilde{T}_1^n(k))$ converges to the `total edge length' or, more precisely, the one-dimensional Hausdorff measure of $\mathcal{T}(k)$, which we have denoted $\Lambda^{(k)}$, see (\ref{Lambdak}). Furthermore, as follows from an observation in Section 2 of \cite{ABG}, the vertex $v^n_i$ is always at a distance 1 from the path from $\rho$ to $u_i^n$ in $T_1^n$. Therefore, we must have that
\[\#E(\tilde{T}_1^n(k))\leq \#E(\mathcal{C}_1^n(k))\leq \#E({T}_1^n(k))+J_n\leq \#E(\tilde{T}_1^n(k))+2J_n.\]
Combining this with the convergence result for $\#E(\tilde{T}_1^n(k))$, the result at (\ref{a2}) follows.

For the proof of (b), fix $k\geq J_n+1$ and let $x\in\mathcal{C}_1^n(k)$. Since $\phi_{\mathcal{C}_1^n,\mathcal{C}_1^n(k)}(y)=y$ for any $y\in \mathcal{C}_1^n(k)$, it must be the case that $\phi_{\mathcal{C}_1^n,\mathcal{C}_1^n(k)}^{-1}(\{x\})\cap \mathcal{C}_1^n(k)=\{x\}$. From this and the obvious fact that the sets $\phi_{\mathcal{C}_1^n,\mathcal{C}_1^n(k)}^{-1}(\{y\})$, $y\in\mathcal{C}_1^n(k)$, are disjoint, it follows that the component of $\mathcal{C}_1^n\backslash E^x$ containing $x$ contains precisely the vertices $\phi_{\mathcal{C}_1^n,\mathcal{C}_1^n(k)}^{-1}(\{x\})$, of which there are $\mu_{\mathcal{C}_1^n(k)}(\{x\})$ and, of these, only $x$ is in $\mathcal{C}_1^n(k)$. Moreover, since the edges in $E^n$ only connect together vertices of $\mathcal{C}_1^n(k)$, it also follows that the component of interest is a graph tree.
\end{proof}

To describe the processes $X^{\mathcal{C}^n_1(k)}$, we start by giving an alternative construction of $\mathcal{C}_1^n(k)$. Let $\tilde{T}_1^n(k)$ be defined as in the proof of Lemma \ref{previous} and introduce a new graph $\hat{T}_1^n(k)$ with vertex set
\begin{equation}\label{vhat}
V(\hat{T}_1^n(k)):=V(\tilde{T}_1^n(k))\cup\left(\cup_{i=1}^{J_n}\{v_i^n,\tilde{w}_i^n\}\right),
\end{equation}
where $\tilde{w}_i^n$, $i=1,\dots,J_n$, are a collection of ``new'' vertices not already contained in $V({T}_1^n(k))$ (we are not concerned with the order of vertices here), and edge set
\begin{equation}\label{ehat}
E(\hat{T}_1^n(k)):=E(\tilde{T}_1^n(k))\cup E^n\cup \left\{\{v_i^n,\tilde{w}_i^n\}:i=1,\dots,J_n\}\right\}.
\end{equation}
We will demonstrate in the following proof that, under Assumption \ref{as1} and for large $n$, the graph $\hat{T}_1^n(k)$ is simply that obtained from $\tilde{T}_1^n(k)$ by connecting to the vertices $u_i^n$, which are leaves of $\tilde{T}_1^n(k)$, the disjoint length two line-segments $\{v_i^n,\tilde{w}_i^n\}$ with an edge from $u_i^n$ to $v_i^n$.

\begin{figure}[t]
\begin{center}
\vspace{0pt}
\scalebox{0.45}{\includegraphics{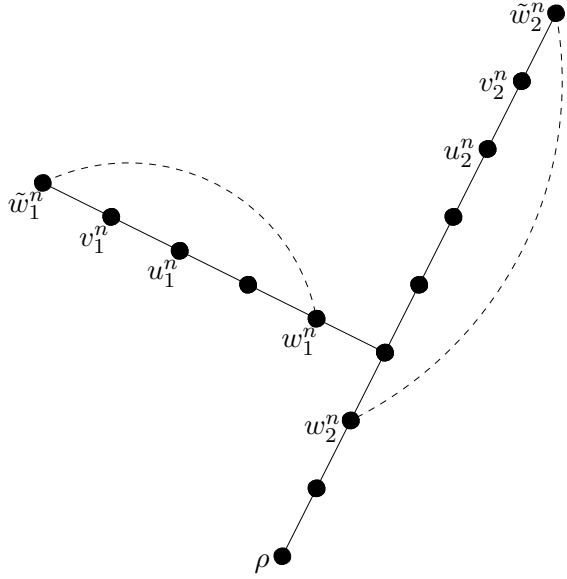}}
\rput(-5.4,3.85){$u_1^n$}
\rput(-6.3,4.3){$v_1^n$}
\rput(-7.2,4.75){$\tilde{w}_1^n$}
\rput(-3.6,2.95){$w_1^n$}
\rput(-1.5,5.45){$u_2^n$}
\rput(-1.05,6.35){$v_2^n$}
\rput(-.6, 7.25){$\tilde{w}_2^n$}
\rput(-3.3,1.8){$w_2^n$}
\rput(-4.1,0){$\rho$}
\vspace{5pt}
\end{center}
\caption{Realisation of $\hat{T}_1^n(2)$ corresponding to the $T_1^n(2)$ shown in Figure \ref{t2nkfig} (solid lines). The graph $\hat{T}_1^n(k)/ \sim_{n,k}$ is obtained from this by identifying the vertices joined by dotted lines.}\label{t2nkhatfig}
\vspace{-10pt}
\end{figure}

On the vertices of $\hat{T}_1^n(k)$, consider the vertex equivalence relation
\[x\sim_{n,k} y\hspace{20pt}\Leftrightarrow\hspace{20pt}x=y\mbox{ or }\{x,y\}=\{w_i^n,\tilde{w}_i^n\}\mbox{ for some }i=1,\dots,J_n,\]
where $w_i^n$ is the unique vertex in the path from $\rho$ to $u_i^n$ connected to $v_i^n$ by an edge in $E(T_1^n)$ (for the existence of such a vertex, see \cite{ABG}, Section 2). Let $\hat{T}_1^n(k)/ \sim_{n,k}$ be the graph obtained from $\hat{T}_1^n(k)$ by identifying vertices in equivalence classes and then replacing any resulting multiple edges with a single one. (See Figure \ref{t2nkhatfig}.) In fact, since $\tilde{w}_i^n$ is not connected in $\hat{T}_1^n(k)$ to any other vertex than $v_i^n$, it is possible to deduce that a graph identical to $\hat{T}_1^n(k)/ \sim_{n,k}$ can be constructed by simply deleting the vertices $\tilde{w}_i^n$ and edges $\{v_i^n,\tilde{w}_i^n\}$ from $\hat{T}_1^n(k)$, and then adding the edges $\{v_i^n,{w}_i^n\}$. More precisely, this new graph has vertex set
\[V(\hat{T}_1^n(k))\backslash\cup_{i=1}^{J_n}\{\tilde{w}_i^n\}=V(\tilde{T}_1^n(k))\cup\left(\cup_{i=1}^{J_n}\{v_i^n\}\right)=V(\mathcal{C}_1^n(k)),\]
and edge set
\begin{eqnarray*}
\lefteqn{E(\hat{T}_1^n(k))\cup\left\{\{v_i^n,{w}_i^n\}:i=1,\dots,J_n\right\}\backslash\left\{\{v_i^n,\tilde{w}_i^n\}:i=1,\dots,J_n\right\}}\\
&&=E({T}_1^n(k))\cup E^n = E(\mathcal{C}_1^n(k)).\hspace{160pt}
\end{eqnarray*}
Thus $\mathcal{C}_1^n(k)$ and $\hat{T}_1^n(k)/ \sim_{n,k}$ have exactly the same graph structure. Why this picture of $\mathcal{C}_1^n(k)$ is helpful is because of the following asymptotic description of $\hat{T}_1^n(k)$.

{\lem \label{l63} Suppose Assumption \ref{as1} holds and $k\geq J+1$. For large enough $n$, we have that $\hat{T}_1^n(k)$ is a graph tree with $\#E(\mathcal{C}_1^n(k))$ edges. Moreover, if $(\hat{T}_1^n(k),d_{\hat{T}^n_1(k)})$ is considered as a real tree by including unit line segments along edges, then there exists a homeomorphism \[\Upsilon_{\hat{T}_1^n(k),\mathcal{T}(k)}:\hat{T}_1^n(k)\rightarrow\mathcal{T}(k)\] such that $\Upsilon_{\hat{T}_1^n(k),\mathcal{T}(k)}(\rho)=\rho$, $\Upsilon_{\hat{T}_1^n(k),\mathcal{T}(k)}(\tilde{w}_i^n)=u_i$ and $\Upsilon_{\hat{T}_1^n(k),\mathcal{T}(k)}(w_i^n)=v_i$ for $i\leq J_n$, $\Upsilon_{\hat{T}_1^n(k),\mathcal{T}(k)}(u_i^n)=u_i$ for $i=J_n+1,\dots,k$, and
the sequence of metrics $d^n_\mathcal{T}$ on $\mathcal{T}(k)$ defined by
\[d^n_\mathcal{T}(x,y):=n^{-1/3}d_{\hat{T}^n_1(k)}\left(\Upsilon_{\hat{T}_1^n(k),\mathcal{T}(k)}^{-1}(x),\Upsilon_{\hat{T}_1^n(k),\mathcal{T}(k)}^{-1}(y)\right)\] satisfy (\ref{compare}). Finally, defining $\phi_{T_1^n,\tilde{T}_1^n(k)}$ as in the proof of Lemma \ref{previous}, it is possible to define $\Upsilon_{\hat{T}_1^n(k),\mathcal{T}(k)}$ so it is further the case that
\[n^{-2/3}\mu_{\mathcal{C}_1^n}\circ \phi_{T_1^n,\tilde{T}_1^n(k)}^{-1}\circ
\Upsilon_{\hat{T}_1^n(k),\mathcal{T}(k)}^{-1}\rightarrow \mu_{\mathcal{T}(k)}\]
weakly as measures on $\mathcal{T}(k)$. Note that the left-hand side above is well-defined because $\tilde{T}_1^n(k)\subseteq \hat{T}_1^n(k)$, and the right-hand side was defined below (\ref{mumkdef}).}
\begin{proof} Under Assumption \ref{as1}, by considering the convergence of contour functions and selection of vertices using $\tilde{\xi}^n$ and $\tilde{\xi}$ as in the proof of the previous result, it is possible to deduce that there exists a homeomorphism
\[\Upsilon_{\tilde{T}_1^n(k),\mathcal{T}(k)}:\tilde{T}_1^n(k)\rightarrow\mathcal{T}(k),\] where the domain here is the real tree version of $\tilde{T}_1^n(k)$ obtained by including unit line segments along edges, such that $\Upsilon_{\tilde{T}_1^n(k),\mathcal{T}(k)}(\rho)=\rho$, $\Upsilon_{\tilde{T}_1^n(k),\mathcal{T}(k)}({u}_i^n)=u_i$ for $i\leq k$, the sequence of metrics $\tilde{d}^n_\mathcal{T}$ on $\mathcal{T}(k)$ defined by
\[\tilde{d}^n_\mathcal{T}(x,y):=n^{-1/3}d_{\tilde{T}^n_1(k)}\left(\Upsilon_{\tilde{T}_1^n(k),\mathcal{T}(k)}^{-1}(x),\Upsilon_{\tilde{T}_1^n(k),\mathcal{T}(k)}^{-1}(y)\right)\] satisfy (\ref{compare}), and also
\[n^{-2/3}\mu_{\mathcal{C}_1^n}\circ \phi_{T_1^n,\tilde{T}_1^n(k)}^{-1}\circ
\Upsilon_{\tilde{T}_1^n(k),\mathcal{T}(k)}^{-1}\rightarrow \mu_{\mathcal{T}(k)}\]
weakly as Borel measures on $\mathcal{T}(k)$. To do this, one needs to make only a very minor modification to the conclusion of \cite{Croydoncbp}, Lemma 4.1, which is essentially the same result expressed in a slightly different way (cf. \cite{Aldous3}, Theorem 20).

Now, by applying ideas from the proof of \cite{ABG}, Theorem 22, it is possible to check that, for $i\leq J_n$,
\begin{eqnarray*}
\left|n^{-1/3}d_{\tilde{T}_1^n(k)}(\rho,w_i^n)-d_\mathcal{T}(\rho,v_i)\right|&\rightarrow&0,
\end{eqnarray*}
where $j_i$ was defined above (\ref{ed}) and $x_i$ is defined in the proof of Proposition \ref{meas}. Hence $\Upsilon_{\tilde{T}_1^n(k),\mathcal{T}(k)}(w_i^n)\rightarrow v_i$ in $\mathcal{T}(k)$ as $n\rightarrow\infty$. Since, by assumption, the vertices $v_i\not \in V:=\{b^\mathcal{T}(x,y,z):x,y,z\in\{\rho, u_1,\dots,u_J\}\}$ and are distinct, it must be the case that there exist disjoint neighbourhoods of the vertices $v_i$, each isometric to a line segment. It readily follows that, by suitably distorting $\Upsilon_{\tilde{T}_1^n(k),\mathcal{T}(k)}$ in neighbourhoods of each of the $w_i^n$, which can be chosen to be disjoint for large $n$, we can redefine the map $\Upsilon_{\tilde{T}_1^n(k),\mathcal{T}(k)}$ so that, in addition to the above properties, it is also the case that $\Upsilon_{\tilde{T}_1^n(k),\mathcal{T}(k)}(w_i^n)=v_i$.

In fact, since $\Upsilon_{\tilde{T}_1^n(k),\mathcal{T}(k)}$ is a homeomorphism satisfying $\Upsilon_{\tilde{T}_1^n(k),\mathcal{T}(k)}(w_i^n)=v_i$, it transpires that, for large $n$,
\[w_i^n\not\in V^n:=\{b^{{T}_1^n}(x,y,z):x,y,z\in\{\rho, u_1^n,\dots,u_{J_n}^n\}\}.\]
In particular, this implies that $w_i^n$ is only connected by an edge in $\tilde{T}_1^n(k)$ to the two adjacent vertices in the path from $\rho$ to $u_i^n$. As a consequence, $\{v_i^n,w_i^n\}$ is not an edge of $\tilde{T}_1^n(k)$, so $v_i^n\not\in\tilde{T}_1^n(k)$ and, moreover, since the $w_i^n$ are necessarily distinct for large $n$ (because the vertices $v_i$ are distinct), then so are the $v_i^n$.  Again applying the fact that $\Upsilon_{\tilde{T}_1^n(k),\mathcal{T}(k)}$ is a homeomorphism, we also have that $u_i^n$, $i=1,\dots,J_n$, are distinct leaves of $\tilde{T}_1^n(k)$. Piecing the above observations together, it is now straightforward to confirm the picture that $\hat{T}_1^n(k)$ is obtained from $\tilde{T}_1^n(k)$ by simply adding disjoint paths of length 2 to each of the leaves $u_i^n$. Clearly this results in a graph tree, and to define the required map $\Upsilon_{\hat{T}_1^n(k),\mathcal{T}(k)}$ on the associated real tree, we can simply ``stretch'' the domain of the map $\Upsilon_{\tilde{T}_1^n(k),\mathcal{T}(k)}$ near the relevant $J_n$ leaves.

Finally, the one remaining claim is that $\#E(\hat{T}_1^n(k))=\#E(\mathcal{C}_1^n(k))$. For large $n$, the description of $\hat{T}_1^n(k)$ from the preceding paragraph implies that $\#E(\hat{T}_1^n(k))=\#E(\tilde{T}_1^n(k))+2J_n$, and we will deduce that $\#E(\mathcal{C}_1^n(k))$ is equal to the same expression.  First, since we can assume that $\{v_i^n,w_i^n\}$ is not an edge of $\tilde{T}_1^n(k)$ and the vertices $w_i^n$, $i=1,\dots,J_n$, are distinct, we have that $\#E(T_1^n(k))=\#E(\tilde{T}_1^n(k))+J_n$. Moreover, when $w_i^n\not\in V^n$ for any $i$, it is the case that $\#E(\mathcal{C}_1^n(k))=\#E({T}_1^n(k))+\#E^n$. Thus, for large $n$, $\#E(\mathcal{C}_1^n(k))=\#E(\tilde{T}_1^n(k))+J_n+\#E^n$. That the $u_i^n$, $i=1,\dots,J_n$ are distinct for large $n$ was noted above, and therefore $\#E^n=J_n$ for large $n$, which completes the proof.
\end{proof}

Whenever $\hat{T}_1^n(k)$ is a graph tree, by \cite{Kigamidendrite}, Theorem 5.4, we can define a resistance form $(\mathcal{E}_{\hat{T}_1^n(k)},\mathcal{F}_{\hat{T}_1^n(k)})$ corresponding to the real tree version of $\hat{T}_1^n(k)$ that satisfies (\ref{dtchar}). If it is also the case that $\#E(\hat{T}_1^n(k))=\#E(\mathcal{C}_1^n(k))$, then by applying the graph equivalence of $\mathcal{C}_1^n(k)$ and $\hat{T}_1^n(k)/ \sim_{n,k}$, it is readily checked that gluing together the real tree version of $\hat{T}_1^n(k)$ at $\{w_i^n,\tilde{w}_i^n\}$, $i=1,\dots,J_n$, yields the same compact length space as would be arrived at by including unit line segments along edges of $\mathcal{C}_1^n(k)$. Consequently, if both of these results are applicable, then by proceeding as in Section \ref{sectree}, we can construct a resistance form on $(\mathcal{E}_{\mathcal{C}_1^n(k)},\mathcal{F}_{\mathcal{C}_1^n(k)})$ on the compact length space version of $\mathcal{C}_1^n(k)$. We then define $X^{\mathcal{C}_1^n(k)}$ to be the process associated with $(\tfrac{1}{2}\mathcal{E}_{\mathcal{C}_1^n(k)},\mathcal{F}_{\mathcal{C}_1^n(k)})$ when considered as a Dirichlet form on $L^2(\mathcal{C}_1^n(k),\lambda_{\mathcal{C}_1^n(k)})$, where $\lambda_{\mathcal{C}_1^n(k)}$ is the one-dimensional Hausdorff measure on the compact length space version of $\mathcal{C}_1^n(k)$, which has total mass $\Lambda_n^{(k)}$. To be consistent with the notation of Section \ref{sectree}, we will write the law of $X^{\mathcal{C}_1^n(k)}$ started from ${\rho}$ as $\mathbf{P}^{\lambda_{\mathcal{C}_1^n(k)}}_{{\rho}}$. Also as in Section \ref{sectree}, it is possible to deduce that there exist jointly continuous local times $(L^{\mathcal{C}_1^n(k)}_t({x}))_{t\geq0,{x}\in\mathcal{C}_1^n(k)}$ for $X^{\mathcal{C}_1^n(k)}$. Importantly, applying results of Section \ref{seccont} and Lemma \ref{l63}, we can conclude the following lemma for $X^{\mathcal{C}_1^n(k)}$, $L^{\mathcal{C}_1^n(k)}$ and the related continuous additive functional
\begin{equation}\label{ac1nk}
\hat{A}^{\mathcal{C}_1^n(k)}_{t}:=\int_{\mathcal{C}_1^n(k)}L^{\mathcal{C}_1^n(k)}_t({x})\mu_{\mathcal{C}_1^n(k)}(d{x}).
\end{equation}

{\lem \label{l64} Suppose Assumption \ref{as1} holds and $k\geq J+1$.\\
(a) For every $c>0$, the joint laws of
\[\left(\left(\mathcal{C}_1^n(k),n^{-1/3}d_{\mathcal{C}_1^n},(X^{\mathcal{C}_1^n(k)}_{ctn^{1/3}\Lambda_n^{(k)}})_{t\in[0,1]}\right),\left(n^{-1}\hat{A}^{\mathcal{C}_1^n(k)}_{tn^{1/3}\Lambda_n^{(k)}}\right)_{t\geq 0}\right)\]
under $\mathbf{P}^{\lambda_{\mathcal{C}_1^n(k)}}_{{\rho}}$ converge to the joint law of
\[\left(\left(\mathcal{M}(k),d_{\mathcal{M}},(X^{\lambda_{\mathcal{M}(k)}}_{ct\Lambda^{(k)}})_{t\in [0,1]}\right),\left(\hat{A}^{\mathcal{M}(k)}_{t\Lambda^{(k)}}\right)_{t\geq 0}\right)\]
under $\mathbf{P}^{\lambda_{\mathcal{M}(k)}}_{\bar{\rho}}$ weakly as probability measures on the space $\mathbb{K}\times C(\mathbb{R}_+,\mathbb{R}_+)$.\\
(b) For $\varepsilon,t_0>0$,
\[\lim_{\delta\rightarrow 0}\limsup_{n\rightarrow\infty}\mathbf{P}^{\lambda_{\mathcal{C}_1^n(k)}}_{{\rho}}\left(n^{-1/3}
\sup_{\substack{{x},{y}\in\mathcal{C}^n_1(k):\\d_{\mathcal{C}^n_1}({x},{y})\leq \delta n^{1/3}}}\sup_{t\leq t_0n^{1/3}\Lambda_n^{(k)}}\left|L^{\mathcal{C}_1^n(k)}_t({x})- L^{\mathcal{C}_1^n(k)}_t({y})\right|>\varepsilon\right)=0.\]}
\begin{proof} Assume that $k\geq J+1$ and $n$ is large enough so that the conclusions of the previous lemma hold. Under these conditions, it is an elementary exercise to check that $\Upsilon_{\hat{T}_1^n(k),\mathcal{T}(k)}$ induces a homeomorphism $\Upsilon_{\mathcal{C}_1^n(k),\mathcal{M}(k)}:\mathcal{C}_1^n(k)\rightarrow \mathcal{M}(k)$ via the relationship, for $x\in\hat{T}_1^n(k)$,
\begin{equation}\label{induce2}
\Upsilon_{\mathcal{C}_1^n(k),\mathcal{M}(k)}(\phi_{n,k}(x))=\phi(\Upsilon_{\hat{T}_1^n(k),\mathcal{T}(k)}(x)),
\end{equation}
where $\phi_{n,k}$ is the canonical projection (with respect to the equivalence $\sim_{n,k}$) from the real tree version of $\hat{T}_1^n(k)$ to the compact length space version of $\mathcal{C}_1^n(k)$  (cf. (\ref{induce})). Since by the definition of $d^n_{\mathcal{T}}$ we have that
\[\mathcal{E}_{\mathcal{T}(k)}^n(f,f)=n^{1/3}\mathcal{E}_{\hat{T}_1^n(k)}(f\circ \Upsilon_{\hat{T}_1^n(k),\mathcal{T}(k)},f\circ\Upsilon_{\hat{T}_1^n(k),\mathcal{T}(k)}),\hspace{20pt}\forall f\in \mathcal{F}_{\mathcal{T}(k)},\]
where $(\mathcal{E}_{\mathcal{T}(k)}^n,\mathcal{F}_{\mathcal{T}(k)})$ is the resistance form associated with $(\mathcal{T}(k),d_{\mathcal{T}}^n)$ (its domain does not depend on $n$ in the range of $n$ that we are considering, see Section \ref{seccont}), it follows that
\[\mathcal{E}_{\mathcal{M}(k)}^n(f,f)=n^{1/3}\mathcal{E}_{\mathcal{C}_1^n(k)}(f\circ \Upsilon_{\mathcal{C}_1^n(k),\mathcal{M}(k)},f\circ\Upsilon_{\mathcal{C}_1^n(k),\mathcal{M}(k)}),\hspace{20pt}\forall f\in \mathcal{F}_{\mathcal{M}(k)},\]
where $(\mathcal{E}_{\mathcal{M}(k)}^n,\mathcal{F}_{\mathcal{M}(k)})$ is defined as in Section \ref{seccont} (cf. the proof of Proposition \ref{meas}). Thus the process $\Upsilon_{\mathcal{C}_1^n(k),\mathcal{M}(k)}(X^{\mathcal{C}_1^n(k)})$ is the diffusion on $\mathcal{M}(k)$ associated with $(\tfrac{1}{2}n^{-1/3}\mathcal{E}_{\mathcal{M}(k)}^n,\mathcal{F}_{\mathcal{M}(k)})$ considered as a Dirichlet form on $L^{2}(\mathcal{M}(k), \lambda_{\mathcal{C}_1^n(k)}\circ \Upsilon_{\mathcal{C}_1^n(k),\mathcal{M}(k)}^{-1})$. Now, from the fact that $\Upsilon_{\mathcal{C}_1^n(k),\mathcal{M}(k)}$ is actually a similitude of contraction ratio $n^{-1/3}$ with respect to the quotient metric $d_{\mathcal{M}}^n$ on $\mathcal{M}(k)$ induced by $d_{\mathcal{T}}^n$, we have that $n^{-1/3}\lambda_{\mathcal{C}_1^n(k)}\circ \Upsilon_{\mathcal{C}_1^n(k),\mathcal{M}(k)}^{-1}$ is equal to the one-dimensional Hausdorff measure $\lambda_{\mathcal{M}(k)}^n$ on $(\mathcal{M}(k),d^n_{\mathcal{M}(k)})$. Therefore two standard reparameterisations of time yield that $(\Upsilon_{\mathcal{C}_1^n(k),\mathcal{M}(k)}(X^{\mathcal{C}_1^n(k)}_{tn^{2/3}}))_{t\geq0}$ is the strong Markov process associated with $(\tfrac{1}{2}\mathcal{E}_{\mathcal{M}(k)}^n,\mathcal{F}_{\mathcal{M}(k)})$ considered as a Dirichlet form on $L^{2}(\mathcal{M}(k), \lambda_{\mathcal{M}(k)}^n)$. Hence, by Proposition \ref{convfin},
\[\left(\Upsilon_{\mathcal{C}_1^n(k),\mathcal{M}(k)}(X^{\mathcal{C}_1^n(k)}_{tn^{2/3}})\right)_{t\geq0}\rightarrow X^{\lambda_{\mathcal{M}(k)}}\]
in distribution in $C(\mathbb{R}_+,\mathcal{M}(k))$, where $X^{\lambda_{\mathcal{M}(k)}}$ was defined in Section \ref{seccont}. The convergence of the $\mathbb{K}$-coordinate in part (a) of the lemma is an easy consequence of this result and Lemma \ref{previous}(a).

For the other coordinate, first write the jointly continuous local times of the process \[(\Upsilon_{\mathcal{C}_1^n(k),\mathcal{M}(k)}(X^{\mathcal{C}_1^n(k)}_{tn^{2/3}}))_{t\geq0}\]
as
$(L^\Upsilon_t(\bar{x}))_{t\geq0,\bar{x}\in\mathcal{M}(k)}$, which are equal in distribution to the local times $L^{\mathcal{M}(k),n}$ defined in Section \ref{seccont}. Then, for any continuous function $f:\mathcal{M}(k)\rightarrow\mathbb{R}$,
\begin{eqnarray*}
\lefteqn{\int_{\mathcal{M}(k)}f(\bar{x})L^\Upsilon_t(\bar{x})\lambda_{\mathcal{M}(k)}^n(d\bar{x})}\\
&=& \int_{0}^t f\left(\Upsilon_{\mathcal{C}_1^n(k),\mathcal{M}(k)}(X^{\mathcal{C}_1^n(k)}_{sn^{2/3}})\right)ds\\
&=&n^{-2/3}\int_{0}^{tn^{2/3}} f\left(\Upsilon_{\mathcal{C}_1^n(k),\mathcal{M}(k)}(X^{\mathcal{C}_1^n(k)}_{s})\right)ds\\
&=&n^{-2/3}\int_{\mathcal{C}_1^n(k)}f\left(\Upsilon_{\mathcal{C}_1^n(k),\mathcal{M}(k)}({x})\right) L^{\mathcal{C}_1^n(k)}_{tn^{2/3}}({x})\lambda_{\mathcal{C}_1^n(k)}(d{x})\\
&=&n^{-1/3}\int_{\mathcal{M}(k)}f(\bar{x}) L^{\mathcal{C}_1^n(k)}_{tn^{2/3}}\left(\Upsilon_{\mathcal{C}_1^n(k),\mathcal{M}(k)}^{-1}(\bar{x})\right)\lambda_{\mathcal{M}(k)}^n(d\bar{x}).
\end{eqnarray*}
Consequently, we have that
\begin{equation}\label{ltlink}
n^{-1/3}L^{\mathcal{C}_1^n(k)}_{tn^{2/3}}({x})=
L^\Upsilon_t(\Upsilon_{\mathcal{C}_1^n(k),\mathcal{M}(k)}({x})),\hspace{20pt}\forall t\geq0,x\in\mathcal{C}_1^n(k),
\end{equation}
which implies in turn that
\begin{equation}\label{aexpression}
n^{-1}\hat{A}^{\mathcal{C}_1^n(k)}_{tn^{2/3}}=\int_{\mathcal{M}(k)}
L^\Upsilon_t(\bar{x})n^{-2/3}\mu_{\mathcal{C}_1^n(k)}\circ\Upsilon_{\mathcal{C}_1^n(k),\mathcal{M}(k)}^{-1}(d\bar{x}).
\end{equation}

Observe now that $\phi_{n,k}(x)=x$ for every $x\in V(\tilde{T}_1^n(k))$. Hence, by Lemma \ref{l63} and the identity at (\ref{induce2}),
\begin{eqnarray*}
\lefteqn{n^{-2/3}\mu_{\mathcal{C}_1^n}\circ \phi_{T_1^n,\tilde{T}_1^n(k)}^{-1}\circ\Upsilon_{\mathcal{C}_1^n(k),\mathcal{M}(k)}^{-1}} \\
 &=&n^{-2/3}\mu_{\mathcal{C}_1^n}\circ \phi_{T_1^n,\tilde{T}_1^n(k)}^{-1}\circ\Upsilon_{\hat{T}_1^n(k),\mathcal{T}(k)}^{-1}\circ\phi^{-1} \\
&\rightarrow& \mu_{\mathcal{T}(k)}\circ \phi^{-1}\\
&=&\mu_{\mathcal{M}(k)}
\end{eqnarray*}
weakly as measures on $\mathcal{M}(k)$. Since \[\sup_{x\in\mathcal{C}_1^n}d_{\mathcal{C}_1^n(k)}(\phi_{T_1^n,\tilde{T}_1^n(k)}(x),
\phi_{T_1^n,{T}_1^n(k)}(x))\leq 1,\]
we also have that
\[\sup_{x\in\mathcal{C}_1^n}d_{\mathcal{M}(k)}^n\left(
\Upsilon_{\mathcal{C}_1^n(k),\mathcal{M}(k)}(\phi_{T_1^n,\tilde{T}_1^n(k)}(x)),
\Upsilon_{\mathcal{C}_1^n(k),\mathcal{M}(k)}(\phi_{T_1^n,{T}_1^n(k)}(x))\right)\leq n^{-1/3},\]
and therefore $n^{-2/3}\mu_{\mathcal{C}_1^n(k)}\circ\Upsilon_{\mathcal{C}_1^n(k),\mathcal{M}(k)}^{-1} \rightarrow \mu_{\mathcal{M}(k)}$ weakly as measures on $\mathcal{M}(k)$.
Thus, recalling the expression for $\hat{A}^{\mathcal{C}_1^n(k)}$ at (\ref{aexpression}), the proof of part (a) is completed by applying Corollary \ref{twoproc} and Lemma \ref{previous}(a). Given the characterisation of $L^{\mathcal{C}_1^n(k)}$ at (\ref{ltlink}), part (b) is an immediate consequence of Lemmas \ref{tightness} and \ref{previous}(a).
\end{proof}

\section{Random walk scaling limit}\label{rwc}

We start this section by introducing notation that will allow us to state a quenched version of the main conclusion of this article: the simple random walk on $\mathcal{C}^1_n$ converges to the Brownian motion on $\mathcal{M}$. As an easy consequence of this result, which appears as Theorem \ref{quench}, we establish the corresponding annealed result, see Theorem \ref{anne}.

Let $\mathcal{C}_1^n$ be a fixed realisation of the largest connected component, rooted at its first ordered vertex, ${\rho}$ say. Let $X^{\mathcal{C}_1^n}=(X^{\mathcal{C}_1^n}_m)_{m\geq 0}$ be the discrete time simple random walk on $\mathcal{C}_1^n$ started from ${\rho}$, and denote by $\mathbf{P}_{{\rho}}^{\mathcal{C}_1^n}$ its law. For convenience, we will sometimes consider $\mathcal{C}_1^n$ as a compact length space by including unit line segments along edges. Moreover, by extending the definition of $X^{\mathcal{C}_1^n}$ to all positive times by linearly interpolating between integers, $\mathbf{P}_{{\rho}}^{\mathcal{C}_1^n}$ will sometimes be considered as a probability measure on $C(\mathbb{R}_+,\mathcal{C}_1^n)$.

Throughout this section, the rescaling operator $\Theta_n$ is defined on triples of the form $\mathcal{K}=(K,d_K,X^K)$, where $(K,d_K)$ is a non-empty compact length space and $X^K$ is a path in $C(\mathbb{R}_+,K)$, by setting
\[\Theta_n(\mathcal{K}):=\left(K,n^{-1/3}d_K,(X^K_{tn})_{t\in[0,1]}\right),\]
which can be considered as an element of $\mathbb{K}$. In particular, $\mathbf{P}^{\mathcal{C}_1^n}_{{\rho}}\circ\Theta_n^{-1}$ can be considered as a probability measure on $\mathbb{K}$, and we can prove the following limit for it as $n\rightarrow\infty$.

{\thm \label{quench} Under Assumption \ref{as1}, the law of the discrete time simple random walk on $\mathcal{C}_1^n$ started from ${\rho}$ satisfies
\[\mathbf{P}^{\mathcal{C}_1^n}_{{\rho}}\circ\Theta_n^{-1}\rightarrow\mathbf{P}^\mathcal{M}_{\bar{\rho}}\]
weakly in the space of probability measures on $\mathbb{K}$, where $\mathbf{P}^\mathcal{M}_{\bar{\rho}}$ is the law of the Brownian motion on $\mathcal{M}$ started from $\bar{\rho}$.}
\bigskip

To prove this result, we will proceed by a sequence of lemmas. Although these are relatively straightforward adaptations of the corresponding results for simple random walks on graph trees proved in \cite{Croydoncbp}, see also \cite{Croydoninf}, we include many of the details in an attempt to keep this article reasonably self-contained.

For the time being, suppose that Assumption \ref{as1} is satisfied and fix $k\geq J_n+1$. Consider the subgraph $\mathcal{C}_1^n(k)$ of $\mathcal{C}_1^n$ defined at (\ref{c1nk}). We set $X^{n,k}:=\phi_{\mathcal{C}_1^n,\mathcal{C}_1^n(k)}(X^{\mathcal{C}_1^n})$, where the projection map $\phi_{\mathcal{C}_1^n,\mathcal{C}_1^n(k)}$ was introduced at (\ref{cproj}). Let $A^{n,k}=(A^{n,k}_m)_{m\geq 0}$ be the jump times of $X^{n,k}$, or more precisely set $A^{n,k}_0=0$ and, for $m\geq 1$,
\[A^{n,k}_m:=\min\left\{l\geq A_{m-1}^{n,k}:X^{\mathcal{C}_1^n}_l\in \mathcal{C}_1^n(k)\backslash\left\{X^{\mathcal{C}_1^n}_{A_{m-1}^{n,k}}\right\}\right\}.\]
The jump-chain associated with $X^{n,k}$ is then given by $J^{n,k}=(J^{n,k}_m)_{m\geq0}$, where $J^{n,k}_m:=X^{n,k}_{A^{n,k}_m}$. Note that, by Lemma \ref{previous}(b), $J^{n,k}$ is the simple random walk on the vertices on $\mathcal{C}_1^n(k)$ started from ${\rho}$.
The discrete time inverse $\tau^{n,k}=(\tau^{n,k}(m))_{m\geq 0}$ of $A^{n,k}$ is defined by
\begin{equation}\label{taudef}
\tau^{n,k}(m):=\max\left\{l:A^{n,k}_l\leq m\right\},
\end{equation}
and we can check that $X^{n,k}$ can be recovered from $J^{n,k}$ and $\tau^{n,k}$ through the relationship
\begin{equation}\label{xnkdef}
X^{n,k}_m=J^{n,k}_{\tau^{n,k}(m)}.
\end{equation}
We define the local times of $J^{n,k}$ by setting
\[L^{n,k}_m({x}):=\frac{2}{{\rm deg}_{n,k}({x})}\sum_{l=0}^m\mathbf{1}_{\{{x}\}}\left(J^{n,k}_l\right),\]
for ${x}\in\mathcal{C}_1^n(k)$, where ${\rm deg}_{n,k}({x})$ is the usual graph degree of ${x}$ in $\mathcal{C}_1^n(k)$. We use these to define an additive functional $\hat{A}^{n,k}=(\hat{A}^{n,k}_m)_{m\geq 0}$ by setting $\hat{A}^{n,k}_0=0$ and, for $m\geq 1$,
\begin{equation}\label{hatank}
\hat{A}^{n,k}_m:=\int_{\mathcal{C}^n_1(k)}L^{n,k}_{m-1}({x})\mu_{\mathcal{C}_1^n(k)}(d{x}),
\end{equation}
where $\mu_{\mathcal{C}_1^n(k)}$ is the measure on $\mathcal{C}_1^n(k)$ introduced at (\ref{muc1nk}). The discrete time inverse of $\hat{A}^{n,k}$ will be denoted by $\hat{\tau}^{n,k}=(\hat{\tau}^{n,k}(m))_{m\geq 0}$ and defined similarly to (\ref{taudef}). This process is used to construct a time-changed version of $J^{n,k}$, $\hat{X}^{n,k}=(\hat{X}^{n,k}_m)_{m\geq0}$ say, by setting
\begin{equation}
\hat{X}^{n,k}_m:=J^{n,k}_{\hat{\tau}^{n,k}(m)}.\label{hatxnkdef}
\end{equation}
Similarly to the arguments of \cite{Croydoncbp} and \cite{Croydoninf}, it will be a goal to show that $X^{n,k}$ and $\hat{X}^{n,k}$ are close, which we do by demonstrating that the additive functionals $A^{n,k}$ and $\hat{A}^{n,k}$ are close, see Lemma \ref{akclose}. The first step is proving a tightness result for the local times $L^{n,k}$. In the statement of the result, we include the scaling factor $\Lambda_n^{(k)}:=\#E(\mathcal{C}_1^n(k))$, which was defined at (\ref{lambdankdef}) and will be useful later. Recall from Lemma \ref{previous}(a) that, under Assumption \ref{as1}, $n^{-1/3}\Lambda_n^{(k)}\rightarrow\Lambda^{(k)}$, where $\Lambda^{(k)}$ was defined at (\ref{Lambdak}).

{\lem \label{ltight} Suppose Assumption \ref{as1} holds. For $k\geq J+1$ and $\varepsilon,t_0>0$,
\[\lim_{\delta\rightarrow 0}\limsup_{n\rightarrow\infty}\mathbf{P}_{{\rho}}^{\mathcal{C}_1^n}\left(n^{-1/3}
\sup_{\substack{{x},{y}\in\mathcal{C}^n_1(k):\\d_{\mathcal{C}^n_1}({x},{y})\leq \delta n^{1/3}}}\sup_{m\leq t_0n^{1/3}\Lambda_n^{(k)}}\left|L^{n,k}_m({x})- L^{n,k}_m({y})\right|>\varepsilon\right)=0.\]}
\begin{proof} Fix $k\geq J+1$. Define $\hat{T}_1^n(k)$ as at (\ref{vhat}) and (\ref{ehat}) and suppose $n$ is large enough so that the conclusions of Lemma \ref{l63} hold. Analogously to the proof of Lemma \ref{tightness}, let $V^{n,k}\subseteq \hat{T}_1^n(k)$ be the set \[\{b^{\hat{T}_1^n(k)}(x,y,z):x,y,z\in\{\rho,\tilde{w}^n_1,\dots,\tilde{w}^n_{J_n},u^n_{J_n+1},\dots,u_k^n,w^n_1,\dots,w^n_{J_n}\}\},\] where the branch-point function $b^{\hat{T}_1^n(k)}$ is defined for the graph tree $\hat{T}_1^n(k)$ as at (\ref{bt}). Under Assumption \ref{as1}, by applying Lemma \ref{l63} it is possible to check that there exists a constant $\varepsilon_0>0$ such that
\begin{equation}\label{lowerlen}
n^{-1/3}\min_{{x},{y}\in V^{n,k}:x\neq y}d_{\hat{T}_1^n(k)}(x,y)\geq 2\varepsilon_0
\end{equation}
for large $n$. For the remainder of the proof, we assume that $n$ is large enough so that this bound holds and $J_n=J$. It follows from (\ref{lowerlen}) that for every ${x}\in\mathcal{C}_1^n(k)$ there exists an injective path $\Pi\subseteq \mathcal{C}_1^n(k)$ of length at least $n^{1/3}\varepsilon_0-1$ that has $x$ as an endpoint and satisfies $\Pi\cap \phi_{n,k}(V^{n,k})\subseteq\{{x}\}$. By considering the random walk $J^{n,k}$ observed on $\Pi$ and applying Lemma \ref{ltexp}, we can deduce that, for $t>0$,
\[\mathbf{P}_{{\rho}}^{\mathcal{C}_1^n}\left(n^{-1/3}L^{n,k}_{t_0n^{1/3}\Lambda_n^{(k)}}({x})\geq t\right)\leq c_1e^{-c_2pt{\rm deg}_{n,k}({x})},\]
where $c_1,c_2$ are constants that do not depend on $n$ or ${x}$, $p$ is the probability that $J^{n,k}_{m+1}\in\Pi$ given $J^{n,k}_m={x}$, and the ${\rm deg}_{n,k}({x})$ term arises as a result of the normalisation of the local times $L^{n,k}$. Since we trivially have $p={\rm deg}_{n,k}({x})^{-1}$, it is therefore the case that, for $t>0$,
\begin{equation}\label{exptail}
\sup_{{x}\in\mathcal{C}_1^n(k)}
\mathbf{P}_{{\rho}}^{\mathcal{C}_1^n}\left(n^{-1/3}
L^{n,k}_{t_0n^{1/3}\Lambda_n^{(k)}}({x})\geq t\right)\leq c_1e^{-c_2t},
\end{equation}
uniformly in $n$.

From this upper bound we will deduce that there exists a constant $c_3$ such that
\begin{equation}\label{diff}
\sup_{\substack{{x},{y}\in\mathcal{C}^n_1(k):\\d_{\mathcal{C}^n_1}({x},{y})\leq \delta n^{1/3}}}\mathbf{P}_{{\rho}}^{\mathcal{C}_1^n}\left(n^{-1/3}
\sup_{m\leq t_0n^{1/3}\Lambda_n^{(k)}}\left|L^{n,k}_m({x})- L^{n,k}_m({y})\right|>\varepsilon\right)\leq c_3\delta^2,
\end{equation}
uniformly in $n$, by adapting the proof of \cite{Croydoncbp}, Lemma 4.5, which is modelled in turn upon an argument from \cite{Borodin}. First, let ${x}\neq{y}\in\mathcal{C}_1^n(k)$ satisfy $d_{\mathcal{C}^n_1}({x},{y})\leq \delta n^{1/3}$ (note that in what follows we may assume that $\delta n^{1/3}\geq 1$, else the probability we are trying to bound is trivially equal to 0), and set $t_1:=\lfloor t_0n^{1/3}\Lambda_n^{(k)}\rfloor$. Conditional on the event that $J^{n,k}$ hits ${x}$ before ${y}$, we can write
\begin{equation}\label{twoterms2}
\sup_{m\leq t_1}\left|L^{n,k}_m({x})- L^{n,k}_m({y})+2\sum_{i=1}^{\ell_{m}^{n,k}({x})}\eta_i\right|\leq
\sup_{i\leq {\ell_{t_1}^{n,k}({x})}}2N_i{\rm deg}_{n,k}({y})^{-1},
\end{equation}
where $\ell_{m}^{n,k}({x}):=L^{n,k}_m({x}){\rm deg}_{n,k}({x})/2$ is the number of visits by $J^{n,k}$ to ${x}$ up to time $m$, $N_i$ is the number of visits by $J^{n,k}$ to ${y}$ between the $i$th and $(i+1)$st visits to ${x}$, and $\eta_i:=N_i{\rm deg}_{n,k}({y})^{-1}-{\rm deg}_{n,k}({x})^{-1}$. Clearly $(\eta_i)_{i\geq 1}$ is an independent, identically-distributed collection of random variables with zero mean (for verification of this final claim, see Section \ref{rwe}). We start by dealing with the sum on the left-hand side of (\ref{twoterms2}). Since $(\sum_{i=1}^m\eta_i)_{m\geq 1}$ is a martingale with respect to the filtration $(\mathcal{F}_m)_{m\geq 1}$, where $\mathcal{F}_m$ is the $\sigma$-algebra generated by $J^{n,k}$ up to the $(m+1)$st hitting time of ${x}$, and $\ell_{t_1}^{n,k}({x})$ is a stopping time for this filtration, Doob's martingale inequality implies that
\[\mathbf{P}_{{\rho}}^{\mathcal{C}_1^n}\left(\sup_{m\leq t_1}\left|\sum_{i=1}^{\ell_{m}^{n,k}({x})}\eta_i\right|>\varepsilon n^{1/3}\right)\leq c_4 n^{-4/3} \mathbf{E}_{{\rho}}^{\mathcal{C}_1^n}\left(\left|\sum_{i=1}^{\ell_{t_1}^{n,k}({x})}\eta_i\right|^4\right),\]
where $c_4$ is a constant that does not depend on the specific choice of ${x}$, ${y}$ or $n$. An  upper bound for the right-hand side in terms of the moments of $\ell_{t_1}^{n,k}({x})$ and $\eta_i$ can be obtained by a simple reworking of the argument that yields \cite{Borodin}, equation (1.29), which is a corresponding bound for simple random walk on the line. In particular, applying (\ref{exptail}) and Lemma \ref{etamoments}, the right-hand side can be bounded above by $c_5 \delta^2$, uniformly in ${x}$, ${y}$ and $n$. Still conditioning on the event that $J^{n,k}$ hits ${x}$ before ${y}$, for the right-hand side of (\ref{twoterms2}) we deduce that
\begin{eqnarray*}
\lefteqn{\mathbf{P}_{{\rho}}^{\mathcal{C}_1^n}\left(\sup_{i\leq {\ell_{t_1}^{n,k}({x})}}2N_i{\rm deg}_{n,k}({y})^{-1}>n^{1/3}\varepsilon\right)}\\
& \leq &
\mathbf{E}_{{\rho}}^{\mathcal{C}_1^n}\left(\sum_{i=1}^{{\ell_{t_1}^{n,k}({x})}}\mathbf{1}_{\{2N_i{\rm deg}_{n,k}({y})^{-1}>n^{1/3}\varepsilon\}}\right)\\
&\leq&\mathbf{E}_{{\rho}}^{\mathcal{C}_1^n}\left( \ell_{t_1}^{n,k}({x}) \right)\mathbf{P}_{{\rho}}^{\mathcal{C}_1^n}\left(2N_1{\rm deg}_{n,k}({y})^{-1}>n^{1/3}\varepsilon\right)\\
&\leq& c_6n^{-4/3}\mathbf{E}_{{\rho}}^{\mathcal{C}_1^n}\left( \ell_{t_1}^{n,k}({x}) \right) \left(1+\mathbf{E}_{{\rho}}^{\mathcal{C}_1^n}\left( \eta_1^4\right)\right)\\
&\leq & c_7 \delta^3,
\end{eqnarray*}
for some constant $c_7$ not depending on the specific choice of ${x}$, ${y}$ or $n$. To obtain the final bound here we have applied Lemma \ref{etamoments} to bound $\mathbf{E}_{{\rho}}^{\mathcal{C}_1^n}(\eta_1^4)$, and also used (\ref{exptail}) to deduce that $\mathbf{E}_{{\rho}}^{\mathcal{C}_1^n}( \ell_{t_1}^{n,k}({x}))\leq \deg_{n,k}(x)\mathbf{E}_{{\rho}}^{\mathcal{C}_1^n}( L_{t_1}^{n,k}({x}))\leq c_8 n^{1/3}$ (where we note that $\deg_{n,k}(x)$ can be crudely bounded above by $(k+1)(J_n+1)$). Putting these pieces together, it follows that
\[\mathbf{P}_{{\rho}}^{\mathcal{C}_1^n}\left(n^{-1/3}
\sup_{m\leq t_0n^{1/3}\Lambda_n^{(k)}}\left|L^{n,k}_m({x})- L^{n,k}_m({y})\right|>\varepsilon|J^{n,k}\mbox{ hits ${x}$ before ${y}$}\right)\leq c_9\delta^2,\]
uniformly in ${x}$, ${y}$ and $n$, and (\ref{diff}) easily follows.

Finally, appealing to the fact that the graph $\mathcal{C}_1^n(k)$ is a graph consisting of a collection of line segments, the number of which is bounded uniformly in $n$, the lemma follows from (\ref{diff}) by applying a standard maximal inequality from \cite{Bill2}, Section 10, exactly as in the proof of \cite{Croydoncbp}, Lemma 4.6.
\end{proof}

The next lemma demonstrates that the jump chains $J^{n,k}$ and additive functionals $\hat{A}^{n,k}$ converge when rescaled appropriately. We note that \[\left(\mathcal{C}_1^n(k),n^{-1/3}d_{\mathcal{C}_1^n},(J^{n,k}_{ctn^{1/3}\Lambda^{(k)}_n})_{t\geq 0}\right)\]
can be considered as an element of $\mathbb{K}$ by including unit line segments along edges in $\mathcal{C}_1^n(k)$ and linearly interpolating $J^{n,k}$. The definition of $\hat{A}^{n,k}$ is also extended to all positive times by linear interpolation.

{\lem \label{kconv} Suppose Assumption \ref{as1} holds. For every $k\geq J+1$ and $c>0$, the joint laws of
\[\left(\left(\mathcal{C}_1^n(k),n^{-1/3}d_{\mathcal{C}_1^n},(J^{n,k}_{ctn^{1/3}\Lambda_n^{(k)}})_{t\in[0,1]}\right),\left(n^{-1}\hat{A}^{n,k}_{tn^{1/3}\Lambda_n^{(k)}}\right)_{t\geq 0}\right)\]
under $\mathbf{P}^{\mathcal{C}_1^n}_{{\rho}}$ converge to the joint law of
\[\left(\left(\mathcal{M}(k),d_{\mathcal{M}},(X^{\lambda_{\mathcal{M}(k)}}_{ct\Lambda^{(k)}})_{t\in [0,1]}\right),\left(\hat{A}^{\mathcal{M}(k)}_{t\Lambda^{(k)}}\right)_{t\geq 0}\right)\]
under $\mathbf{P}^{\lambda_{\mathcal{M}(k)}}_{\bar{\rho}}$ weakly as probability measures on the space $\mathbb{K}\times C(\mathbb{R}_+,\mathbb{R}_+)$.}
\begin{proof} Recall the definition of $X^{\mathcal{C}_1^n(k)}$ from above Lemma \ref{l64}, and let
\[(h^{n,k}(m))_{m\geq 0}\]
be the hitting times of the graph vertices of $\mathcal{C}_1^n(k)$ by $X^{\mathcal{C}_1^n(k)}$.  Since we are assuming that $X^{\mathcal{C}_1^n(k)}_0={\rho}$, then $h^{n,k}(0)=0$. By the Markov property and trace theorem for Dirichlet forms, we can readily check that, conditional on $X^{\mathcal{C}_1^n(k)}_{h^{n,k}(m)}={x}$,
\[\left(X^{\mathcal{C}_1^n(k)}_{(h^{n,k}(m)+t)\wedge h^{n,k}(m+1)}\right)_{t\geq 0}\]
behaves exactly as a Brownian motion on a real tree star started from its internal vertex and stopped on hitting one of the ${\rm deg}_{n,k}({x})$ external vertices (see Section \ref{star} for a precise definition of such a process). Thus, by Lemma \ref{starlabel}, it is possible to assume that $J^{n,k}$ and $X^{\mathcal{C}_1^n(k)}$ are coupled so that
\[J^{n,k}_m=X^{\mathcal{C}_1^n(k)}_{h^{n,k}(m)}.\]
Hence, in light of Lemmas \ref{previous}(a) and \ref{l64}, to prove the convergence of the first coordinate in the statement of the lemma, it will be enough to establish that, for $t>0$,
\begin{equation}\label{hest}
\lim_{n\rightarrow\infty}n^{-2/3}\sup_{m\leq tn^{2/3}}\left|h^{n,k}(m)-m\right|=0
\end{equation}
in $\mathbf{P}^{\lambda_{\mathcal{C}_1^n(k)}}_{{\rho}}$-probability. Again applying Lemma \ref{starlabel}, under $\mathbf{P}^{\lambda_{\mathcal{C}_1^n(k)}}_{{\rho}}$ we have that the random variables in the sequence $(h^{n,k}(m+1)-h^{n,k}(m))_{m\geq0}$ are independent and identically distributed as the hitting time of $\{\pm 1\}$ by a standard Brownian motion on $\mathbb{R}$ started from 0, which is a random variable with mean 1 and finite fourth moments. Consequently, a standard martingale argument (cf. the proof of \cite{Croydoncbp}, Lemma 4.2) implies the desired result.

For the convergence of the second coordinate, we will start by showing that, for $t,\varepsilon>0$,
\begin{equation}\label{locest}
\lim_{n\rightarrow\infty}\sup_{{x}\in\mathcal{C}_1^n(k)}\mathbf{P}^{\lambda_{\mathcal{C}_1^n(k)}}_{{\rho}}\left( n^{-1/3}\sup_{m\leq t n^{1/3}\Lambda_n^{(k)}}\left|L^{n,k}_m({x})-L^{\mathcal{C}_1^n(k)}_{h^{n,k}(m)}({x})\right|>\varepsilon\right)=0,
\end{equation}
where $L^{\mathcal{C}_1^n(k)}$ are the local times of $X^{\mathcal{C}_1^n(k)}$ and the supremum is taken over the graph vertices of $\mathcal{C}_1^n(k)$. Fix an ${x}$ which is a graph vertex of $\mathcal{C}_1^n(k)$, let $\varsigma_i$ be the $i$th hitting time of ${x}$ by $J^{n,k}$ and set
\[\eta_i:=L^{\mathcal{C}_1^n(k)}_{h^{n,k}(\varsigma_i+1)}(x)-L^{\mathcal{C}_1^n(k)}_{h^{n,k}(\varsigma_i)}(x).\]
From Lemma \ref{starlabel} we find that $({\rm deg}_{n,k}({x})\eta_i/2)_{i\geq0}$ are independent and identically distributed as the local time at 0 of a standard Brownian motion started from 0 evaluated at the hitting time of $\{\pm 1\}$; a random variable with this distribution will be referred to as $Z$. By conditioning on $L^{n,k}_{\lfloor tn^{1/3}\Lambda_n^{(k)}\rfloor}(x)$, applying the fact that $Z$ is a random variable with mean 1 and finite fourth moments, and recalling (\ref{exptail}), it is possible to check that
\begin{equation}\label{locest1}
\lim_{n\rightarrow\infty} \sup_{x\in \mathcal{C}_1^n(k)}
\mathbf{P}^{\lambda_{\mathcal{C}_1^n(k)}}_{{\rho}}\left(
n^{-1/3}\sup_{m\leq t n^{1/3}\Lambda_n^{(k)}}
\left|\eta_1+\dots+\eta_{\deg_{n,k}(x)L^{n,k}_m({x})/2}-L^{n,k}_{m}({x})\right|>\varepsilon\right)=0,
\end{equation}
(cf. \cite{Croydoncbp}, equation (40)). Now, if $J^{n,k}_m=x$, then
\[\eta_1+\dots+\eta_{\deg_{n,k}(x)L^{n,k}_m({x})/2}=L^{\mathcal{C}_1^n(k)}_{h^{n,k}(m+1)}(x),\]
otherwise the sum is equal to $L^{\mathcal{C}_1^n(k)}_{h^{n,k}(m)}(x)$. Noting that the random variable
\[{\rm deg}_{n,k}({x})\left(L^{\mathcal{C}_1^n(k)}_{h^{n,k}(m+1)}(x)-L^{\mathcal{C}_1^n(k)}_{h^{n,k}(m)}(x)\right)/2\] has the same distribution as the $Z$ described above, it is thus clear that (\ref{locest1}) holds when $L^{n,k}_{m}({x})$ is replaced by $L^{\mathcal{C}_1^n(k)}_{h^{n,k}(m)}(x)$. That (\ref{locest}) holds follows. Moreover, by the tightness results of Lemmas \ref{l64}(b) and \ref{ltight}, it is possible to move the supremum inside the probability. Since, by the definitions of $\hat{A}^{\mathcal{C}_1^n(k)}$ and $\hat{A}^{n,k}$ (see (\ref{ac1nk}) and (\ref{hatank})),
\[\sup_{m\leq tn^{1/3}\Lambda_n^{(k)}}\left|\hat{A}^{n,k}_m-\hat{A}^{\mathcal{C}_1^n(k)}_{h^{n,k}(m)}\right|\leq Z^n_1\sup_{m\leq t n^{1/3}\Lambda_n^{(k)}}\sup_{{x}\in\mathcal{C}_1^n(k)}\left|L^{n,k}_m({x})-L^{\mathcal{C}_1^n(k)}_{h^{n,k}(m)}({x})\right|,\]
it follows that, for $t,\varepsilon>0$,
\[\lim_{n\rightarrow\infty}\mathbf{P}^{\lambda_{\mathcal{C}_1^n(k)}}_{{\rho}}\left( n^{-1}\sup_{m\leq tn^{1/3}\Lambda_n^{(k)}}\left|\hat{A}^{n,k}_m-\hat{A}^{\mathcal{C}_1^n(k)}_{h^{n,k}(m)}\right|>\varepsilon\right)=0,\]
where we recall that, under Assumption 1, $n^{-2/3}Z_1^n\rightarrow Z_1$. Thus, by Lemma \ref{l64} and (\ref{hest}), the proof is complete.
\end{proof}

We now establish a tightness result for $A^{n,k}$ and $\hat{A}^{n,k}$.

{\lem \label{akclose} Suppose Assumption \ref{as1} holds. For $t_0,\varepsilon>0$,
\[\lim_{k\rightarrow\infty}\limsup_{n\rightarrow\infty}\mathbf{P}^{\mathcal{C}_1^n}_{{\rho}}\left(n^{-1}\sup_{m\leq t_0n^{1/3}\Lambda_n^{(k)}}\left|A^{n,k}_m-\hat{A}^{n,k}_m\right|>\varepsilon\right)=0.\]}
\begin{proof} Fix $k\geq J+1$. We can write
\begin{eqnarray}
\left|A^{n,k}_m-\hat{A}^{n,k}_m\right|&\leq&\left|\sum_{l=0}^{m-1}\left(A^{n,k}_{l+1}-A^{n,k}_l-\mathbf{E}^{\mathcal{C}^n_1}_{{\rho}}(A^{n,k}_{l+1}-A^{n,k}_l|J^{n,k})\right)\right|\nonumber\\
&&+\sum_{l=0}^{m-1}\left|\mathbf{E}^{\mathcal{C}^n_1}_{{\rho}}(A^{n,k}_{l+1}-A^{n,k}_l|J^{n,k})-\frac{2\mu_{\mathcal{C}_1^n(k)}(\{J^{n,k}_l\})}{{\rm deg}_{n,k}(J^{n,k}_l)}\right|.\label{twoterms3}
\end{eqnarray}
By Lemma \ref{previous}(b), conditional on $J^{n,k}$, the random variable $A^{n,k}_{l+1}-A^{n,k}_l$ is precisely the time taken by a random walk started at the root of a tree with $\mu_{\mathcal{C}_1^n(k)}(\{{J^{n,k}_l}\})$ vertices to leave via one of ${\rm deg}_{n,k}(J^{n,k}_l)$ additional vertices, each attached to the root of this tree via a single edge. Consequently, we can deduce from  \cite{Croydoncbp}, Lemma B.3, that
\begin{equation}\label{mean}
\mathbf{E}^{\mathcal{C}^n_1}_{{\rho}}(A^{n,k}_{l+1}-A^{n,k}_l|J^{n,k})= \frac{2\mu_{\mathcal{C}_1^n(k)}(\{J^{n,k}_l\})-2+{\rm deg}_{n,k}(J^{n,k}_l)}{{\rm deg}_{n,k}(J^{n,k})},
\end{equation}
\begin{equation}\label{meansquare}
\mathbf{E}^{\mathcal{C}^n_1}_{{\rho}}((A^{n,k}_{l+1}-A^{n,k}_l)^2|J^{n,k}_l)\leq 36\left({\rm deg}_{n,k}(J^{n,k}_l)+\Delta_n^{(k)}\right)\frac{\mu_{\mathcal{C}_1^n(k)}(\{J^{n,k}_l\})^2}{{\rm deg}_{n,k}(J^{n,k}_l)},
\end{equation}
where $\Delta_n^{(k)}$ was defined at (\ref{deltankdef}). For the second term of (\ref{twoterms3}), substituting in the expression at (\ref{mean}) yields
\[\sup_{m\leq t_0n^{1/3}\Lambda_n^{(k)}}\sum_{l=0}^{m-1}\left|\mathbf{E}^{\mathcal{C}^n_1}_{{\rho}}(A^{n,k}_{l+1}-A^{n,k}_l|J^{n,k})-\frac{2\mu_{\mathcal{C}_1^n(k)}(\{J^{n,k}_l\})}{{\rm deg}_{n,k}(J^{n,k}_l)}\right|\leq t_0n^{1/3}\Lambda_n^{(k)}.\]
That the upper bound converges to 0 when multiplied by $n^{-1}$ follows from Lemma \ref{previous}(a). For the first term of (\ref{twoterms3}), we apply Kolmogorov's maximal inequality (see \cite{Kallenberg}, Lemma 4.15) and (\ref{meansquare}) to obtain that
\begin{eqnarray}
&&\label{condprob}\\
\lefteqn{\mathbf{P}^{\mathcal{C}_1^n}_{{\rho}}\left(n^{-1}\sup_{m\leq t_0n^{1/3}\Lambda_n^{(k)}}\left|\sum_{l=0}^{m-1}\left(A^{n,k}_{l+1}-A^{n,k}_l-\mathbf{E}^{\mathcal{C}^n_1}_{{\rho}}(A^{n,k}_{l+1}- A^{n,k}_l|J^{n,k}_l)\right)\right|>\varepsilon\:|\:J^{n,k}\right)}\nonumber\\
&\leq&\frac{1}{n^{2}\varepsilon^2}\sum_{l=0}^{\lfloor t_0n^{1/3}\Lambda_n^{(k)}\rfloor-1}\mathbf{E}^{\mathcal{C}^n_1}_{{\rho}}((A^{n,k}_{l+1}-A^{n,k}_l)^2\:|\:J^{n,k}_l)\nonumber\\
&\leq& \frac{18 Z_1^n}{n^{2}\varepsilon^2}\left(\max_{{x}\in\mathcal{C}_1^n(k)}{\rm deg}_{n,k}({x})+\Delta_n^{(k)}\right)\hat{A}^{n,k}_{\lfloor t_0n^{1/3}\Lambda_n^{(k)}\rfloor}.\nonumber\hspace{140pt}
\end{eqnarray}
For $t,\delta>0$,
\begin{eqnarray}
&&\label{hello}\\
\lefteqn{\limsup_{k\rightarrow\infty}\limsup_{n\rightarrow\infty}
\mathbf{P}_{{\rho}}^{\mathcal{C}_1^n}\left(n^{-2}Z_1^n
\left(\max_{{x}\in\mathcal{C}_1^n(k)}{\rm deg}_{n,k}({x})+\Delta_n^{(k)}\right)\hat{A}^{n,k}_{\lfloor t_0n^{1/3}\Lambda_n^{(k)}\rfloor}>\delta\right)}\nonumber\\
&\leq& \limsup_{k\rightarrow\infty}\limsup_{n\rightarrow\infty} t\delta^{-1} n^{-1} Z^n_1 \left(\max_{{x}\in\mathcal{C}_1^n(k)}{\rm deg}_{n,k}({x})+\Delta_n^{(k)}\right)\nonumber\\
&&+\limsup_{k\rightarrow\infty}\limsup_{n\rightarrow\infty}\mathbf{P}_{{\rho}}^{\mathcal{C}_1^n}\left(n^{-1}\hat{A}^{n,k}_{\lfloor t_0n^{1/3}\Lambda_n^{(k)}\rfloor}>t\right).\hspace{120pt}\nonumber
\end{eqnarray}
As noted in the proof of Lemma \ref{ltight}, the maximum degree of a vertex in $\mathcal{C}_1^n(k)$ can be bounded above by $(k+1)(J_n+1)$. Hence, that the first term is equal to 0 is a consequence of Assumption \ref{as1} and Lemma \ref{previous}(a). From Lemmas \ref{akprops} and \ref{kconv}, the second term can be made arbitrarily small by choosing $t$ suitably large. Hence the expression at (\ref{hello}) is equal to 0. In particular, we have shown that the conditional probability at (\ref{condprob}) converges to 0 as $n$ and then $k$ tend to infinity. Proving from this that the unconditional probability satisfies the same result is elementary (cf. \cite{Kallenberg}, Exercise 6.11).
\end{proof}

\begin{proof}[Proof of Theorem \ref{quench}] The proof of this result is almost identical to that of \cite{Croydoninf}, Proposition 4.1. First note that since $\hat{A}^{\mathcal{M}(k)}$ is $\mathbf{P}_{\bar{\rho}}^{\lambda_{\mathcal{M}(k)}}$-a.s. continuous and strictly increasing (see Lemma \ref{akprops}), the continuous mapping theorem implies that Lemma \ref{kconv} holds when $n^{-1}\hat{A}^{n,k}_{tn^{1/3}\Lambda_n^{(k)}}$ is replaced by its inverse $({n^{1/3}\Lambda_n^{(k)}})^{-1}\hat{\tau}^{n,k}(tn)$ and $\hat{A}^{\mathcal{M}(k)}_{t\Lambda^{(k)}}$ is replaced by $(\Lambda^{(k)})^{-1}\hat{\tau}^{\mathcal{M}(k)}(t)$. Thus, by applying the representations of $\hat{X}^{n,k}$ and $X^{\mu_\mathcal{M}(k)}$ at (\ref{hatxnkdef}) and Lemma \ref{timechange} respectively, we obtain that
\[\Theta_n\left(\mathcal{C}_1^n(k),d_{\mathcal{C}_1^n},\hat{X}^{n,k}\right)\rightarrow \left(\mathcal{M}(k),d_\mathcal{M},X^{\mu_\mathcal{M}(k)}\right)\]
weakly in $\mathbb{K}$ as $n\rightarrow\infty$ for every $k\geq J+1$. Furthermore, from Proposition \ref{kinfty} we have that
\[\left(\mathcal{M}(k),d_\mathcal{M},X^{\mu_\mathcal{M}(k)}\right)\rightarrow \left(\mathcal{M},d_\mathcal{M},X^{\mathcal{M}}\right)\]
weakly in $\mathbb{K}$ as $k\rightarrow\infty$. Therefore to complete the proof, by applying \cite{Bill2}, Theorem 3.2, for example, it will suffice to demonstrate that: for $\varepsilon>0$,
\[\lim_{k\rightarrow\infty}\limsup_{n\rightarrow\infty}\mathbf{P}_{{\rho}}^{\mathcal{C}_1^n}\left( n^{-1/3}\sup_{m\leq n}d_{\mathcal{C}_1^n}\left(X^{\mathcal{C}_1^n}_m,\hat{X}^{n,k}_m\right)>\varepsilon\right)=0.\]
By Lemma \ref{previous}(a), we can immediately replace $X^{\mathcal{C}_1^n}$ by $X^{n,k}$ in this requirement. On recalling that the latter process can be expressed as at (\ref{xnkdef}), the result is a simple exercise in analysis involving the application of the above convergence results and Lemma \ref{akclose},  (cf. \cite{Croydoninf}, Proposition 4.1).
\end{proof}

We now prove the annealed version of Theorem \ref{quench}. For the remainder of this section, it is supposed that $\mathcal{C}_1^n$ is the largest connected  component of the random graph $G(n,p)$, where $p=n^{-1}+\lambda n^{-4/3}$. Defining the annealed law of the rescaled simple random walk on $\mathcal{C}_1^n$ started from ${\rho}$ by setting
\[\mathbb{P}_n(A):=\int \mathbf{P}^{\mathcal{C}_1^n}_{{\rho}}\circ\Theta_n^{-1}(A)\mathbf{P}\left(d\mathcal{C}_1^n\right),\]
for measurable $A\subseteq\mathbb{K}$, the above quenched theorem easily yields the following result.

{\thm \label{anne} The annealed law of the rescaled simple random walk on $\mathcal{C}_1^n$ started from ${\rho}$ satisfies
\[\mathbb{P}_n\rightarrow\mathbb{P}^\mathcal{M}_{\bar{\rho}}\]
weakly in the space of probability measures on $\mathbb{K}$, where $\mathbb{P}^\mathcal{M}_{\bar{\rho}}$ is the annealed law of the Brownian motion on $\mathcal{M}$ started from $\bar{\rho}$, as defined at (\ref{ann}).}
\begin{proof} From Theorem \ref{dhc} and the fact that $\mathbf{P}((\tilde{e}^{(Z_1)},\mathcal{P},\xi)\in\Gamma)=1$, we can assume that the random graphs $(\mathcal{C}_1^n)_{n\geq 1}$ have been constructed in such a way that Assumption \ref{as1} holds $\mathbf{P}$-a.s. Consequently, applying the dominated convergence theorem in combination with Theorem \ref{quench} yields the result.
\end{proof}

To complete this section, let us remark that it is possible to prove a corresponding result for the sequence of simple random walks on the collection of components of the random graph $G(n,p)$. In \cite{ABG}, Theorem 24, it is described how the sequence of connected components of $G(n,p)$ in the critical window arranged so that the number of vertices of the components are non-increasing, $(\mathcal{C}_1^n,\mathcal{C}_2^n,\dots)$ say, can be rescaled to converge in distribution with respect to a fourth-order Gromov Hausdorff distance for sequences of compact metric spaces. The limit object, $(\mathcal{M}_1,\mathcal{M}_2,\dots)$ say, consists of a collection of compact metric spaces $\mathcal{M}_i=(\mathcal{M}_i,d_{\mathcal{M}_i})$, $i\geq 1$, that are each distributed as $\mathcal{M}$ up to a random scaling factor. Since $d_\mathbb{K}(\mathcal{K},\mathcal{K}')\leq {\rm diam}(K,d_K)+{\rm diam}(K',d_{K'})$, we can readily adapt the proof of this result to deduce the analogous conclusion for the associated random walks. Specifically, suppose that, conditional on $(\mathcal{C}_1^n,\mathcal{C}_2^n,\dots)$, the processes $(X^{\mathcal{C}_i^n})_{i\geq 1}$ are independent simple random walks on the components $\mathcal{C}_i^n$, each started from the first ordered vertex of the relevant component. Note, in the case that $\mathcal{C}_i^n$ is empty, we replace it with a metric space consisting of a single point and $X^{\mathcal{C}_i^n}$ by a constant process. Then, by applying Theorem \ref{anne}, one can show that the annealed law of the sequence $\{\Theta_n(\mathcal{C}_i^n,d_{\mathcal{C}_i^n},X^{\mathcal{C}_i^n})\}_{i\geq 1}$ converges to the annealed law of $\{(\mathcal{M}_i,d_{\mathcal{M}_i},X^{\mathcal{M}_i})\}_{i\geq 1}$, where $d_{\mathcal{M}_i}$ is the metric on $\mathcal{M}_i$ and, conditional on the metric spaces, $(X^{\mathcal{M}_i})_{i\geq 1}$ are independent Brownian motions on the spaces $(\mathcal{M}_i,d_{\mathcal{M}_i})$ started from a root vertex, with respect to the topology induced by the metric
\[d^{(4)}_{\mathbb{K}}\left((\mathcal{K}_i)_{i\geq 1},(\mathcal{K}'_i)_{i\geq 1}\right):=\left(\sum_{i=1}^\infty d_\mathbb{K}(\mathcal{K}_i,\mathcal{K}_i')^4\right)^{1/4}\]
on sequences of elements of $\mathbb{K}$.

\section{Properties of the limiting process}\label{propsec}

We have already described how the Brownian motion on the scaling limit of the largest connected component of the critical random graph, $X^\mathcal{M}$, can be constructed as a $\mu_\mathcal{M}$-symmetric Markov diffusion for $\mathbf{P}$-a.e. realisation of $\mathcal{M}$. To complete this article we will explain how to transfer the short-time asymptotic results for the quenched heat kernel of the Brownian motion on the continuum random tree obtained in \cite{Croydoncrt} to the transition density of $X^\mathcal{M}$. The key to doing this is the subsequent lemma, which is proved in the deterministic setting of Section \ref{sectree} and allows comparison of volume and resistance properties of $\mathcal{M}$ and $\mathcal{T}$. Note that, for a resistance form $(\mathcal{E},\mathcal{F})$ on a set $X$, the associated resistance operator $R$ can be extended to disjoint subsets $A,B\subseteq X$ by setting \[R(A,B):=\{\mathcal{E}(f,f):f\in\mathcal{F},f|_A=0,f|_B=1\}^{-1}.\]

{\lem In the setting of Section \ref{sectree}, the following results hold.\\
(a) For every $x\in\mathcal{T}$ and $r>0$,
\[\mu_\mathcal{T}\left(B_\mathcal{T}(x,r)\right)\leq \mu_\mathcal{M}\left(B_{(\mathcal{M},d_\mathcal{M})}(\bar{x},r)\right)\leq (2J+1)\sup_{y\in\mathcal{T}}\mu_\mathcal{T}\left(B_\mathcal{T}(y,r)\right).\]
(b) For $x\not\in\cup_{i=1}^J\{u_i,v_i\}$,
\begin{equation}\label{volasymp}
\mu_\mathcal{M}\left(B_{(\mathcal{M},d_\mathcal{M})}(\bar{x},r)\right)\sim \mu_\mathcal{T}\left(B_\mathcal{T}(x,r)\right),
\end{equation}
\begin{equation}\label{resasymp}
R_\mathcal{M}\left(\bar{x},B_{(\mathcal{M},d_\mathcal{M})}(\bar{x},r)^c\right)\sim R_\mathcal{T}\left(x,B_\mathcal{T}(x,r)^c\right),
\end{equation}
as $r\rightarrow 0$, where $R_\mathcal{T}$ is the resistance operator associated with $(\mathcal{E}_\mathcal{T},\mathcal{F}_\mathcal{T})$.\\
(c) The metric $R_\mathcal{M}$ satisfies the chaining condition: there exists a finite constant $c$ such that for all $\bar{x},\bar{y}\in\mathcal{M}$ and $n\in\mathbb{N}$, there exist $\bar{x}_0,\dots,\bar{x}_n\in\mathcal{M}$ with $\bar{x}_0=\bar{x}$ and $\bar{x}_n=\bar{y}$ such that
\[R_\mathcal{M}(\bar{x}_{i-1},\bar{x}_i)\leq \frac{cR_\mathcal{M}(\bar{x},\bar{y})}{n},\]
for every $i=1,\dots,n$.}
\begin{proof} The left-hand inequality of (a) is an immediate consequence of the readily deduced fact that $B_\mathcal{T}(x,r)\subseteq \phi^{-1}(B_{(\mathcal{M},d_\mathcal{M})}(\bar{x},r))$. Now let $x,y\in\mathcal{T}$. By the definition of $d_\mathcal{M}$, for every $\varepsilon>0$ there exist vertices $x_i,y_i\in\mathcal{T}$, $i=1,\dots,k$, satisfying $\bar{x}_1=\bar{x}$, $\bar{y}_i=\bar{x}_{i+1}$, $\bar{y}_k=\bar{y}$ and $\sum_{i=1}^kd_\mathcal{T}(x_i,y_i)\leq d_\mathcal{M}(\bar{x},\bar{y})+\varepsilon$. Note that if $x_{i+1}\not\in\cup_{j=1}^J\{u_j,v_j\}$ for some $i=1,\dots,k-1$, then $\bar{x}_{i+1}=\{x_{i+1}\}$ and so $y_i=x_{i+1}$. Therefore, by the triangle inequality for $d_\mathcal{T}$,
\[d_\mathcal{T}(x_i,y_{i+1})\leq d_\mathcal{T}(x_i,y_{i})+d_\mathcal{T}(y_i,y_{i+1})=d_\mathcal{T}(x_i,y_{i})+d_\mathcal{T}(x_{i+1},y_{i+1}),\]
which means that we can obtain a shorter sequence of pairs of vertices
with the same properties as above by replacing the two pairs $(x_i,y_i)$ and $(x_{i+1},y_{i+1})$ by the single pair $(x_i,y_{i+1})$. In particular, it follows that we can assume the vertices have been chosen to satisfy $x_{i+1},y_i\in\cup_{j=1}^J\{u_j,v_j\}$ for $i=1,\dots,k-1$. Using this observation, it is elementary to establish that
\[\phi^{-1}(B_{(\mathcal{M},d_\mathcal{M})}(\bar{x},r))\subseteq B_\mathcal{T}(x,r)\cup \left(\cup_{y\in\cup_{i=1}^J\{u_i,v_i\}}B_\mathcal{T}(y,r)\right),\]
from which the right-hand side of (a) follows.

For part (b), since $\cup_{i=1}^J\{u_i,v_i\}$ is a closed set, for every $x\not \in\cup_{i=1}^J\{u_i,v_i\}$, there exists an $r_0>0$ such that $B_\mathcal{T}(x,r_0)\subseteq \mathcal{T}\backslash \cup_{i=1}^J\{u_i,v_i\}$, and we can therefore use the argument of the previous paragraph to check that \[\phi^{-1}(B_{(\mathcal{M},d_\mathcal{M})}(\bar{x},r))= B_\mathcal{T}(x,r)\] for $r<r_0$. Hence $\mu_\mathcal{M}(B_{(\mathcal{M},d_\mathcal{M})}(\bar{x},r))=\mu_\mathcal{T}(B_\mathcal{T}(x,r))$ for $r<r_0$, which proves (\ref{volasymp}). Moreover, if $r<r_0$ and $f$ is a function on $\mathcal{T}$ satisfying $f|_{B_\mathcal{T}(x,r)^c}=1$, then $f$ immediately satisfies $f(u_i)=f(v_i)$, $i=1,\dots,J$. Thus, for $r<r_0$,
\begin{eqnarray*}
\lefteqn{R_\mathcal{T}\left(x,B_\mathcal{T}(x,r)^c\right)^{-1}}\\
&=& \inf\{\mathcal{E}_\mathcal{T}(f,f):f\in\mathcal{F}_\mathcal{T},f(x)=0,f|_{B_\mathcal{T}(x,r)^c}=1,f(u_i)=f(v_i),i=1,\dots,J\}\\
&=&\inf\{\mathcal{E}_\mathcal{M}(f,f):f\in\mathcal{F}_\mathcal{M},f(\bar{x})=0,f|_{B_\mathcal{M}(\bar{x},r)^c}=1\}\\
&=&R_\mathcal{M}\left(\bar{x},B_{(\mathcal{M},d_\mathcal{M})}(\bar{x},r)^c\right)^{-1},
\end{eqnarray*}
where the second equality again applies the result that $\phi^{-1}(B_{(\mathcal{M},d_\mathcal{M})}(\bar{x},r))= B_\mathcal{T}(x,r)$ for $r<r_0$, and this establishes (\ref{resasymp}).

Finally, let $\bar{x},\bar{y}\in\mathcal{M}$. By the definition of $d_\mathcal{M}$, there exist vertices $x_i,y_i\in\mathcal{T}$, $i=1,\dots,k$, such that $\bar{x}_1=\bar{x}$, $\bar{y}_i=\bar{x}_{i+1}$, $\bar{y}_k=\bar{y}$, and also $\sum_{i=1}^kd_\mathcal{T}(x_i,y_i)\leq 2d_\mathcal{M}(\bar{x},\bar{y})$. Define $t_k$ and $\phi\circ \gamma:[0,t_k]\rightarrow \mathcal{M}$ from these points as in Lemma \ref{lengthspace}. Given $n\in\mathbb{N}$, if we let $\bar{z}_i:=\phi\circ\gamma(it_k/n)$, $i=0,1,\dots, n$, then
\[R_\mathcal{M}(\bar{z}_{i-1},\bar{z}_i)\leq d_\mathcal{M}(\bar{z}_{i-1},\bar{z}_i)\leq \frac{t_k}{n}\leq \frac{2d_\mathcal{M}(\bar{x},\bar{y})}{n}\leq \frac{2c^{-1}R_\mathcal{M}(\bar{x},\bar{y})}{n},\]
where $c$ is the constant of Lemma \ref{compmet}, which proves (c).
\end{proof}

Now let $\mathcal{M}$ be the random scaling limit of the largest connected  component of the critical random graph. By the description of $\mathcal{M}$ in terms of the tilted continuum random tree and the above result, it is possible to deduce that the measure $\mu_\mathcal{M}$ satisfies precisely the same $\mathbf{P}$-a.s. global (uniform) and local (pointwise) estimates that were proved for $\mu_\mathcal{T}$ in the case that $\mathcal{T}$ is the continuum random tree in \cite{Croydoncrt}, Theorems 1.2 and 1.3. In particular, in the $\mu_\mathcal{M}$-measure of balls of radius $r$, there are global logarithmic fluctuations about a leading order $r^2$ term, and local fluctuations of log-logarithmic order as $r\rightarrow 0$. Consequently, noting the bounds for $R_\mathcal{M}$ in terms of $d_\mathcal{M}$ obtained in Lemma \ref{compmet} and the chaining condition of the previous lemma, we can proceed as in \cite{Croydoncrt} to deduce bounds and fluctuation results for the quenched transition density of $X^\mathcal{M}$ by applying the conclusions regarding transition densities of processes associated with general resistance forms obtained in \cite{Croydon}. For example, it is $\mathbf{P}$-a.s. the case that
$(p^\mathcal{M}_t(\bar{x},\bar{y}))_{\bar{x},\bar{y}\in\mathcal{M},t>0}$ satisfies, for some random constants $c_1,c_2,c_3,c_4$, $t_0>0$ and deterministic $\theta_1,\theta_2,\theta_3\in(0,\infty)$,
\[p^\mathcal{M}_t(\bar{x},\bar{y})\geq c_1 t^{-\frac{2}{3}}(\ln_1 t^{-1})^{-\theta_1}\exp \left\{-c_2 \left(\frac{d^{3}}{t}\right)^{1/2}\ln_1 \left(\frac{d}{t}\right)^{\theta_2}\right\},\]
and
\[p^\mathcal{M}_t(\bar{x},\bar{y})\leq c_3 t^{-\frac{2}{3}}(\ln_1 t^{-1})^{1/3}\exp \left\{-c_4 \left(\frac{d^{3}}{t}\right)^{1/2}\ln_1 \left(\frac{d}{t}\right)^{-\theta_3}\right\},\]
for all $\bar{x},\bar{y}\in\mathcal{M}$, $t\in(0,t_0)$, where $d:=d_{\mathcal{M}}(\bar{x},\bar{y})$ and $\ln_1 x := 1\vee\ln x$ (cf. \cite{Croydoncrt}, Theorem 1.4), which confirms the spectral dimension result of (\ref{specdim}). Moreover, if we consider the on-diagonal part of the transition density $p^\mathcal{M}_t(\bar{x},\bar{x})$, then exactly as in \cite{Croydoncrt}, Theorems 1.5 and 1.6, we obtain that global logarithmic fluctuations occur, and no more than log-logarithmic fluctuations occur locally, $\mathbf{P}$-a.s. Note that it is in proving the result analogous to \cite{Croydoncrt}, Theorem 1.6, which demonstrates that the local transition density fluctuations are of at most log-logarithmic order, that the resistance asymptotics of (\ref{resasymp}) are needed. Finally, observe that \cite{Croydoncrt}, Proposition 1.7, gives estimates for the annealed transition density asymptotics at the root of the continuum random tree, which demonstrate that no fluctuations occur when the environment is averaged out. Since proving the corresponding bounds for $\mathcal{M}$ will require a careful consideration of the distribution of the random variable $J$, or at least the number of ``glued'' points in a neighbourhood of the root, we leave the pursuit of such results as a project for the future.

\begin{appendices}
\section{Appendix}

This section collects together a number of technical results that are applied in the proofs of our main theorems. Although they are important to our arguments, we consider that they rather disrupt the flow of the article to appear where they are used.

\subsection{Resistance lower bound on finite graphs}

Let $G=(V,E)$ be a connected graph with finite vertex and edge sets. We allow the possibility of multiple edges between pairs of vertices and loops (edges whose endpoints are equal). To each edge $e\in E$, assign a resistance $r_e\in(0,\infty)$, and denote by $R_{G}$ the effective resistance metric on the resulting electrical network. In the following lemma we provide a lower bound for $R_G$ in terms of the shortest path metric $d_G$ for the graph $G$ with edges weighted according to $(r_e)_{e\in E}$, which is applied in Lemma \ref{compmet}.

{\lem \label{reslower} For any $x,y\in V$,
\[R_G(x,y)\geq \frac{d_G(x,y)}{\#E!}.\]}
\begin{proof} Fix $x\neq y\in V$. A path $\pi$ from $x$ to $y$ is a finite sequence of vertices $v_0,v_1,\dots,v_k$ such that $v_0=x$, $v_k=y$ and $\{v_{j-1},v_j\}\in E$ for every $j=1,\dots,k$. For such a path $\pi$, we will write $e\in \pi$ if $e=\{v_{j-1},v_j\}$ for some $j=1,\dots,k$. If $\Pi$ is the collection of all paths from $x$ to $y$, then it is known that
\begin{equation}\label{rsum}
R_{G}(x,y)=\inf_{(p_\pi)_{\pi\in\Pi}}\sum_{e\in E} r_e \left(\sum_{\pi\in\Pi}p_\pi \mathbf{1}_{\{e\in\pi\}}\right)^2,
\end{equation}
where the infimum is taken over all probability measures $(p_\pi)_{\pi\in\Pi}$ on $\Pi$ (see \cite{LP}, Exercise 3.20, for example). In fact, it is readily checked that the same result holds when $\Pi$ is restricted to only edge-simple paths (i.e. those paths which do not pass the same edge in either direction more than once), and so we will henceforth assume that this is the case. Since all the terms in the final summand of (\ref{rsum}) are positive, it follows that
\[R_{G}(x,y)\geq\inf_{(p_\pi)_{\pi\in\Pi}}\sum_{\pi\in\Pi}p_\pi^2 \sum_{e\in E} r_e \mathbf{1}_{\{e\in\pi\}}\geq d_G(x,y)\inf_{(p_\pi)_{\pi\in\Pi}}\sum_{\pi\in\Pi}p_\pi^2= \frac{d_{G}(x,y)}{\#\Pi},\]
where the second inequality is immediate from the definition of the shortest path metric, and the equality is obtained by solving an elementary constrained minimisation problem. That $\#\Pi\leq \#E!$ is a simple counting exercise, and this completes the proof.
\end{proof}

\subsection{Local times of Brownian motion on a circle}\label{ltsec}

In this section, we deduce the estimate for the local times of Brownian motion on a circle that is used in proving the tightness result of Lemma \ref{tightness}. By Brownian motion on a circle of perimeter $r>0$, we mean the process $X^r=(X^r_t)_{t\geq 0}$ obtained by setting $X^r_t=\phi_r(X_t)$, where $X=(X_t)_{t\geq 0}$ is a standard Brownian motion on $\mathbb{R}$ started from 0, which is assumed to be built on a probability space with probability measure $\mathbf{P}$, and $\phi_r$ is the canonical projection from $\mathbb{R}$ to the circle $\mathbb{T}_r:=\mathbb{R}/r\mathbb{Z}$. This process corresponds to the local regular Dirichlet form $(\tfrac{1}{2}\mathcal{E}_r,\mathcal{F}_r)$ on $L^2(\mathbb{T},\lambda_r)$, where $\lambda_r$ is the one-dimensional Hausdorff measure on $\mathbb{T}_r$, which is assumed to be equipped with the quotient metric $d_r$ corresponding to Euclidean distance on $\mathbb{R}$,
\[\mathcal{E}_r(f,f):=\int_{\mathbb{T}_r}f'(x)^2\lambda_r(dx)\]
and
\[\mathcal{F}_r:=\{f\in L^2(\mathbb{T}_r,\lambda_r):f\mbox{ is absolutely continuous},\mathcal{E}_r(f,f)<\infty\}.\]
We can check $(\mathcal{E}_r,\mathcal{F}_r)$ is also a resistance form on $\mathbb{T}_r$, whose resistance metric satisfies
\begin{equation}\label{resr}
R_r(x,y)=\left(d_r(x,y)^{-1}+(r-d_r(x,y))^{-1}\right)^{-1},\hspace{20pt}\forall x,y\in\mathbb{T}_r.
\end{equation}
As a consequence of this, exactly as for the processes of Sections \ref{sectree} and \ref{seccont}, we obtain the existence of jointly continuous local times $(L^r_t(x))_{t\geq0,x\in\mathbb{T}_r}$ for $X^r$, and it is straightforward to deduce that we can write
\begin{equation}\label{ltsum}
L^r_t(x)=\sum_{y\in\phi_r^{-1}(x)}L_t(y),
\end{equation}
where $(L_t(x))_{t\geq 0,x\in\mathbb{R}}$ are the jointly continuous local times for $X$ (cf. \cite{Bolt}). We use these local times to construct the trace of the process $X$ on an arc of length $\varepsilon_0\leq r/2$. In particular, identifying $\mathbb{T}_r$ with the interval $[0,r)$ in the natural way, define
\begin{equation}\label{ar}
A^{r,\varepsilon_0}_t:=\int_{[0,\varepsilon_0]}L^r_t(x)\lambda_r(dx),
\end{equation}
its inverse $\tau^{r,\varepsilon_0}(t):=\inf\{s:{A}^{r,\varepsilon_0}_s>t\}$, and set $X^{r,\varepsilon_0}_t:=X^r_{\tau^{r,\varepsilon_0}(t)}$, which is the process associated with $\tfrac{1}{2}{\rm Tr}(\mathcal{E}_r|[0,\varepsilon_0])$ considered as a Dirichlet form on $L^2([0,\varepsilon_0],\lambda_r([0,\varepsilon_0]\cap\cdot)$, and has local times $L^{r,\varepsilon_0}_t(x):=L^r_{\tau^{r,\varepsilon_0}(t)}(x)$. The particular tightness estimate we require for these local times is the following.

{\lem \label{ltlem}Fix $\varepsilon_0>0$. For every $\varepsilon>0$ and $t_0<\infty$,
\begin{equation}\label{lttight}
\lim_{\delta\rightarrow 0}\sup_{r\geq2\varepsilon_0}\mathbf{P}\left(\sup_{\substack{s,t\in[0,t_0]:\\|s-t|\leq \delta}}\sup_{\substack{x,y\in[0,\varepsilon_0]:\\|x-y|\leq \delta}}\left|L^{r,\varepsilon_0}_s(x)-L^{r,\varepsilon_0}_t(y)\right|>\varepsilon\right)=0.
\end{equation}}
\begin{proof}
Note that the probability in (\ref{lttight}) can be bounded above by
\begin{eqnarray}
\lefteqn{\mathbf{P}\left(\sup_{\substack{t\in[0,t_0]}}\sup_{\substack{x,y\in[0,\varepsilon_0]:\\|x-y|\leq \delta}}\left|L^{r,\varepsilon_0}_t(x)-L^{r,\varepsilon_0}_t(y)\right|>\varepsilon/2\right)}\nonumber\\
&&+\mathbf{P}\left(\sup_{\substack{s,t\in[0,t_0]:\\|s-t|\leq \delta}}\sup_{x\in[0,\varepsilon_0]}\left|L^{r,\varepsilon_0}_s(x)-L^{r,\varepsilon_0}_t(x)\right|>\varepsilon/2\right).\label{twoterms}
\end{eqnarray}
We will consider each of these terms separately. To deal with the first of these terms, we first observe that, by (\ref{ltsum}) and (\ref{ar}), $\mathbf{P}$-a.s.,
\[\inf_{r\geq2\varepsilon_0}A_t^{r,\varepsilon_0}=\inf_{r\geq2\varepsilon_0}\int_{[0,\varepsilon_0]}L^r_t(x)\lambda_r(dx)
\geq \int_{[0,\varepsilon_0]}L_t(x)dx\rightarrow \infty,\]
as $t\rightarrow\infty$, where the limit is an application of the strong Markov property for Brownian motion (see \cite{Kallenberg}, Chapter 22, Exercise 9, for example). Consequently, it is the case that $t_1:=\sup_{r\geq2\varepsilon_0}\tau^{r,\varepsilon_0}(t_0)$ is finite, $\mathbf{P}$-a.s., and so, for $\delta\leq \varepsilon_0$,
\begin{eqnarray}
\lefteqn{\sup_{r\geq2\varepsilon_0}\sup_{\substack{t\in[0,t_0]}}\sup_{\substack{x,y\in[0,\varepsilon_0]:\\|x-y|\leq \delta}}\left|L^{r,\varepsilon_0}_t(x)-L^{r,\varepsilon_0}_t(y)\right|}\nonumber\\
&\leq& \sup_{r\geq2\varepsilon_0}\sup_{\substack{t\in[0,t_1]}}\sup_{\substack{x,y\in[0,\varepsilon_0]:\\|x-y|\leq \delta}}\left|L^{r}_t(x)-L^{r}_t(y)\right|\nonumber\\
&\leq&\sup_{r\geq 2\varepsilon_0}\sup_{\substack{t\in[0,t_1]}}\sup_{\substack{x,y\in[0,\varepsilon_0]:\\|x-y|\leq \delta}}\sum_{n\in\mathbb{Z}}\left|L_t(x+nr)-L_t(y+nr)\right|\nonumber\\
&\leq&\sup_{\substack{t\in[0,t_1]}}\sup_{\substack{x,y\in\mathbb{R}:\\|x-y|\leq \delta}}\left|L_t(x)-L_t(y)\right|\sup_{r\geq 2\varepsilon_0}2(1+|X_t|/r)\label{terms}\\
&\rightarrow&0\nonumber
\end{eqnarray}
where the final supremum in (\ref{terms}) is an upper bound for the number of non-zero terms in sum in the previous line, and the limit holds $\mathbf{P}$-a.s. as a result of the joint continuity of the local times of standard Brownian motion. That the first term of (\ref{twoterms}) decays uniformly over $r\geq 2\varepsilon_0$ as $\delta\rightarrow 0$ readily follows from this.

For the second term, we will compare the local times of the process on a circle with those on a line segment. Applying the identification of $\mathbb{T}_r$ with $[0,r)$, define a map $\psi_r:\mathbb{T}_r=[0,r)\rightarrow [0,r/2]$ by setting
\[\psi_r(x)=\left\{\begin{array}{ll}
                   \varepsilon_0/2-x, & \mbox{if }x\in[0,\varepsilon_0/2], \\
                   x-\varepsilon_0/2, & \mbox{if }x\in[\varepsilon_0/2,\varepsilon_0/2+r/2], \\
                   r+\varepsilon_0/2-x, & \mbox{otherwise.}
                 \end{array}
\right.\]
We then have that $\tilde{X}^r:=\psi_r(X^r)$ is reflected Brownian motion on $[0,r/2]$ started from $\varepsilon_0/2$, and has jointly continuous local times $\tilde{L}^r_t(x)=\sum_{y\in\psi_r^{-1}(x)}L_t^r(y)$. If we set $\tilde{X}^{r,\varepsilon_0}_t:=\tilde{X}^r_{\tau^{r,\varepsilon_0}(t)}$, then, since
\[A^{r,\varepsilon_0}_t=\int_{[0,\varepsilon_0]}L^r_t(x)\lambda_r(dx)=\int_{[0,\varepsilon_0/2]}\tilde{L}^r_t(x)dx,\]
the process $\tilde{X}^{r,\varepsilon_0}$ is the trace of $\tilde{X}^r$ on $[0,\varepsilon_0/2]$, which is simply reflected Brownian motion on $[0,\varepsilon_0/2]$ started from $\varepsilon_0/2$, regardless of $r$. Moreover, the local times of $\tilde{X}^{r,\varepsilon_0}$ are given by $\tilde{L}^{r,\varepsilon_0}_t(x)=\tilde{L}^{r}_{\tau^{r,\varepsilon_0}(t)}(x)=\sum_{y\in\psi_r^{-1}(x)}{L}^{r}_{\tau^{r,\varepsilon_0}(t)}(y)=\sum_{y\in\psi_r^{-1}(x)}{L}^{r,\varepsilon_0}_{t}(y)$. Consequently,
\[{\sup_{\substack{s,t\in[0,t_0]:\\|s-t|\leq \delta}}\sup_{x\in[0,\varepsilon_0]}\left|L^{r,\varepsilon_0}_s(x)-L^{r,\varepsilon_0}_t(x)\right|}\leq\sup_{\substack{s,t\in[0,t_0]:\\|s-t|\leq \delta}}\sup_{x\in[0,\varepsilon_0/2]}\left|\tilde{L}^{r,\varepsilon_0}_s(x)-\tilde{L}^{r,\varepsilon_0}_t(x)\right|,\]
which implies that  the second term of (\ref{twoterms}) is bounded above by
\[\mathbf{P}\left(\sup_{\substack{s,t\in[0,t_0]:\\|s-t|\leq \delta}}\sup_{x\in[0,\varepsilon_0/2]}\left|\tilde{L}^{2\varepsilon_0,\varepsilon_0}_s(x)-\tilde{L}^{2\varepsilon_0,\varepsilon_0}_t(x)\right|>\varepsilon\right),\]
uniformly over $r\geq 2\varepsilon_0$. Since the local times $\tilde{L}^{2\varepsilon_0,\varepsilon_0}$ are jointly continuous, this probability converges to zero as $\delta\rightarrow 0$.
\end{proof}

\subsection{Random walk estimates}\label{rwe}

In proving the tightness of the rescaled local times of $X^{\mathcal{C}_1^n}$ in Lemma \ref{ltight}, we apply a pair of estimates for simple random walks on graphs that are proved in this section. The first is a tail bound for the occupation time of a simple random walk on an interval satisfying a certain boundary condition. The second involves the moments of the number of visits a random walk makes to a particular vertex before returning to its starting point. We note that both results are adaptations of estimates that appear in \cite{Croydoncbp}, Appendix B.
\medskip

\subsubsection{Occupation time for simple random walk on an interval}

Let $L$ be an integer and $p\in(0,1]$. Suppose $X^{(L,p)}=(X^{(L,p)}_n)_{n\geq 0}$ is a Markov chain on $\{0,1,\dots,2L\}$ such that $X^{(L,p)}$ behaves like a symmetric simple random walk on the vertices $\{1,2,\dots,2L-1\}$ (i.e. jumps up one or down one with probability $\tfrac{1}{2}$). We do not specify the transition probabilities for $X^{(L,p)}_n\in\{0,2L\}$, apart from assuming a jump from $0$ to $1$ occurs with probability of $p$. Write $\mathbf{P}_x^{(L,p)}$ for the law of $X^{(L,p)}$ started from $x$. For $\xi^{(L,p)}_n$, the number of visits to 0 up to an including time $n$, the following tail bound holds.

{\lem \label{ltexp} There exists a constant $c>0$, not depending on $(L,p)$, such that
\[\mathbf{P}_0^{(L,p)}\left(\xi_{sL^2}^{(L,p)}\geq tL\right)\leq 2^{2+s}e^{-cpt},\]
for every $s,t>0$.}
\begin{proof} Let $\sigma$ be the time of the first return of $X^{(L,p)}$ to 0. Applying the Markov property at the first step implies that
\begin{eqnarray*}
\lefteqn{\mathbf{P}_0^{(L,p)}\left(\sigma\geq L^2\right)}\\
&\geq& p \mathbf{P}_1^{(L,p)}\left(\mbox{$X^{(L,p)}$ hits $L$ before 0}\right)\\
&&\times \mathbf{P}_L^{(L,p)}\left(\mbox{$X^{(L,p)}$ spends $\geq L^2$ in $\{1,2,\dots,2L-1\}$ before hitting $\{0,2L\}$}\right).
\end{eqnarray*}
Using the standard gambler's ruin estimate, the second term is equal to $L^{-1}$. As $L\rightarrow \infty$, the third term converges to $\mathbf{P}(\tau_{\{-1,1\}}\geq 1)>0$, where $\tau_{\{-1,1\}}$ is the hitting time of $\{-1,1\}$ by a standard Brownian motion on $\mathbb{R}$ started from $0$. Thus there exists a constant $c>0$, not depending on $(L,p)$ such that
\[\mathbf{P}_0^{(L,p)}\left(\sigma\geq L^2\right)\geq \frac{cp}{L}.\]
Consequently if $(\sigma_i)_{i\geq 1}$ are independent identically-distributed copies of $\sigma$, then
\begin{eqnarray*}
\mathbf{P}_0^{(L,p)}\left(\xi^{(L,p)}_{sL^2}\geq tL\right)&\leq &\mathbf{P}\left(\sum_{i=1}^{\lfloor tL\rfloor -1}\sigma_i\leq sL^2\right)\\
&\leq &\mathbf{P}\left(\sum_{i=1}^{\lfloor tL\rfloor -1}\mathbf{1}_{\{\sigma_i\geq L^2\}}\leq s\right)\\
&\leq &\mathbf{P}\left({\rm Bin}(\lfloor tL\rfloor -1,cpL^{-1})\leq s\right)\\
&\leq& \mathbf{E}\left(2^{s- {\rm Bin}(\lfloor tL\rfloor -1,cpL^{-1})}\right)\\
&\leq&2^s\left(1-\frac{cp}{2L}\right)^{tL-2},
\end{eqnarray*}
where ${\rm Bin}(n,p)$ represents a binomial random variable with parameters $n$ and $p$. The result follows.
\end{proof}

\subsubsection{Number of visits to vertices before return}

We now consider an arbitrary locally finite graph $G=(V,E)$. Denote the usual graph degree of a vertex $x\in V$ by ${\rm deg}_{G}(x)$, and the effective resistance metric by $R_G$ (where edges are assumed to have unit conductance). A standard result for the random walk $X^G=(X^G_n)_{n\geq 0}$ on $G$ is that
\[\mathbf{P}_x^G\left(\tau_y^G<\tau_x^G\right)=\frac{1}{{\rm deg}_{G}(x)R_G(x,y)},\]
for every $x\neq y\in V$, where $\mathbf{P}_x^G$ is the law of $X^G$ started from $x$, and $\tau_z^G$ is the first strictly positive time the random walk $X^G$ hits the vertex $z\in V$ (see \cite{LP}, Section 2.2, for example). Hence if $N_G(x,y)$ is the number of visits to $y$ before returning to $x$, then $\mathbf{P}_x^G\left(N_G(x,y)=k\right)$ is equal to
\begin{equation}\label{probk}
\frac{1}{{\rm deg}_{G}(x)R_G(x,y)}\left(1-\frac{1}{{\rm deg}_{G}(y)R_G(x,y)}\right)^{k-1}\frac{1}{{\rm deg}_{G}(y)R_G(x,y)},
\end{equation}
for $k\geq 1$. From this, it is possible to check that $\mathbf{E}^G_xN_G(x,y)={\rm deg}_{G}(y){\rm deg}_{G}(x)^{-1}$. Now, introduce a centred random variable \[\eta_G(x,y):=\frac{N_G(x,y)}{{\rm deg}_{G}(y)}-\frac{1}{{\rm deg}_{G}(x)}.\]
From (\ref{probk}), applying standard results about the moments of a geometric random variable, it is an elementary exercise to deduce the following lemma.

{\lem \label{etamoments} There exists a universal finite constant $c$ such that
\[\left|\mathbf{E}_x^G\left(\eta_G(x,y)^k\right)\right|\leq c d_G(x,y)^{k-1},\]
for every $x\neq y\in V$ and $k\in\{1,2,3,4\}$.}

\subsection{Brownian motion on a real tree star}\label{star}

Suppose $T=(T,d_T)$ is a real tree star formed by including unit line segments along the edges of a star graph with $D+1$, $D\geq 1$, vertices. The root $\rho$ of the real tree will be the internal (degree $D$) vertex of the star graph and we write $L\subseteq T$ to represent the set of external (degree 1) vertices (in the case when $D=1$, $\rho$ is simply assumed to be one end of the line segment and $L$ consists of the other end vertex). Let $(\mathcal{E}_T,\mathcal{F}_T)$ be the resistance form associated with $(T,d_T)$ through (\ref{dtchar}). Write $X^T$ to represent the Brownian motion on $T$; that is, the process associated with $(\tfrac{1}{2}\mathcal{E}_T,\mathcal{F}_T)$ when this quadratic form is considered as a Dirichlet form on $L^2(T,\lambda_T)$, where $\lambda_T$ is the one-dimensional Hausdorff measure on $T$. By \cite{Croydoninf}, Lemma 2.2, $X^T$ admits jointly continuous local times $(L^T_t(x))_{x\in T,t\geq0}$.

In determining how well the simple random walk on $\mathcal{C}_1^n(k)$ is approximated by the Brownian motion on the corresponding metric space with line segments included (see proof of Lemma \ref{kconv}) we apply the following lemma, which gives a pair of simple properties of $X^T$ and $L^T$ run up to the hitting time of $L$, which we denote by $\tau_L$.

{\lem\label{starlabel} Conditional on $X^T_0=\rho$, $X_{\tau_L}$ is distributed uniformly on $L$,
\begin{equation}\label{XT}
\left(d_T\left(\rho,X^T_{t\wedge \tau_L}\right)\right)_{t\geq 0}\buildrel{d}\over {=} \left(|B_{t\wedge \tau_{\{\pm1\}}}|\right)_{t\geq 0},
\end{equation}
\begin{equation}\label{LT}
DL^T_{\tau_L}(\rho)\buildrel{d}\over {=}2L^B_{\tau_{\{\pm1\}}}(0),
\end{equation}
where $B$ is a standard Brownian motion on $\mathbb{R}$ started from $0$, $L^B$ is a jointly continuous version of its local times, and $\tau_{\{\pm1\}}$ is the hitting time of $\{\pm1\}$ by $B$.}
\begin{proof} Let $T_n=\{x\in T:nd_T(\rho,x)\in\mathbb{Z}\}$ and $\lambda_n$ be the measure on $T_n$ defined by supposing that $2n\lambda_n(\{x\})$ is equal to the number of connected components of $T\backslash\{x\}$. Let $\mathcal{E}_{T_n}={\rm Tr}(\mathcal{E}_T|T_n)$ and $X^{T_n}$ be the process associated with $(\tfrac{1}{2}\mathcal{E}_{T_n},\mathcal{F}_{T_n})$ considered as a Dirichlet form on $L^2(T_n,\lambda_n)$. By the trace theorem for Dirichlet forms (\cite{FOT}, Theorem 6.2.1), it is possible to assume that $X^{T_n}$ and $X^T$ are coupled through the relation $X^{T_n}_t=X^T_{\tau^n(t)}$, where $(\tau^n(t))_{t\geq 0}$ is the right-continuous inverse of
\[A^n_t:=\int_{T_n}L^T_t(x)\lambda_n(dx).\]
Clearly $\lambda_n\rightarrow \lambda_T$ weakly as measures on $T$, and therefore $A_t^n\rightarrow t$ uniformly on compacts, almost-surely. From this and the almost-sure continuity of $X^T$, it follows that $X^{T_n}$ converges in distribution to $X^T$ as $n\rightarrow\infty$.

By construction, we also have that $X^{T_n}$ is the continuous time simple random walk on $T_n$ (equipped with the graph structure induced by the line segments of $T$), where the holding time at each vertex is an exponential, mean $n^{-2}$, random variable (cf. \cite{Barlow}, Remark 7.23). Thus, we immediately obtain that $(d_T(\rho,X^{T_n}_t))_{t\geq 0}$ has the same distribution as $(n^{-1}S_{tn^2})_{t\geq 0}$, where $S$ is a continuous time simple random walk on $\{0,1,\dots, n\}$ started from 0 with exponential, mean 1, holding times. Hence $(d_T(\rho,X^{T_n}_t))_{t\geq 0}$ converges to a reflected Brownian motion on $[0,1]$, started from 0. Since, by the continuous mapping theorem, $(d_T(\rho,X^{T_n}_t))_{t\geq 0}$ also converges to $(d_T(\rho,X^{T}_t))_{t\geq 0}$, the distributional equality at (\ref{XT}) follows. The claim that $X^T_{\tau_L}$ is uniformly distributed on $L$ is another easy consequence of this discrete approximation picture (or can be proved directly from the resistance form construction of $X^T$).

To prove (\ref{LT}), we start by supposing that $X^T$ and $B$ are coupled so that the relation at (\ref{XT}) holds almost-surely. Let $f:[0,1]\rightarrow \mathbb{R}$ be a continuous function. One can easily check that
\[\int_{-1}^1f(|x|)L^B_{\tau_{\{\pm 1\}}}(x)dx=\int_Tf(d_T(\rho,x))L^T_{\tau_L}(x)\lambda(dx).\]
Since $f$ is arbitrary, it follows that, for Lebesgue-almost-every $r\in(0,1)$,
\[L^B_{\tau_{\{\pm 1\}}}(r)+L^B_{\tau_{\{\pm 1\}}}(-r)=\sum_{x\in T:d_T(\rho,x)=r}L^T_{\tau_L}(x).\]
Letting $r\rightarrow 0$ and applying the continuity of $L^T$ and $L^B$, we obtain (\ref{LT}).
\end{proof}
\end{appendices}

\noindent
{\bf Acknowledgements.} The author would like to thank Takashi Kumagai for indicating that this might be an interesting problem to work on, and also for suggesting the extension of the random walk convergence result discussed at the end of Section \ref{rwc}.

\def\cprime{$'$}
\providecommand{\bysame}{\leavevmode\hbox to3em{\hrulefill}\thinspace}
\providecommand{\MR}{\relax\ifhmode\unskip\space\fi MR }
\providecommand{\MRhref}[2]{%
  \href{http://www.ams.org/mathscinet-getitem?mr=#1}{#2}
}
\providecommand{\href}[2]{#2}

\end{document}